\documentclass[a4paper,11pt]{amsart} % ,9pt
\usepackage{amsmath}
\usepackage{fullpage}
\usepackage{amssymb}%This gives us symbols like \square
\usepackage{verbatim}%This gives us the \begin{comment} environment.
\usepackage{amsthm}%This gives us Definition style of theorem.
\usepackage{mathrsfs}%This gives us \mathscr#
\usepackage{graphicx}
\usepackage{mathtools}
\usepackage[colorlinks,pdfdisplaydoctitle,linkcolor=blue,citecolor=blue]{hyperref}
\usepackage{caption}
\usepackage{subcaption}
\usepackage{xcolor}
\usepackage{url}
\usepackage{epstopdf}
\usepackage{pgf,tikz}
\usepackage{bbm}
\usepackage{extarrows}
\usetikzlibrary{arrows}
\usepackage{xcolor}
\usepackage{float}
\usepackage{placeins}
\usepackage{amssymb,amsthm}    % amssymb package contains more mathematical symbols
\usepackage{graphicx}          % graphicx package enables you to paste in graphics
\usepackage{amsmath}         % best maths package on offer
\usepackage{caption}
\usepackage{subcaption}
\usepackage{dsfont}
\usepackage[backend=biber,
  style=alphabetic,
  giveninits=true,
  maxnames=3,      % controls citations
  minnames=1,
  maxbibnames=99   % show all names in bibliography
]{biblatex}
\usepackage{listings}
\usepackage{bbm}
\usepackage{color}
\usepackage{enumerate}
\usepackage{todonotes}
\usepackage[ruled,vlined,noline]{algorithm2e}
\usepackage{algorithmic}
\usepackage{parskip}
\usepackage{soul}

\let\bm\boldsymbol

\renewcommand{\Delta}{\triangle}

\definecolor{darkblue}{rgb}{0,0,0.7}

\definecolor{darkgreen}{rgb}{0.01,0.75,0.24}

\def \Ee[#1]{\mathcal{E}^{\text{{#1}}}}

\def\pa[#1,#2]{\frac{\partial {#1}}{\partial {#2}} }

\def\idom[#1,#2,#3]{\int_{#1}\hspace{1pt} {#2} \hspace{1pt} \text{d}{#3}}
\def\res[#1,#2]{\left.{#1}\right|_{#2}}

\def\var[#1,#2]{\langle \delta \mathcal{E}^{\text{{#1}}}({#2}),v\rangle}
\def\vars[#1,#2,#3]{\langle \delta^2\mathcal{E}^{\text{{#1}}}({#2})v,{#3}\rangle}
\def\vard[#1,#2,#3,#4]{\langle \delta\mathcal{E}^{\text{{#1}}}({#2})-\delta\mathcal{E}^{\text{{#3}}}({#4}),v\rangle}

\def\E{\E}

\newcommand{\be}{\begin{equation}}
\newcommand{\en}{\end{equation}}
\newcommand{\ben}{\begin{equation*}}
\newcommand{\enn}{\end{equation*}}
\newcommand{\bea}{\begin{aligned}}
\newcommand{\ena}{\end{aligned}}

\def\ba#1\ena{\begin{align}#1\end{align}}
\def\ban#1\enan{\begin{align*}#1\end{align*}}

\theoremstyle{plain}
\newtheorem{theorem}{Theorem}[section]

\newtheorem{remark}[theorem]{Remark}

\numberwithin{equation}{section}

\theoremstyle{definition}

\renewcommand{\E}{\mathbb{E}}

\begin{document}

\title[Cluster formation for weakly interacting kinetic Langevin dynamics]{\textbf{Cluster formation for weakly interacting kinetic Langevin dynamics}}

\author[B. Leimkuhler] {Benedict Leimkuhler}
\address{Department of Mathematics, University of Edinburgh, EH9 3FD, Edinburgh, United Kingdom}
\email{b.leimkuhler@ed.ac.uk}

\author[R. Lohmann] {Ren\'e Lohmann}
\address{Department of Mathematics, University of Edinburgh, EH9 3FD, Edinburgh, United Kingdom}
\email{r.lohmann@ed.ac.uk}

\author[G. A. Pavliotis] {Grigorios A. Pavliotis}
\address{Department of Mathematics, Imperial College London, SW7 2AZ, London, United Kingdom}
\email{g.pavliotis@imperial.ac.uk}

\author[P. A. Whalley] {Peter A. Whalley}
\address{Seminar for Statistics, ETH Z{\"u}rich, Z{\"u}rich, Switzerland}
\email{pwhalley@ethz.ch}

%\subjclass{94A12, 86A22, 60G35, 62M99.}
%\keywords{steepest descent, stochastic coefficients, Banach spaces, comparison of experiments.}

\maketitle
%\tableofcontents
%\begin{appendix}

\begin{abstract}
In this paper, we study the formation of clusters for stochastic interacting particle systems (IPS) that interact through short-range attractive potentials in a periodic domain. We consider kinetic (underdamped) Langevin dynamics and focus on the low-friction regime. Employing a linear stability analysis for the kinetic McKean-Vlasov equation, we show that, at sufficiently low temperatures, and for sufficiently short-ranged interactions, the particles form clusters that correspond to metastable states of the mean-field dynamics. We derive the friction and particle-count dependent cluster-formation time and numerically measure the friction-dependent times to reach a stationary state (given by a state in which all particles are bound in a single cluster). By providing both theory and numerical methods in the inertial stochastic setting, this work acts as a bridge between cluster formation studies in overdamped Langevin dynamics and the Hamiltonian (microcanonical) limit.
\end{abstract}

\section{Introduction}

Interacting particle systems arise in many applications, ranging from plasma physics \cite{bittencourt2013fundamentals} and stellar dynamics \cite{binney2008galactic} to biology \cite{Suz05} and are also widely used in algorithms for sampling and optimization \cite{carrillo2018analytical,mei2018mean}, and mathematical models in the social sciences \cite{rainer2002opinion,GSW19,MKT18,MaF19,EGP18, GGS21}. We refer to~\cite{Diez_2022a, Diez_2022b} for a recent review and references to the literature. Quite often, interactions between particles {(or `agents')} on the microscale lead to the emergence of collective behavior at the macroscale \cite{NPT10, garnier2017consensus}. This collective behavior manifests itself in various ways: the formation of clusters in stellar dynamics models~\cite{binney2008galactic}, the emergence of consensus in models for opinion formation \cite{rainer2002opinion,garnier2017consensus,GGS21, WSMSPW_2025}, synchronization in systems of interacting nonlinear oscillators~\cite{bertini2010dynamical}, the emergence of multipeak states in models in biology and in dynamics on graphons~\cite{PSV_2009, BBBGGP_2024, bertoli2025phasetransitionsinteractingparticle}, dynamical clustering in active matter~\cite{Caprini2024DynamicalClustering}, pattern formation~\cite{delfau2016pattern},  or the swarming of animal populations \cite{leverentz2009asymptotic}. Recently, mathematical models for the dynamics of transformers have also been shown to exhibit clustering phenomena~\cite{GLPR_2205, shalova2024solutions, bruno2025multiscaleanalysismeanfieldtransformers, bruno2024emergence,balasubramanian2025structure}. In these works, the evolution of tokens, as they travel through transformer layers, is interpreted as a system of weakly interacting particles on the sphere, with a short-ranged attractive interaction potential that is similar to a multidimensional extension of the Hegselmann-Krause model for opinion dynamics in one dimension~\cite{garnier2017consensus}. In these papers it was observed (and under certain assumptions rigorously proved) that tokens tend to organize into clusters. This observation provides a compelling qualitative explanation of how transformer models can develop representations of complex input data.

Many different physical and biological systems exhibit cluster formation. Examples range from the formation of clusters of galaxies to molecular clusters to active matter. The emergence of clusters in multi-agent systems and, in particular, in noise-driven interacting particle systems and their mean-field limit, is often a manifestation of dynamical metastability~\cite{WSMSPW_2025, posch1990dynamics, Butta_al_2003}; clusters, or, more generally, multipeak states, can be dynamically stable over long time intervals, before the dynamics converges to a stationary state~\cite{PSV_2009, bertoli2025phasetransitionsinteractingparticle}.  From the perspective of statistical physics, metastability is the dynamical manifestation of a first-order phase transition. Indeed, the emergence of collective behavior, e.g., synchronization, consensus formation, flocking and swarming in noise-driven multi-agent systems can often be interpreted as a phase transition between a disordered and ordered state~\cite{Carrillo2020}.

In this paper, we study, both numerically and analytically, the kinetics of interacting particle systems that
exhibit collective behavior and, in particular, the formation of clusters. Our model framework is weakly interacting underdamped Langevin dynamics:
\begin{equation}\label{eq:Langevin}
\begin{split}
\Dot{\bm{x}}_{i} &= \bm{v}_{i},\\
\Dot{\bm{v}}_{i} &= -\nabla_{{\boldsymbol{x}_i}} \Phi_{N}(\bm{x}) - \gamma \bm{v}_{i} + \sqrt{2\gamma \beta^{-1}}\dot{\mathbf{B}}_{i},
\end{split}
\end{equation}
where $\bm{x} = (\bm{x}_{1},...,\bm{x}_{N})$ denotes the position vector of $N$ {unit mass} particles, $\bm{v}=(\bm{v}_{1},...,\bm{v}_{N})$ the corresponding momentum vector, $\gamma$ the friction coefficient, $\beta$ the inverse temperature and $\mathbf{B} = (\mathbf{B}_1, . . . ,\mathbf{B}_N)$ a collection of $N$ independent
$d$-dimensional standard Brownian motions. The potential $\Phi_{N}$ typically consists of a confining potential and pairwise interactions, 
\[
\Phi_{N}(\bm{x}) = \sum^{N}_{j = 1}V(\bm{x}_{j}) + \frac{1}{2 N}\sum^{N}_{i,j=1} W(\bm{x}_{i}-\bm{x}_{j}),
\]
where the pairwise interaction potential $W$ only depends on the particle distances, i.e., $W(\bm{x}_{i}-\bm{x}_{j})=\widetilde{W}(\|\bm{x}_{i}-\bm{x}_{j}\|)$. 
We will consider the dynamics \eqref{eq:Langevin} in $\mathcal{D}_{\bm{x}} \times {\mathcal{D}_{\bm{v}}}$, where the {velocity domain is $\mathcal{D}_{\bm{v}}:=\mathbb{R}^{d}$}, and the configurational domain consists of a box with periodic boundary conditions, i.e., the $d$-dimensional torus of side length $L$, $\mathcal{D}_{\bm{x}}:=\mathbb{T}^d_L=[0,L)^d$. {Since the configurational domain is bounded, we do not require a confining potential and set $V\equiv 0$.} The purpose of this paper is to understand the formation of clustered metastable states for the dynamics~\eqref{eq:Langevin} and for attractive interaction potentials. We study the problem by means of numerical simulations and by a stability analysis of the mean-field $N \to +\infty$. {While our analysis focuses on the one-dimensional case, $d=1$, we believe the results are generalizable to higher dimensions, and we also consider the case $d=2$ in our simulations.} The interacting particle system~\eqref{eq:Langevin} has been extensively studied, both in Hamiltonian and stochastic settings.  In the Hamiltonian case one is led to consider the Vlasov mean-field PDE~\cite{posch1990dynamics, bavaud1991equilibrium, Spohn91}.  For treatment of the model with noise and dissipation ($\gamma >0$ in \eqref{eq:Langevin}), see, e.g.,~\cite{dressler1987, Dolbeault_1999, guillin2021uniform} and the references therein. A brief summary of the relevant results on the kinetic Langevin dynamics~\eqref{eq:Langevin} may be found in Section~\ref{sec:LD_theory_intro}. 

\begin{figure}[H]
\includegraphics[width=1.0\textwidth]{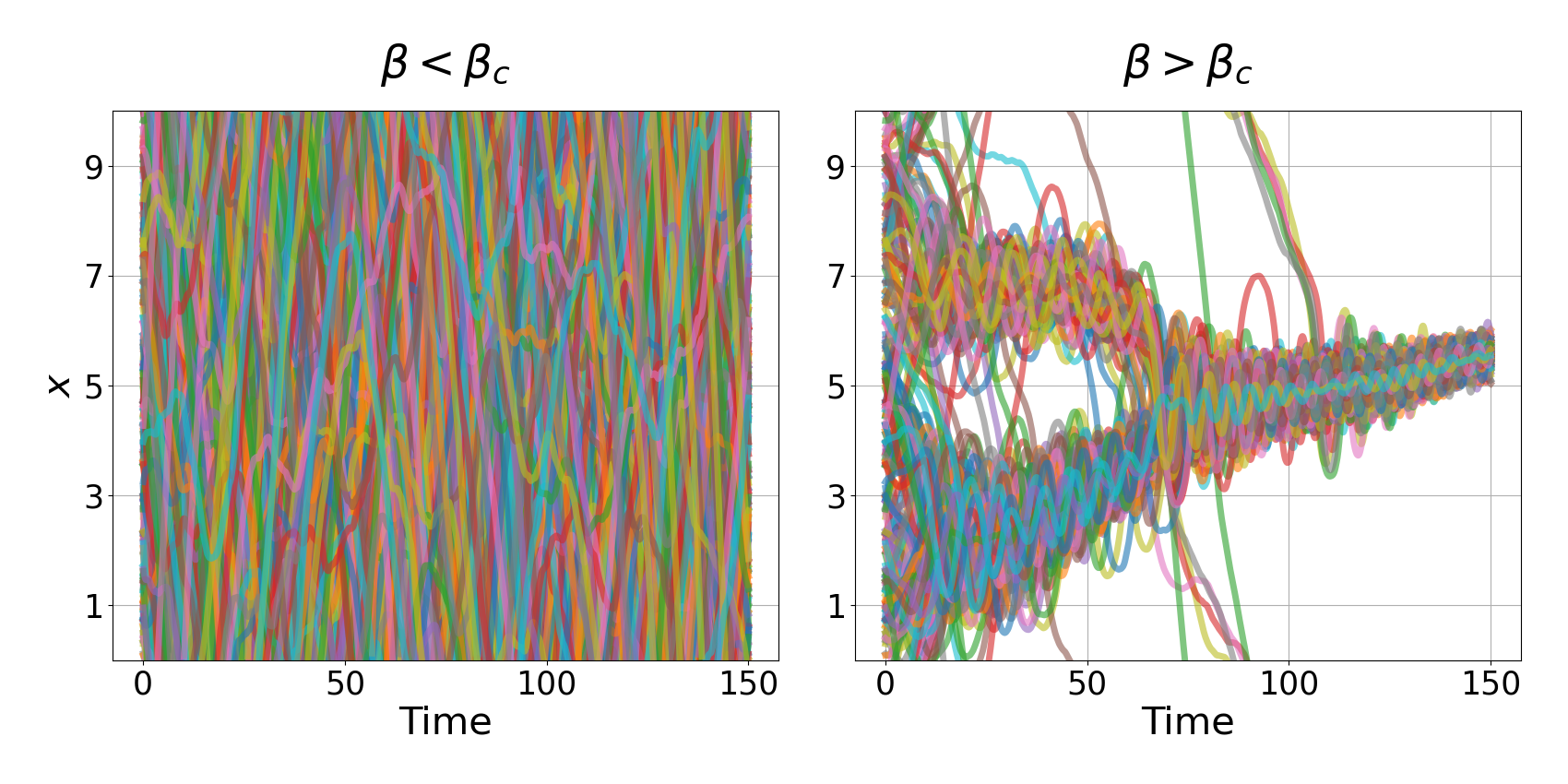}
\caption{\label{fig:example_trajectories1D} Example trajectories for the one-dimensional interacting particle system starting from a uniform distribution. \textbf{Left:} for reciprocal temperature $\beta<\beta_c$, {with $\beta_c$ the critical value,} the uniform distribution remains stable. \textbf{Right:} $\beta>\beta_c$ (i.e., small enough temperatures) the particle trajectories form two clusters that later merge into one cluster, yielding the new stationary state.}
\end{figure}

\begin{figure}[H]
\includegraphics[width=1.0\textwidth]{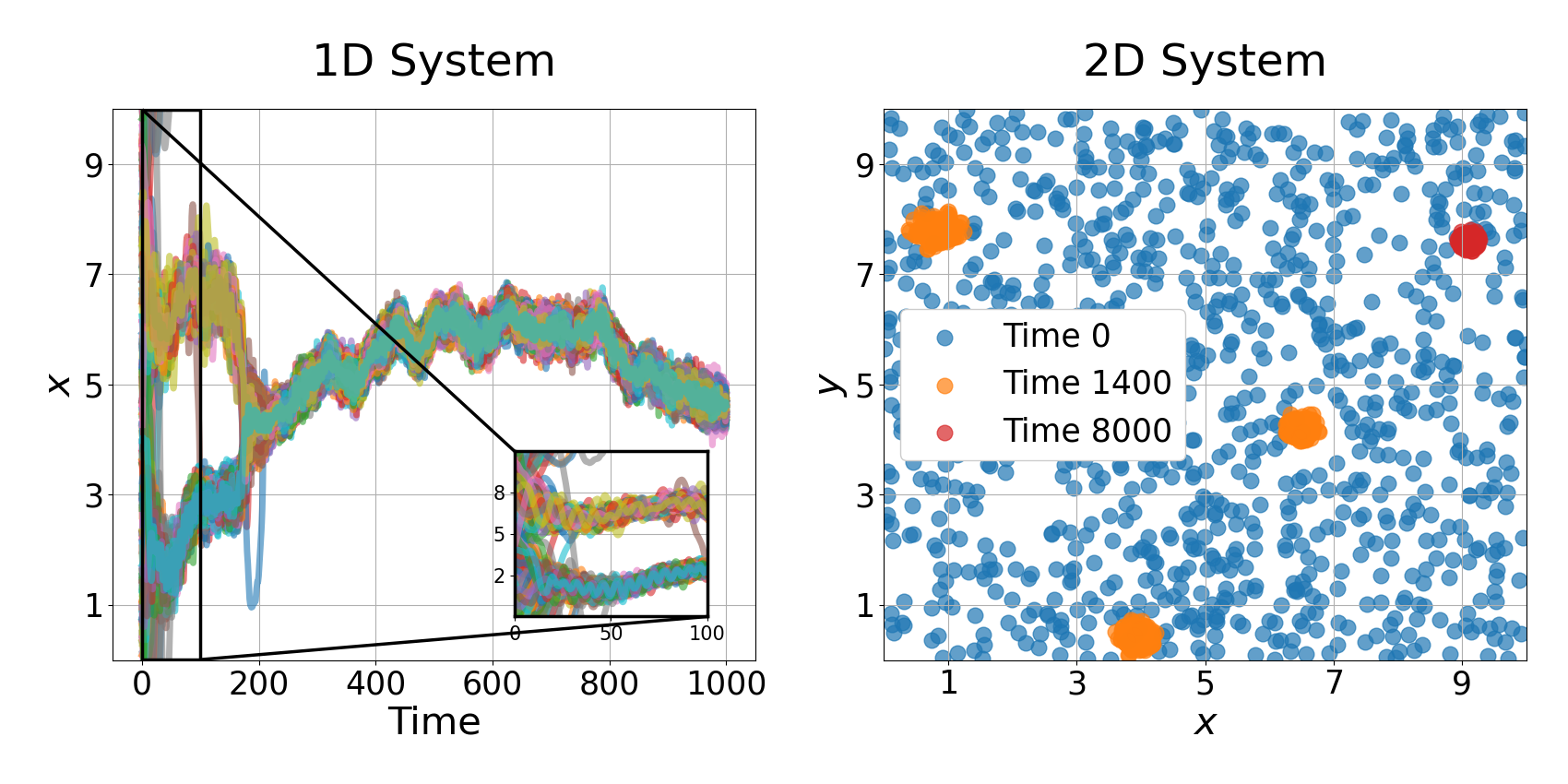}
\caption{\label{fig:example_trajectories1D_2D} Visualization of the IPS evolution in time in one dimension (\textbf{left}) and two dimensions (\textbf{right}) at $\beta > \beta_c$, leading to cluster formation. The clusters merge until only one cluster is left.}
\end{figure}

Cluster formation for Hamiltonian dynamics, $\gamma =0$ in~\eqref{eq:Langevin}, and for Gaussian attractive interactions, was studied in~\cite{posch1990dynamics}. The emergence of clustered states for the overdamped dynamics, obtained in the large $\gamma$ limit~\cite[Sec. 6.5]{PavliotisStochBook} and for Gaussian attractive interactions was investigated numerically in~\cite{martzel2001mean}. In particular, it was observed that the clusters coalesce and merge until a single cluster has formed. Similar numerical experiments, together with linear stability analysis of the uniform state, were performed for the Hagselmann-Krause bounded confidence model in~\cite{garnier2017consensus}. A study of cluster formation for overdamped dynamics was performed using the Dean-Kawasaki stochastic PDE in~\cite{WSMSPW_2025}. Naturally, the formation of clusters occurs only at sufficiently low temperatures. At high temperatures, diffusion due to Brownian motion wins over the attractive interaction between particles, and the uniform spatial state (corresponding to the disordered state) is stable. This is illustrated in Figures \ref{fig:example_trajectories1D} and \ref{fig:example_trajectories1D_2D} in one and two spatial dimensions.

The goal of this paper is to analyze cluster formation for the underdamped weakly interacting Langevin dynamics~\eqref{eq:Langevin} {across different attractive interaction potentials, exploring a wide range of friction coefficients (from the overdamped to the zero-friction limit) and establishing the critical temperature at which this collective behaviour emerges. }{Our investigations extend those presented in~\cite{martzel2001mean, garnier2017consensus, WSMSPW_2025} to the kinetic case. The presence of inertia leads to more interesting dynamical phenomena as well as additional analytical and numerical difficulties. The linear stability analysis is more involved due to the fact that the distribution function, the solution of the stationary McKean-Vlasov PDE, now depends on both position and velocity; see~\cite{Garnier_2019} for a similar linear stability analysis. Finally, the accurate numerical simulation of the interacting particle system in the underdamped regime requires the use of splitting numerical schemes for the solution of the corresponding stochastic differential equations~\cite{benMDbook}. }

The main contributions of this paper can be summarized as follows.

\begin{itemize}
\item We study the emergence of clusters~\eqref{eq:Langevin} for several attractive interaction potentials in the physically relevant weak-friction limit. Our analysis enables us, among other things, to derive an expression for the (inverse) critical temperature  $\beta_c$ below which clustering occurs. Furthermore, we show that the statistics of the clusters depend crucially on the friction coefficient, and we explore the limit as $\gamma$ goes to zero.

\item  We analyze the clustering time, i.e., the time it takes for clusters to form as a function of the friction coefficient and particle count. 

\item We perform a linear stability analysis of the kinetic (hypoelliptic) McKean-Vlasov PDE, including a study of the fluctuations around the mean-field limit described by the deterministic PDE. Similar to the findings of ~\cite{garnier2017consensus, Garnier_2019}, we conclude that the formation of clusters is related to the breakdown of the central limit theorem that describes fluctuations around the mean-field limit.

\item Motivated by~\cite{GomesSpectral, fok2002combined}, we develop a spectral numerical method to solve the linearized kinetic McKean-Vlasov PDE.

\item {We implement customized, high-accuracy integration schemes for the weak-friction regime, enabling the simulation of significantly larger interacting particle systems ($N \sim 10^3$) than previously analyzed in this context.} The simulation code is published on our \href{https://github.com/SchroedingersLion/Cluster-Formation-in-Diffusive-Systems}{GitHub} repository \cite{GitHub_IPS}. The latter also holds a notebook that demonstrates the numerical solution of the linearized mean-field PDE, as well as illustrative video animations.

\end{itemize}

The rest of the paper is structured as follows. In Section~\ref{sec:LD_theory_intro} we present the interacting particle system, the interaction potentials and the mean-field PDE that we will study in this paper, together with their mathematical properties that will be used later on. In Section~\ref{sec:analysis} we present a linear stability analysis of the mean-field dynamics, around the uniform and Maxwellian distributions in space and velocity, respectively. We derive expressions for the critical temperature and the time until clusters begin to form (the `onset of clustering') which depends on both friction and particle count. Section~\ref{sec:num_analysis} presents our numerical simulations of the $N$-particle system, which confirm our analytical results and yield additional insights into the system behavior. In Section ~\ref{sec:conclusions} we give conclusions and outlooks on future research directions.

\section{Kinetic Langevin dynamics and its mean-field limit}\label{sec:LD_theory_intro}
\subsection{The interacting particle system} \label{sec:IPS_and_mean_field}

We consider a system of weakly interacting hypoelliptic diffusions of the form~\cite{guillin2021uniform, iguchi2025parameterestimationweaklyinteracting}
\begin{equation}\label{eq: Langevin Dynamics}
\begin{split}
d\bm{x}_{i} &= \bm{v}_{i} dt,\\
d\bm{v}_{i} &= -\frac{1}{N}\sum^{N}_{j=1}\nabla W(\bm{x}_{i}-\bm{x}_{j})dt - \gamma \bm{v}_{i}dt + \sqrt{2\gamma\beta^{-1}}d\mathbf{B}^{i}_{t}.
\end{split}
\end{equation}
We consider {random} initial conditions, i.e., the particles are initially independent, distributed according to ${\left(\bm{x}_{i}(0), \, \bm{v}_i(0) \right) \sim \mu_0, \, i=1, \dots, N}$. We will consider the system in a box of side length $L$ with periodic boundary conditions. We assume a Maxwellian distribution for the initial velocities. The initial positions are distributed uniformly. Details are presented in Section~\ref{sec:num_analysis}.

In this work, we consider three pairwise interaction potentials in \eqref{eq: Langevin Dynamics}, given by functions $W(\boldsymbol{x})$, with $\boldsymbol{x} \in \mathbb{R}^d$, $d\in\{1,2\}$, defined as follows:
\begin{equation}\label{eq:Potentials}
\begin{aligned} 
W(\boldsymbol{x}) &= - e^{-\frac{\|\boldsymbol{x} \|^2}{2\sigma^2}},  & \text{Gaussian Potential}, \\
W(\boldsymbol{x}) &= D_e \Big(e^{-2a\|\boldsymbol{x}\|} -2e^{-a\|\boldsymbol{x}\|}\Big),  & \text{Morse potential}, \\
W(\boldsymbol{x}) &= -e^{-\Big(\frac{\|\boldsymbol{x}\|}{\sqrt{2\sigma^2}}\Big)^\alpha},  & \text{GEM-$\alpha$}.
\end{aligned}
\end{equation}
We fix $\sigma^2=0.5$, $a=2$, $D_e=1$, and $\alpha=4$. Note that the generalized exponential model (GEM-$\alpha$) (see \cite{mladek2006formation,delfau2016pattern,miyazaki2019slow}) reduces to the Gaussian potential for $\alpha=2$ and becomes supergaussian for $\alpha=4$. For more details and a plot of the potentials, see Sec. \ref{sup:sec:potentials} in the appendix. {Note that \eqref{eq:Potentials} gives the formulation of the potentials in an unbounded domain. Since the configurational domain we consider is periodic, we would normally require periodic extensions of the potentials instead. However, if the effective interaction range is short enough (which is given for our choices of $\sigma^2$, $a$ and torus period $L$), the distinction between the periodic and unbounded versions is negligible, in both analysis and simulation.}
We remark that, due to the absence of an external potential and {due to the periodic boundary conditions}, the IPS and its mean-field limit are translation invariant~\cite{bertini2010dynamical}. We could, as well, consider the dynamics in the whole space and with a sufficiently strong confining potential, resulting in dynamics that is not translation-invariant. This problem will be taken up elsewhere. We also remark that the potentials in \eqref{eq:Potentials} are all purely attractive, with equal depth and comparable characteristic width (see Sec. \ref{sup:sec:potentials} in the appendix for an illustration). The latter point ensures that their thermodynamic properties such as critical temperatures are all of the same order of magnitude. {We emphasise that discontinuous phase transitions and dynamical metastability are only present when we have sufficiently strong and sufficiently localised attractive interactions, see \cite[Proposition 3.2]{gerber2025formationclusterscoarseningweakly}.} One would expect that clustering phenomena also occur for interaction potentials that have both attractive and repulsive parts, provided that the attractive part of the interaction is sufficiently strong. We will not consider this case in this paper.

\subsection{Propagation of chaos, mean-field PDE, free energy, and invariant measures}

Consider the underdamped Langevin dynamics~\eqref{eq: Langevin Dynamics} and define the empirical measure
\[
\rho^{N}(t,d\bm{x},d\bm{v}) = \frac{1}{N}\sum^{N}_{j=1}\delta_{(\bm{x}_{j}(t),\bm{v}_{j}(t))}(d\bm{x},d\bm{v}).
\]
We take the mean-field limit as $N \to \infty$. Under appropriate assumptions on the interaction potential (that are satisfied by the potentials that we consider in this paper), the empirical measure converges to $\rho(t,\bm{x},\bm{v}) \, d \bm{x} d \bm{v}$ for $t \in [0,T]$. The probability density function $\rho(t,\bm{x},\bm{v})$ satisfies the nonlinear Fokker-Planck equation (the McKean-Vlasov PDE)
\begin{equation}\label{eq: Fokker-Planck}
\begin{split}
\frac{\partial \rho}{\partial t} = &-\bm{v}\cdot \nabla_{\bm{x}}\rho + \nabla_{\bm{v}}\rho\cdot\left[\int  \rho(t,\bm{x}-\bm{y},\bm{v}')\nabla W(\bm{y})d\bm{y}d\bm{v}'\right] \\
&+ \gamma \left[\nabla_{\bm{v}}\cdot\left(\bm{v}\rho\right) + \beta^{-1}\Delta_{\bm{v}}\rho\right].
\end{split}
\end{equation}
We refer to~\cite{Monmarche_2017, guillin2021uniform} and the references therein for the rigorous derivation and for extensions. See also~\cite{duong2018mean, GomesSpectral} for extensions to non-Markovian problems. 

As is well known (see, e.g., \cite{Dolbeault_1999, DuPeZi_2013, bavaud1991equilibrium}), there is a free energy associated with the McKean-Vlasov PDE~\eqref{eq: Fokker-Planck},
\begin{equation}\label{e:free_energy}
\mathcal{F}(\rho) = \beta^{-1} \int \rho \log(\rho) \, d \bm{x} d \bm{v} +\frac{1}{2} \int  \Big( | \bm{v}|^2 + W \star \rho \Big) \rho \,  d \bm{x} d \bm{v},
\end{equation}
{where $\star$ denotes the spatial convolution acting on the spatial marginal density, $(W \star \rho)(\bm{x}) = \int W(\bm{x}-\bm{y}) \left( \int \rho(\bm{y}, \bm{v}') \, d\bm{v}' \right) d\bm{y}$.} Stationary states of the McKean-Vlasov PDE are critical points of the free energy \eqref{e:free_energy}. Invariant measures for the mean-field dynamics and equivalently critical points of~\eqref{e:free_energy} have been thoroughly studied; see, e.g., \cite{bavaud1991equilibrium} and references therein. For example, it is well known that stationary states for the mean-field PDE~\eqref{eq: Fokker-Planck} have a product structure~\cite[Proposition 2]{dressler1987} (see also \cite[Prop. 1]{duong2016stationary} and \cite[Sec. 3]{duong2018mean}). The (possibly non-unique) invariant densities satisfy the integral equation
\begin{equation}\label{eq:dressler}
\rho_{\infty}(\bm{x},\bm{v}) = \frac{1}{Z}\exp\left[-\beta\left(W \star  \rho_{\infty} (\bm{x}) + \frac{1}{2}\|\bm{v}\|^{2}\right)\right],
\end{equation}
with normalization constant $Z = \int \exp\left[-\beta\left(W \star  \rho_{\infty} (\bm{x}) + \frac{1}{2}\|\bm{v}\|^{2}\right)\right]d\bm{x}d\bm{v}$. From Dressler's theorem ~\cite[Proposition 2]{dressler1987}, the invariant density $\rho_{\infty}(\bm{x},\bm{v})$ can be written in the form $\rho_{\infty}(\bm{x},\bm{v}) = \mathcal{M}(\bm{v}) \tilde{\rho}_{\infty}(\bm{x})$, where $\mathcal{M}(\cdot)$ denotes the Maxwellian distribution and $\tilde{\rho}_{\infty}(\bm{x})$ is the solution to the Kirkwood-Monroe integral equation~\cite{KirkwoodMonroe}
\begin{equation}\label{e:kirk_monr}
\tilde{\rho}_{\infty}(\bm{x}) = \frac{1}{\tilde{Z}} e^{- \beta W \star \tilde{\rho}_{\infty}(\bm{x})},
\end{equation}
where $\tilde{Z}$ denotes the configurational partition function. This integral equation is the main object of study in~\cite{ChayesPanferov2010,Carrillo2020}; it provides information about the number of and nature of stationary states.

\subsection{Phase transition for the kinetic McKean-Vlasov PDE}

Since invariant measures of the kinetic McKean-Vlasov PDE are product measures, steady states are independent of $\gamma$. In particular, considering the characterization of the steady states~\eqref{eq:dressler} and the free energy \eqref{e:free_energy}, we have that critical points of the free energy~\eqref{eq:dressler} are also critical points of the free energy associated with the overdamped system, which is given by
\begin{equation}\label{e:free_energy_ov}
\mathcal{F}_{\textnormal{ov}}(\rho) = \beta^{-1} \int \rho \log(\rho) \, d \bm{x} +\frac{1}{2} \int  \Big(  W \star \rho \Big) \rho \,  d \bm{x}.
\end{equation}
The Euler-Lagrange equation corresponding to this functional is precisely the integral equation~\eqref{e:kirk_monr}. An important consequence of this observation is that the overdamped and  kinetic/underdamped dynamics exhibit the same phase transitions; furthermore, the critical temperature at which the transition occurs is independent of the friction coefficient $\gamma$. {This is expected from a statistical mechanics perspective: both underdamped and overdamped Langevin equations sample the canonical ensemble, where the equilibrium state is defined by the Hamiltonian and the temperature $\beta^{-1}$. The friction coefficient $\gamma$ merely facilitates the coupling to the heat bath and determines the dynamical relaxation rate; it does not shift the thermodynamic boundaries (such as $\beta_c^{-1}$) where the phase transition occurs.}  On the other hand, and this is crucial for the work presented in this paper, the transient dynamics, that leads from the initial uniform-in-space configuration to the formation of a single cluster, depends on the value of the friction coefficient.  

A detailed analysis of phase transitions for the overdamped dynamics, i.e., of the number and nature of critical points of the free energy~\eqref{e:free_energy_ov}, is presented in~\cite{Carrillo2020}. From the analysis presented in this article and, in particular from~\cite[Cor. 5.14]{Carrillo2020}, it follows that for sufficiently localized attractive potentials, like the ones that we consider in this paper, the free energy exhibits a discontinuous phase transition.\footnote{In fact, in this result it is required that the interaction potential is normalized, $\|W \|_{L^1} = 1$. However, as shown in~\cite{GvalSch_2020} this normalization assumption may be avoided.}  This means that, at sufficiently low temperatures, and for all values of the friction coefficient, the uniform distribution in position becomes unstable and a localized stationary state emerges. In addition, it is well known that systems exhibiting discontinuous phase transitions exhibit dynamical metastability. For the kinetic Langevin dynamics considered in this paper, this dynamical metastability manifests itself in the formation of clusters. It is our goal to understand the breakdown of the uniform distribution and the formation of clusters. We mention in passing that, for the proof of ~\cite[Cor. 5.14]{Carrillo2020}, only the attractive part of the interaction, i.e. the negative Fourier coefficients, need to be considered. In particular, we can extend the analysis also to potentials with both attractive and repulsive parts. We will leave this for future work.
\begin{comment}
The non-trivial critical points of \eqref{e:free_energy_ov} of the corresponding overdamped system have been studied rigorously in \cite{Carrillo2020}.

Now, considering the critical points of \eqref{e:free_energy_ov}, we have due to \cite[Theorem 1.3]{Carrillo2020}, non-uniqueness of critical points of the corresponding free-energy under appropriate assumptions on the interaction potential (which hold for the considered potentials in this paper) and for sufficiently localized attractive potentials we have existence of discontinuous phase transitions. In combination with \cite[Proposition 2]{dressler1987}, as we have discussed, we have non-uniqueness of critical points of \eqref{e:free_energy} and, for sufficiently localized attractive potentials, existence of discontinuous phase transitions for \eqref{eq: Fokker-Planck}.

\begin{remark}
 Although the potentials we consider do not have repulsive components, only the attractive part of the potentials (the negative Fourier modes) are relevant for the proof of the existence of a discontinuous phase transition in \cite{Carrillo2020}.
\end{remark}

\end{comment}

\section{Analysis of the kinetic McKean-Vlasov PDE}\label{sec:analysis}
\subsection{Linear stability analysis} \label{sec: linear_stability}
Let the particle position domain be given by $\mathcal{D}_{\boldsymbol{x}}=\mathbb{T}^d_L=[0,L)^d$. The velocity domain is simply $\mathcal{D}_{\boldsymbol{v}}=(-\infty, \infty)^d$.
We know that there is a discontinuous phase transition in the mean-field limit, governed by the temperature parameter $\beta^{-1}$. We may then perform linear stability analysis of ($\ref{eq: Fokker-Planck}$) to provide insight into the critical temperature {$\beta^{-1}$ at which the uniform state reaches its linear instability threshold.} 
For this aim, we first derive the Fokker-Planck equation for the perturbation in $d$ dimensions, but then restrict our analysis to the one-dimensional case. 
We consider a small perturbation of the equilibrium measure, the latter given by the uniform measure in position and Gaussian in velocity. More precisely, we consider a density of the form 
\begin{equation} \label{eq: linearized_density}
\rho(t,\bm{x},\bm{v}) = \frac{1}{L^d}F(\bm{v}) + \rho_{1}(t,\bm{x},\bm{v}),
\end{equation}
with $F(\bm{v})=(2\pi \beta^{-1})^{-\frac{d}{2}}e^{-\beta \frac{\| \bm{v} \|^2 }{2}}$, and where the perturbation $\rho_1$ is assumed to be small such that terms in  $\mathcal{O}(\rho^{2}_{1})$ are negligible. 

Linear stability analysis has been performed in the overdamped setting in \cite{garnier2017consensus}, where they perform a Fourier expansion of the perturbation, $\rho_{1}$, diagonalizing the infinite-dimensional linear system. Consequently, they are able to derive a precise formula for the maximal growth rate across the Fourier modes. If the maximal growth rate is positive, the system is unstable. Their formula allows them to provide an expression for the critical temperature $\beta_{c} > 0$ for which the system is unstable for any $\beta > \beta_{c}$.

In the kinetic Langevin setting, the linear stability analysis is substantially more difficult, requiring us to perform a Fourier expansion in position, $\bm{x}$, and a Hermite expansion in velocity, $\bm{v}$. Rather than diagonalizing the linear system, we are able to tridiagonalize the linear system via this approach, allowing for efficient estimation of the maximal growth rate via numerical methods  (among Fourier and Hermite coefficients). This allows us to compute the critical temperature $\beta_{c} > 0$ for which the system is unstable for any $\beta > \beta_{c}$, as in the overdamped system. 

We remark that, in the development of spectral methods for the numerical simulation of the kinetic Fokker-Planck equation,  it is common  to perform Hermite expansion in the velocity component (see \cite{fok2002combined,GomesSpectral,avelin2023galerkin}).

We proceed by inserting the ansatz $\rho=\frac{1}{L^d}F(\bm{v})+\rho_1(t,\bm{x},  \bm{v})$ into the kinetic Fokker-Planck equation
(\ref{eq: Fokker-Planck}). {When neglecting all terms in $\mathcal{O}(\rho^{2}_{1})$, this leads to}
\begin{equation}\label{eq: Fokker-Planck-linear}
\begin{split}
\frac{\partial \rho_1}{\partial t} = &-\bm{v}\cdot \nabla_{\bm{x}}\rho_1 +  \frac{1}{L^d}\nabla_{\bm{v}}F \cdot \left[\int_{\mathcal{D}_{\bm{v}'}} \int_{\mathcal{D}_{\bm{y}}} \rho_1(t,\bm{x}-\bm{y},\bm{v}')\nabla W(\bm{y})d\bm{y}d\bm{v}'\right] \\
&+ \gamma \nabla_{\bm{v}} \cdot \left( \bm{v}\rho_1 \right) + \gamma \beta^{-1} \Delta_{\bm{v}} \rho_1,
\end{split}
\end{equation}
where the integration domains are given by $\mathcal{D}_{\bm{v}'}={\mathbb{R}}^d$ and $\mathcal{D}_{\bm{y}}:=\left[-\frac{L}{2}, \frac{L}{2}\right)^d$, see Sec. \ref{sup:sec:linear_stability1} in the appendix for a detailed calculation. Note that $\mathcal{D}_{\boldsymbol{y}} \neq \mathcal{D}_{\boldsymbol{x}}=[0,L)^d$, since $\bm{y}$ plays the role of a {difference} between two positions, which lies in $\left[-\frac{L}{2}, \frac{L}{2} \right]^d$ by the minimum image convention.

The spatial spectral components of $\rho_1$ are given by
$$\hat{\rho}_1(t,\bm{k},\bm{v}) = \int_{\mathcal{D}_{\boldsymbol{x}}}\rho_1(t, \bm{x}, \bm{v}) e^{-i\bm{k}\cdot \bm{x}} d\bm{x}, $$ with $\bm{k}\in\left\{\left(\frac{2\pi n_1}{L}, ...\frac{2\pi n_d}{L}\right)\middle| n_1,...,n_d \in \mathbb{Z} \right\}$. {Taking the time-derivative and using (\ref{eq: Fokker-Planck-linear}), we follow the steps outlined in Sec. \ref{sup:sec:linear_stability2} of the appendix to obtain}
\begin{equation}\label{eq: linear_stability_spectral}
\begin{split}
\frac{\partial \hat{\rho}_{1}}{\partial t} &= -i\bm{k}\cdot\left[\bm{v} \hat{\rho}_1 -  {\frac{1}{L^{d}}}\nabla_{\bm{v}}F \,\hat{W} \int_{\mathcal{D}_{\bm{v'}}} \hat{\rho}_1d\bm{v'}\right] + \gamma \nabla_{\bm{v}} \cdot \left(\hat{\rho}_1 \bm{v} \right) + \gamma \beta^{-1} \Delta_{\bm{v}} \hat{\rho}_1,
\end{split}
\end{equation}
with $\hat{W}(\bm{k}):=\int_{\mathcal{D}_{\bm{y}}} W(\bm{y})e^{-i\bm{k}\cdot \bm{y}}d\bm{y}$ the Fourier transform of $W(\bm{y})$. 

In the following, we restrict ourselves to the one-dimensional case to simplify computation. 
We define $\Tilde{\rho}_{1} := \hat{\rho}_{1}/F$ and consider its expansion in terms of an orthonormal basis of (normalized) Hermite polynomials $(h_{n})_{n \in \mathbb{N}}$ (see Sec. \ref{sup:sec: Hermite_polynomials} for more information), i.e.,
\begin{equation} \label{eq:hermite_expansion}
\Tilde{\rho}_{1} := \sum^{\infty}_{n=0}c_{n}(t,k)h_{n}(v).
\end{equation}
The functions $h_n$ satisfy an orthonormality condition 
\begin{equation}\label{eq: orthonormality}
\int_{-\infty}^{\infty} F h_n h_m dv=\delta_{n,m},
\end{equation}
for all $n,m \in \mathbb{N}_0$.

Inserting the expansion into (\ref{eq: linear_stability_spectral}), we have
\begin{equation}\label{eq: hermite expansion}
\begin{split}
&\sum^{\infty}_{n=0}F(v)h_{n}(v)\frac{\partial c_{n}(t,k)}{\partial t} =\\
&\sum^{\infty}_{n=0} \Bigg\{-ik \left( v F(v)h_{n}(v) +  \frac{\beta}{{L}} \hat{W}(k) v F(v) \int_{-\infty}^{\infty} F(v')h_{n}(v')dv' \right)\\
& \qquad \quad+ \gamma\left[\frac{\partial}{\partial v} \left(v F(v)h_{n}(v)\right) + \beta^{-1}\frac{\partial^2}{\partial v^2} (F(v)h_{n}(v))\right]\Bigg\} \, c_{n}(t,k),
\end{split}
\end{equation}
where in the last step we inserted $\partial F/ \partial v =-\beta v F(v)$.
We now observe that the last term in (\ref{eq: hermite expansion}) is the Fokker-Planck operator $\mathcal{L}_{OU}^*$ of the OU process applied to $Fh_n${, where $\mathcal{L}_{OU}^*$ applied to test function $f(v)$ is given by
\begin{equation}
\mathcal{L}_{OU}^*f(v)=\frac{\partial}{\partial v} \big(v f(v)\big)+ \beta^{-1}\frac{\partial^2}{\partial v^2} f(v).
\end{equation}}
We use the identity $\mathcal{L}_{OU}^*(Fg)=F\mathcal{L}_{OU}g $ together with the fact that the normalized Hermite polynomials $\left(h_{n}\right)_{n\in\mathbb{N}}$ are eigenfunctions of the generator of the OU process, 
\begin{equation}
-\mathcal{L}_{OU}h_n = \lambda_n h_n, \; \text{with $\lambda_n=n$},
\end{equation}
to obtain
\begin{equation*}
\begin{split}
&\sum^{\infty}_{n=0}F(v)h_{n}(v)\frac{\partial c_{n}(t,k)}{\partial t} \\
&= \sum^{\infty}_{n=0}\Bigg\{-ik \left( v F(v)h_{n}(v) +  \frac{\sqrt{\beta}}{{L}}\hat{W}(k)F(v)h_1(v)\delta_{n,0} \right) - n \gamma F(v) h_{n}(v)\Bigg\} \, c_{n}(t,k),\\
\end{split}
\end{equation*}
where we solved the integral by writing the integrand as $F(v)h_n(v)h_0(v)$ and using the orthonormality condition (\ref{eq: orthonormality}). We also replaced $v$ in the prefactor of the integral by $\frac{1}{\sqrt{\beta}}h_1(v)$.
After multiplying both sides by $h_{m}(v)$ and integrating with respect to $v$, we have, again  due to (\ref{eq: orthonormality}),
\begin{equation*}
\begin{split}
&\frac{\partial c_{m}(t,k)}{\partial t} =\\
& -ik \sum^{\infty}_{n=0} \left(\int_{-\infty}^{\infty} vh_{n}(v) h_{m}(v)F(v)dv\right) c_{n}(t,k) - \frac{ik \sqrt{\beta}}{{L}}\hat{W}(k)\delta_{m,1}c_{0}(t,k) -m\gamma c_{m}(t,k ).
\end{split}
\end{equation*}
The integral evaluates to (see Sec. \ref{sup:sec:linear_stability3} in the appendix)
$$
\int_{-\infty}^{\infty} vh_{n}(v) h_{m}(v)F(v)dv=\sqrt{\frac{m!}{n!\beta}}\Big\{ \delta_{m,n+1} + (m+1)\delta_{m,n-1} \Big\}.
$$
Finally, we arrive at 
\begin{equation*}
\begin{split}
\frac{\partial c_{m}}{\partial t} = - \frac{ik \sqrt{\beta}}{{L}}\hat{W}(k)\delta_{m,1}c_0 -i k \sqrt{\frac{m}{\beta}} c_{m-1} - m\gamma c_m -i k \sqrt{\frac{m+1}{\beta}}c_{m+1} \quad \text{for $m=0,1,2,...$}.
\end{split}
\end{equation*}
In matrix form, this becomes 
\begin{equation}\label{c-ODE}\frac{\partial}{\partial t}\bm{c}(t,k)=\mathbf{A}_{\gamma}(k,\beta)\bm{c}(t,k),
\end{equation}
with $\bm{c}(t,k):=\big(c_0(t,k),\ c_1(t,k),\ c_2(t,k),...\big)^T$ and
\begin{equation*}\mathbf{A}_{\gamma}(k,\beta)=
\resizebox{.8\hsize}{!}{$\begin{bmatrix}
0 & -ik\sqrt{\frac{1}{\beta}} & 0 & 0 & \cdots \\
-i k \sqrt{\frac{1}{\beta}} - \frac{ik \sqrt{\beta}}{{L}}\hat{W}(k)& -\gamma & -ik\sqrt{\frac{2}{\beta}} & 0  & \cdots  \\
0 & -ik\sqrt{\frac{2}{\beta}} & -2\gamma & -ik\sqrt{\frac{3}{\beta}} &  \cdots  \\
0& 0 & -ik\sqrt{\frac{3}{\beta}} & -3\gamma   \ddots \\
0& 0 & 0 & -ik\sqrt{\frac{4}{\beta}} &  \ddots  \\
\vdots & \vdots &  & \ddots & \ddots & \ddots  \\
 &  &   &  & -ik\sqrt{\frac{m}{\beta}} & -m\gamma & -ik\sqrt{\frac{m+1}{\beta}}  \\
 &  &  &  &  & \ddots & \ddots & \ddots 
\end{bmatrix}$}.
\end{equation*}
{The uniform state of the mean-field system described by \eqref{eq: Fokker-Planck} loses local stability when the linear stability analysis is no longer valid \cite{garnier2017consensus}. In particular, this linear instability threshold is crossed if the linearization \eqref{eq: linearized_density} is unstable.} This happens when the eigenvalue of the linear system \eqref{c-ODE}  with largest real part is positive for some frequency $k$, i.e.,
\[
\psi_{\max}(\beta, \gamma) := \max_{k}\left[\psi^{k}_{\max}(\mathbf{A}_{\gamma}(k,\beta))\right] :=  \max_{k}\left[\max\{\textnormal{Re}(\lambda) \mid \lambda \textnormal{ is an eigenvalue of }\mathbf{A}_{\gamma}(k,\beta)\}\right]>0.
\]
Due to the tridiagonal structure of the matrix $\mathbf{A}_{\gamma}(k,\beta)$ in \eqref{c-ODE}, this quantity is difficult to estimate analytically and we use spectral methods to numerically approximate $\psi_{\max}(\beta, \gamma)$ for the considered potentials.
{Although the growth rate \(\psi_{\max}\) depends on \(\gamma\), the linear instability threshold in \(\beta\) is observed, and is expected on thermodynamic grounds, to be independent of \(\gamma\).}
We thus define the critical (inverse) temperature $\beta_{c}$ by
\begin{equation}\label{eq:beta_crit}
{\beta_{c} := \inf\{\beta > 0\,:\,\psi_{\max}(\beta, \gamma)>0\},}
\end{equation}
the smallest $\beta$ for which the linear system becomes unstable.  
\subsection{Spectral approximation for the linearized PDE}\label{sec:spectral_approx}

\begin{figure}
\centering
\includegraphics[width=0.45\textwidth]{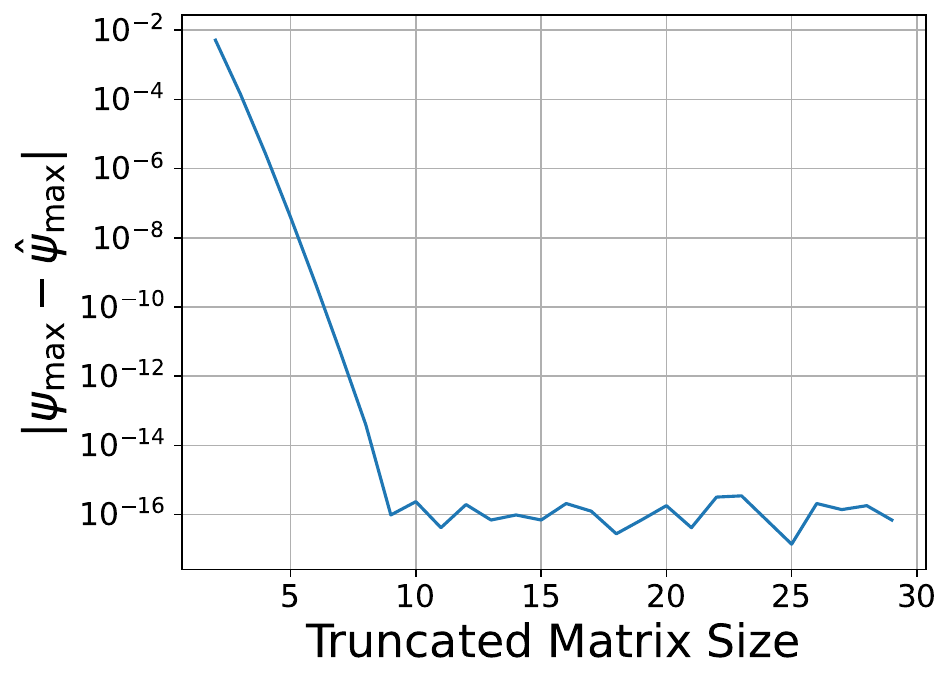}
\hspace{0.5cm}
\includegraphics[width=0.45\textwidth]{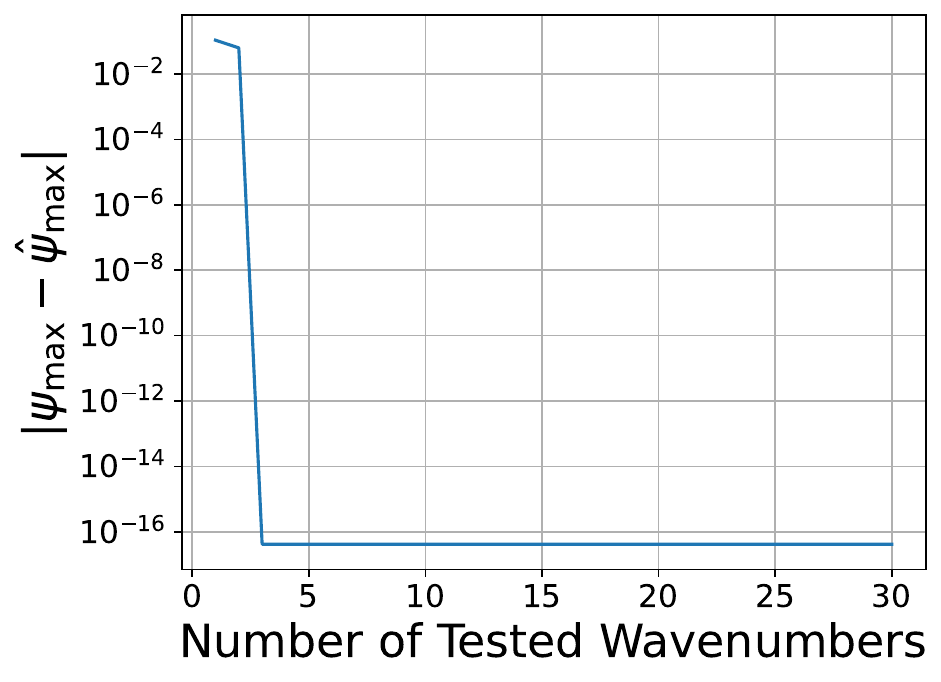} 
\caption{Convergence of $\psi_{\max}$ obtained from $\boldsymbol{A}_{\gamma}(k,\beta)$ when considering increasing numbers of matrix dimensions (left) and wavenumbers (right). Ground truth estimate $\hat{\psi}_{\max}$ was obtained by truncating after 100 matrix rows and using 30 wavenumbers (plus their negative images). Interaction potential \eqref{eq: wrapped_gaussian} for $\sigma = 1/\sqrt{2}$, $L = 10$, $\beta = 25 > \beta_{c}$ and $\gamma = 1.0$.}
\label{fig:spectral_conv}
\end{figure}

\begin{figure}[H]
\centering
\includegraphics[width=0.45\textwidth]{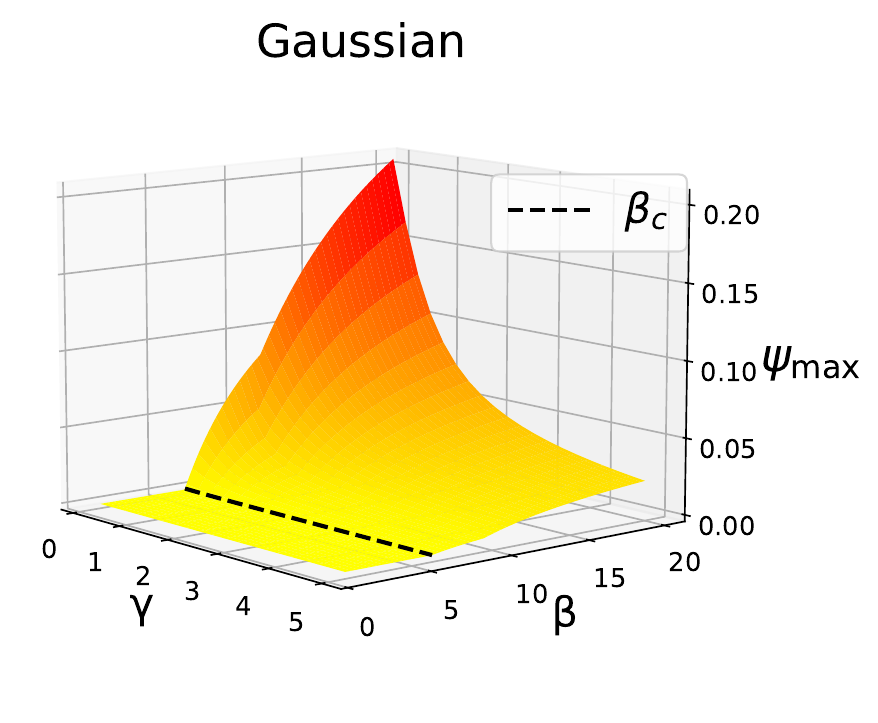} 
\includegraphics[width=0.45\textwidth]{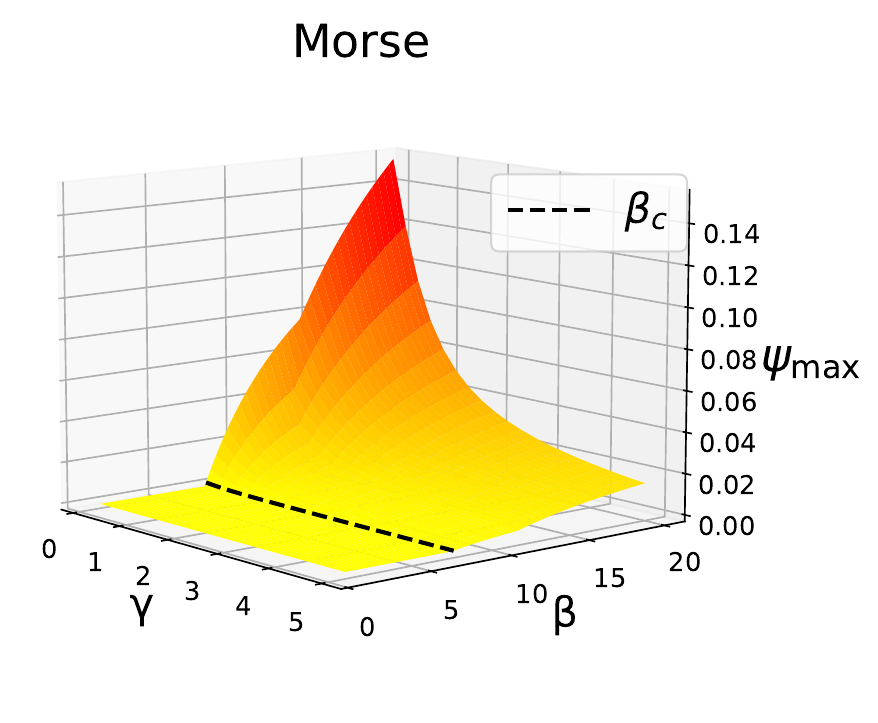} 
\includegraphics[width=0.45\textwidth]{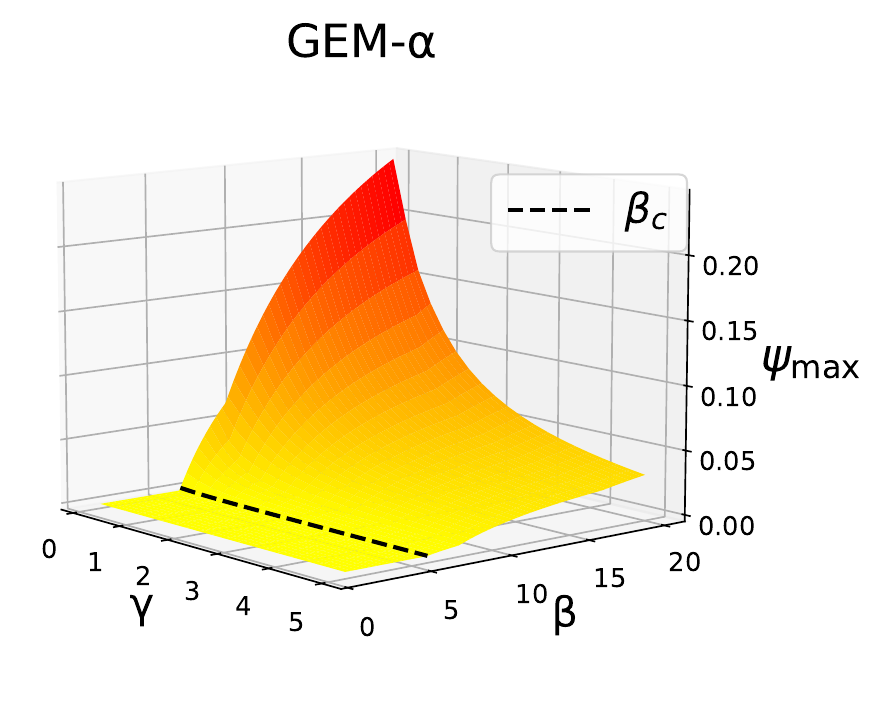} 

\caption{The relationship between friction, inverse temperature and $\psi_{\max}$, illustrating a discontinuous phase transition at $\beta_{c}$. \textbf{{Upper left:}} Gaussian potential for $\sigma = 1/\sqrt{2}$ and $L = 10$. \textbf{{Upper right:}} Morse potential for $a = {2}$, $D_{e} = {1}$ and $L = 10$. \textbf{{Bottom:}} GEM-$\alpha$ potential for $\sigma = 1/\sqrt{2}$, $\alpha  =4$ and $L = 10$.}
\label{fig:spectral}
\end{figure}

To estimate $\psi_{\max}$ for a particular interaction potential, we truncate the infinite-dimensional matrix $\mathbf{A}_{\gamma}(k,\beta)$ and only consider a finite number of wavenumbers $k$. To determine how many {matrix} dimensions and wavenumbers are required for reliable results, we first compute $\psi_{\max}$ for varying numbers of dimensions while keeping the number of wavenumbers large and fixed. We then fix the number of dimensions at a large value and compute $\psi_{\max}$ for varying numbers of wavenumbers. The ground truth $\hat{\psi}_{\max}$ is obtained by using both large numbers of considered dimensions and wavenumbers, where it was experimentally verified that increasing these numbers substantially gave no further alteration to the result. We can then plot the errors $|\psi_{\max}-\hat{\psi}_{\max}|$ against the number of considered matrix rows and wavenumbers. To evaluate the matrix $\mathbf{A}_{\gamma}(k,\beta)$, one has to compute the Fourier transform of the interaction potential $W(r)$, which in general needs to be done numerically. One could, in principle, work with periodic extensions of these potentials, whose Fourier transforms would then be obtainable analytically. For example, one could consider the cyclical Gaussian, given by 
\begin{equation}\label{eq: wrapped_gaussian}
W(r) = -\sum^{\infty}_{n=-\infty}\exp{\left(-\frac{\left(r+Ln\right)^{2}}{2\sigma^{2}}\right)}.
\end{equation}
This potential is periodic on a  torus of length $L$. In Fig. \ref{fig:spectral_conv} we show the results for $|\psi_{\max}-\hat{\psi}_{\max}|$ when using interaction \eqref{eq: wrapped_gaussian}.
We observe rapid decay of the error in both considered matrix rows and wavenumbers, making it numerically feasible to accurately evaluate $\boldsymbol{A}_{\gamma}(k,\beta)$.\\
Since the terms for $|n|>1$ in the wrapped Gaussian \eqref{eq: wrapped_gaussian} are negligible for the parameters we consider, $\sigma^2=0.5$ and $L=10$, the results obtained for the non-periodic Gaussian are identical, and we do not report them here. However, the interested reader may run their own trials with both implemented versions with the Jupyter notebook we provide on our GitHub page \cite{GitHub_IPS}, together with the other interaction potentials, all of which lead to similar convergence as in Fig. \ref{fig:spectral_conv}. {Although wrapped potentials like \eqref{eq: wrapped_gaussian} are mathematically rigorous, $N$-particle simulations achieve greater computational efficiency by evaluating the non-periodic form via the minimum-image convention.} In the following, we use the non-periodic interactions given in Sec. \ref{sec:IPS_and_mean_field} and obtain their Fourier transforms via numerical quadrature. In each of the experiments, we pick sufficient numbers of matrix rows and wavenumbers, based on the error convergences shown in Fig. \ref{fig:spectral_conv}.
We proceed by computing the spectral abscissa $\psi_{\max}$ for a range of friction $\gamma$ and the inverse temperature $\beta$. The results are illustrated in Fig. \ref{fig:spectral}. For all three models our results show $\psi_{\max} = 0$ for $\beta \leq \beta_{c}$, where $\beta_c$ is illustrated by the black dashed lines. {Hence, for $\beta \leq \beta_{c}$, the uniform state remains linearly stable. Then, for $\beta > \beta_{c}$, $\psi_{\max} > 0$ and the uniform state of the Fokker-Planck equation becomes locally unstable, triggering spontaneous cluster formation.} Our results imply that the critical temperature $\beta^{-1}_{c}$ is independent of friction $\gamma$. For fixed $\gamma$, there is a classical pitchfork bifurcation. The critical temperatures are given in Table \ref{tab:beta_crit}. 

\begin{table}[!htbp]
    \centering
    \begin{tabular}{|c||c|c|c|}
        \hline
        Potential & Gaussian & Morse & GEM-4\\
        \hline
        $\beta_c$  &  6.23   &  7.51  &  5.90\\
        \hline
    \end{tabular}
    \caption{The critical (inverse) temperatures for the three interaction potentials in one dimension obtained by \eqref{eq:beta_crit}.}
    \label{tab:beta_crit}
\end{table}

\subsection{Study of fluctuations around the mean-field limit}\label{sec:fluc_analysis}
As in \cite{garnier2017consensus}, when $\beta > \beta_{c}$ we wish to identify the time $t_{\textnormal{cl}}$ when the central limit theorem breaks down, i.e., when the size of the fluctuations becomes comparable to the size of the domain. $t_{\textnormal{cl}}$ will then be the time where the system, starting from a uniform distribution, will begin to form clusters.

We assume that $\beta > \beta_{c}$, and the linear system is unstable. We consider $N$ particles whose initial conditions $({x}_{1}(0),{v}_{1}(0))$, \dots $({x}_{N}(0),{v}_{N}(0))$ are i.i.d and sampled from a measure $\rho_{0}(d{x},d{v})$ with density $\rho_{0}= \frac{1}{L}F({v})$.
{The empirical measure at time $t=0$ can be written as  
$$
\rho^N(0,dx,dv) \approx \rho_{0}(dx,dv) + \frac{1}{\sqrt{N}} \mu(dx,dv),
$$
i.e., from the central limit theorem it follows that the fluctuations converge to a Gaussian-valued measure $\mu$. 
We refer to \cite[Sec. 2.3]{gerber2025formationclusterscoarseningweakly} for a rigorous treatment and extension to non-uniform initialisation in position.
}
We then have that 
\begin{align} \label{eq: fluctuation_measure}
{\mu}^{N}(0,d{x},d{v}) &:= \sqrt{N}\left(\rho^{N}(0,d{x},d{v}) - \rho_{0}(d{x},d{v})\right), \\
&= \sqrt{N}\left(\frac{1}{N}\sum^{N}_{j=1}\delta_{{x}_{j}(0),{v}_{j}(0)}(d{x},d{v}) - \frac{1}{L}F({v})d{x}d{v}\right)
\end{align}
converges in distribution as $N \to \infty$ to a measure ${\mu}(t=0,dx,dv)$ with density ${\mu}(0, x, v)$. Just like in \eqref{eq:hermite_expansion} of Sec. \ref{sec: linear_stability}, the Fourier components of ${\mu}$ may then be expanded in terms of Hermite coefficients $c_{n}(0,k)$, which satisfy
\begin{align*}
c_{n}(t=0,k) &= \lim_{N \to \infty} \int^{\infty}_{-\infty}\int^{L/2}_{-L/2}\sqrt{N}e^{-ikx}h_{n}(v)\left(\frac{1}{N}\sum^{N}_{j=1}\delta_{x_{j}(0)}(dx)\delta_{v_{j}(0)}(dv)-\frac{F(v)dxdv}{L}\right),\\
&= \lim_{N \to \infty }\sqrt{N}\left(\frac{1}{N}\sum^{N}_{j=1}e^{-ikx_{j}(0)}h_{n}(v_{j}{(0)}){-\delta_{k0}\delta_{n0}}\right).
\end{align*}
{For the nonzero modes, the coefficients are centred complex circular Gaussian random variables with covariance}
\begin{equation*}
\mathbb{E}\left[c_{n}(t=0,k){\overline{c_{m}(t=0,k')}}\right] = \delta_{kk'}\delta_{nm}.
\end{equation*}

{For $t\geq 0$, define}
\begin{equation*}
{
\mu^N(t,d{x},d{v})
:=
\sqrt{N}\left(\rho^N(t,d{x},d{v})-\frac{1}{L}F(v)d{x}d{v}\right).
}
\end{equation*}
For any $T<\infty$, the measure-valued process {$\mu^N(t,d{x},d{v})$} converges to ${\mu}(t,d{x},d{v})$ whose density satisfies a stochastic PDE (see \cite{garnier2017consensus,Dawson1983,muller2024well})

\begin{align*}
d{\mu} &= -{v}\cdot \nabla_{{x}}{\mu}\, dt +  \frac{1}{L}\nabla_{{v}}F \cdot \left[\int_{\mathcal{D}_{{v}'}} \int_{\mathcal{D}_{{y}}} {\mu}(t,{x}-{y},{v}')\nabla W({y})d{y}d{v}'\right]dt \\
& \quad + \gamma \nabla_{{v}} \cdot \left( {v}{\mu} \right)dt + \gamma \beta^{-1} \Delta_{{v}} {\mu} \, dt + \sqrt{2\gamma\beta^{-1}}d\eta,
\end{align*}
with $\mathcal{D}_{{v}'}$ and $\mathcal{D}_{{y}}$ as in Sec. \ref{sec: linear_stability}. $\eta(t,{x},{v})$ is a Gaussian process with mean zero and covariance
\begin{equation}\label{eq:Cov_op}
\begin{split}
    &\textnormal{Cov}\left(\int^{L/2}_{-L/2}\int^{\infty}_{-\infty}\eta(s,{x},{v})f_{1}({x},{v})d{x}d{v},\overline{\int^{L/2}_{-L/2}\int^{\infty}_{-\infty}\eta(t,{x},{v})f_{2}({x},{v})d{x}d{v}}\right)\\
    &= \frac{1}{L}\min(s,t)\int^{L/2}_{-L/2}\int^{\infty}_{-\infty}\frac{\partial f_{1}}{\partial {v}}({x},{v})\overline{\frac{\partial f_{2}}{\partial {v}}({x},{v})}F({v})d{x}d{v},
\end{split}
\end{equation}
for any test functions $f_{1}$, $f_{2}$. \\

Following the preceding linear stability analysis, we first take the Fourier transform in $x$ of the SPDE and have for wavenumber $k$
\begin{equation*}\label{eq: linear stability, spectral}
\begin{split}
d \hat{{\mu}} &= -i{k}\cdot\left[{v} \hat{{\mu}} -  {\frac{1}{L}}\nabla_{{v}}F \, \hat{W} \int_{\mathcal{D}_{v'}}  \hat{{\mu}} \, d{v'}\right]dt + \gamma \nabla_{{v}} \cdot \left(\hat{{\mu}} {v} \right)dt\\
&+ \gamma \beta^{-1} \Delta_{{v}} \hat{{\mu}}\, dt + \sqrt{2\gamma \beta^{-1}}\int_{\mathcal{D}_{{y}}} e^{-i{k}\cdot {y}}d\eta(t,{y},{v}).
\end{split}
\end{equation*}

Now consider the Hermite expansion of {$\hat{\mu}$} and multiplying both sides by $h_{m}(v)$ and integrating over $v$ we have for each wavenumber $k$
\begin{equation}\label{eq:hermite_SDE}
\begin{split}
d c_{m} &= - \frac{ik \sqrt{\beta}}{{L}}\hat{W}(k)\delta_{m,1}c_0 dt -i k \sqrt{\frac{m}{\beta}} c_{m-1}dt - m\gamma c_m dt -i k \sqrt{\frac{m+1}{\beta}}c_{m+1} dt \\
&+ \sqrt{2\gamma \beta^{-1}}\int_{\mathcal{D}_{v}}\int_{\mathcal{D}_{{y}}} e^{-i{k}\cdot {y}}h_{m}(v)d\eta(t,{y},{v}) \quad \text{for $m=0,1,2,...$}.
\end{split}
\end{equation}
Then, considering the covariance of the noise \eqref{eq:Cov_op}, we have that the {stochastic integral in \eqref{eq:hermite_SDE}, before multiplication by $\sqrt{2\gamma\beta^{-1}}$,} is a complex-valued Brownian motion in time $t$ scaled by variance
\begin{align*}
    \frac{t}{L}\int^{\infty}_{-\infty}\int^{L/2}_{-L/2}\frac{\partial h_{m}(v)}{\partial v} \frac{\partial h_{m}(v)}{\partial v}F(v)dxdv &=  \frac{t}{L}\int^{\infty}_{-\infty}\int^{L/2}_{-L/2}{\beta m} h_{m-1}(v)h_{m-1}(v)F(v)dxdv \\
    &=  {\beta mt},
\end{align*}
using the orthonormality conditions of the Hermite polynomials. If we consider the covariance of the complex-valued Brownian motion in time $t$, we have for $m \neq n$
\begin{align*}
    \frac{t}{L}\int^{\infty}_{-\infty}\int^{L/2}_{-L/2}\frac{\partial h_{m}(v)}{\partial v} \frac{\partial h_{n}(v)}{\partial v}F(v)dxdv &=  \frac{t}{L}\int^{\infty}_{-\infty}\int^{L/2}_{-L/2}{\beta \sqrt{mn}} h_{m-1}(v)h_{n-1}(v)F(v)dxdv, \\
    &=  0.
\end{align*}
Therefore we have for each $k$ the following infinite-dimensional SDE
\begin{equation*}\label{c-SDE}
d\bm{c}(t,k)=\mathbf{A}_{\gamma}(k,\beta)\bm{c}(t,k)dt + {\sqrt{\gamma}}\mathbf{B}\, d(\bm{W}^{(k)}_t + i\Tilde{\bm{W}}^{(k)}_t),
\end{equation*}
where $\mathbf{B}$ is a diagonal matrix with entries {$\mathbf{B}_{\ell j}=\sqrt{\ell}\,\delta_{\ell j}, \ell,j= 0,1,2,...$}
The solution of this SDE is given by
\begin{equation*}
    \bm{c}(t,k) = e^{t\bm{A}(k)}\bm{c}(0,k) + {\sqrt{\gamma}}\int^{t}_{0}e^{(t-s)\bm{A}(k)}\bm{B}d(\bm{W}^{(k)}_{s} + i\Tilde{\bm{W}}^{(k)}_{s}),
\end{equation*}
{where for brevity we write $\bm{A}(k)$ for $\bm{A}_{\gamma}(k, \beta)$. For the nonzero modes, since $\bm{c}(0,k)$ has mean zero and identity covariance, we have}
\begin{align*}
{\bm{C}_{k}(t)}
&{:= \mathbb{E}\left[\bm{c}(t,k)\bm{c}(t,k)^{*}\right],} \\
&{= e^{t\bm{A}(k)}e^{t\bm{A}(k)^{*}}
+ 2\gamma \int^{t}_{0} e^{(t-s)\bm{A}(k)} \bm{B}\bm{B}^{T}e^{(t-s)\bm{A}(k)^{*}}ds.}
\end{align*}

{Since cluster formation is a spatial phenomenon, we consider the spatial marginal fluctuation}
\begin{equation*}
{\bar{\mu}(t,x):=\int_{\mathcal{D}_v}\mu(t,x,v)dv.}
\end{equation*}
{Using the Hermite expansion,}
\begin{align*}
    {\hat{\bar{\mu}}(t,k)}
    &{= \int_{\mathcal{D}_v}\hat{\mu}(t,k,v)dv,} \\
    &{= \int_{\mathcal{D}_v}\sum_{n}c_{n}(t,k)h_{n}(v)F(v)dv,} \\
    &{= c_0(t,k),}
\end{align*}
{where we used $\int_{\mathcal{D}_v}h_n(v)F(v)dv=\delta_{n,0}$. Hence}
\begin{equation*}
{
\mathbb{E}\left[\hat{\bar{\mu}}(t,k)\overline{\hat{\bar{\mu}}(t,k)}\right]
=
\mathbb{E}\left[c_0(t,k)\overline{c_0(t,k)}\right]
=
\left(\bm{C}_{k}(t)\right)_{00}.
}
\end{equation*}
{Equivalently,}
\begin{align*}
{\mathbb{E}\left[c_0(t,k)\overline{c_0(t,k)}\right]}
&{=
\left(e^{t\bm{A}(k)}e^{t\bm{A}(k)^{*}}\right)_{00}} \\
&{\quad + 2\gamma \int^{t}_{0}
\left(e^{(t-s)\bm{A}(k)} \bm{B}\bm{B}^{T}e^{(t-s)\bm{A}(k)^{*}}\right)_{00}ds.}
\end{align*}
Therefore
\begin{equation}\label{eq:fourier_sum}
\begin{split}
    \mathbb{E}\left[{\bar{\mu}}(t,x)\overline{{\bar{\mu}}(t,x')}\right]
    &= \sum_{k}\mathbb{E}\left[{\hat{\bar{\mu}}}(t,k)\frac{e^{ikx}}{L}\overline{{\hat{\bar{\mu}}}(t,k)}\frac{e^{-ikx'}}{L}\right], \\
    &=\frac{1}{L^{2}}\sum_{k \neq 0}e^{ik(x-x')}\mathbb{E}\left[{c_0(t,k)\overline{c_0(t,k)}}\right].
\end{split}
\end{equation}

\begin{figure}[htb!]
\centering
\includegraphics[width=0.725\textwidth]{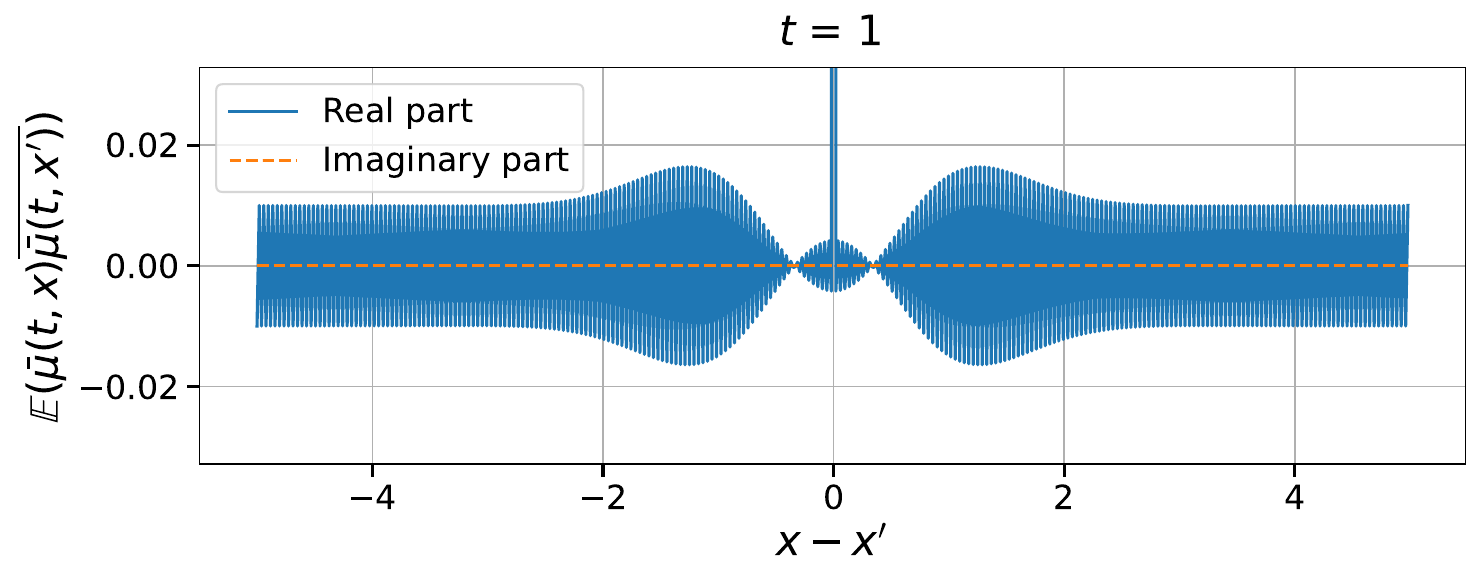}
\includegraphics[width=0.725\textwidth]{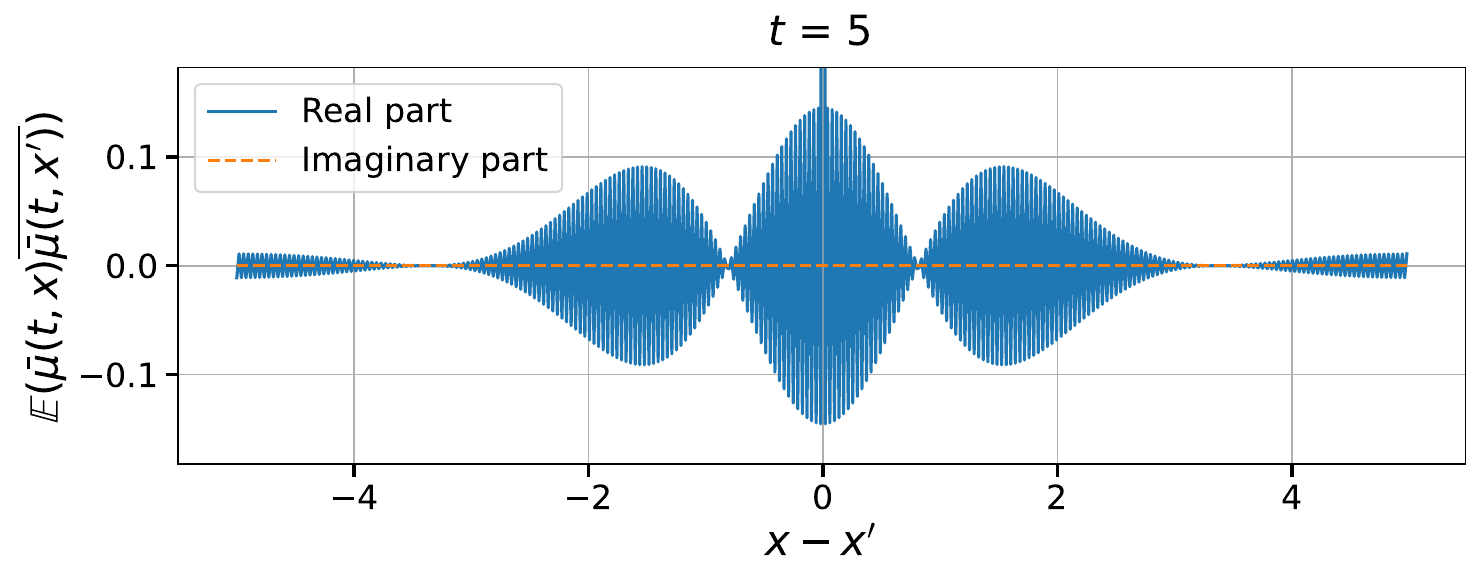} 
\includegraphics[width=0.725\textwidth]{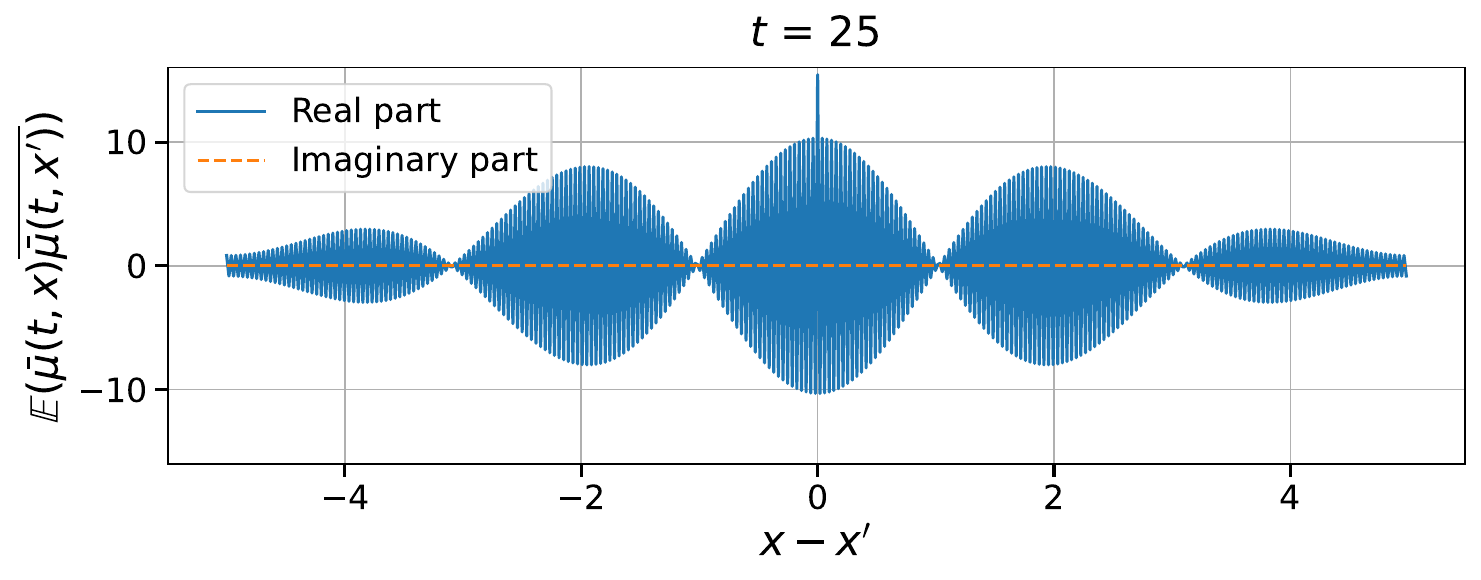} 
\includegraphics[width=0.725\textwidth]{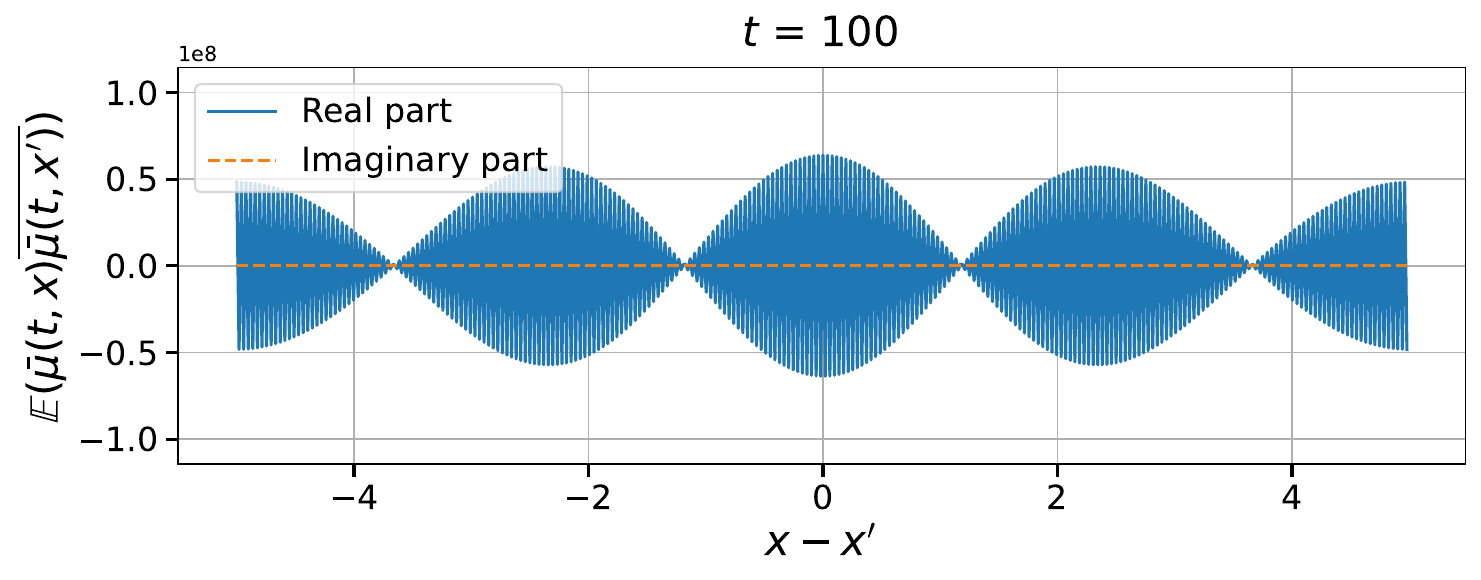} 
\caption{ ${\mathbb{E}[\bar{\mu}(t,x)\overline{\bar{\mu}(t,x')}]}$ from \eqref{eq:fourier_sum} for the Gaussian potential for $\sigma = 1/\sqrt{2}$, $L = 10$, $\gamma = 1.0$ and $\beta = 25 > \beta_{c}$.}
\label{fig:gaussian_fluc_snapshots}
\end{figure}

\begin{figure}[htb!]
\centering
\includegraphics[width=0.325\textwidth]{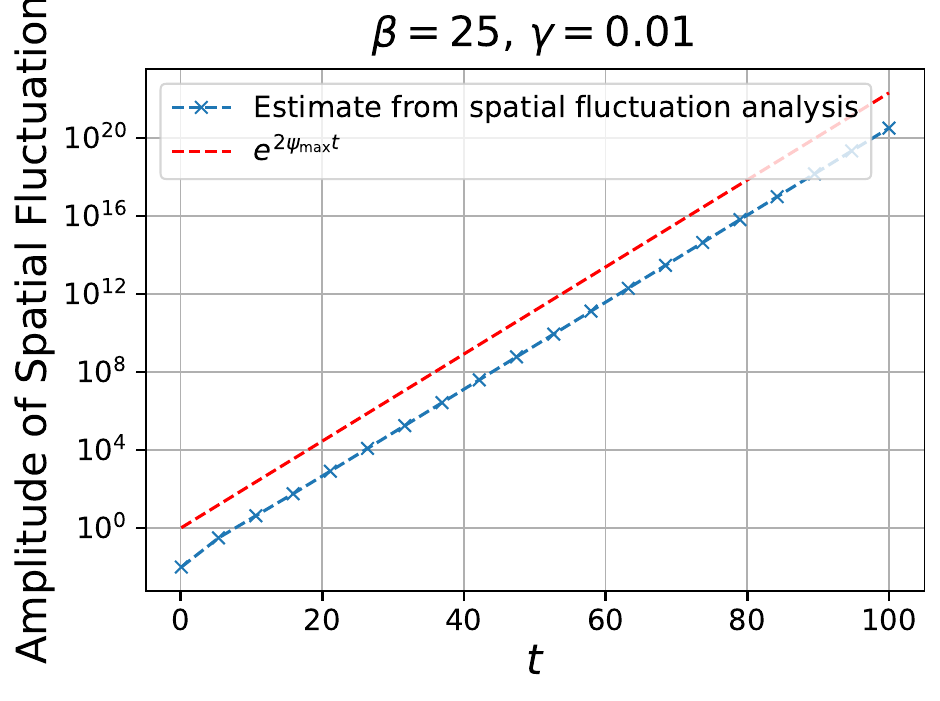}
\includegraphics[width=0.325\textwidth]{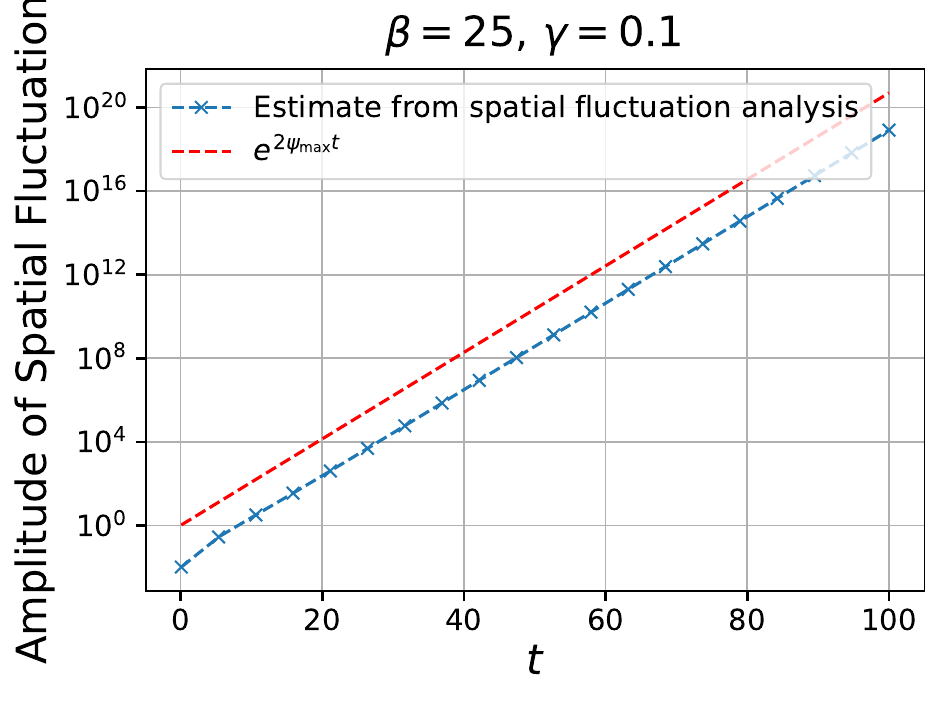} 
\includegraphics[width=0.325\textwidth]{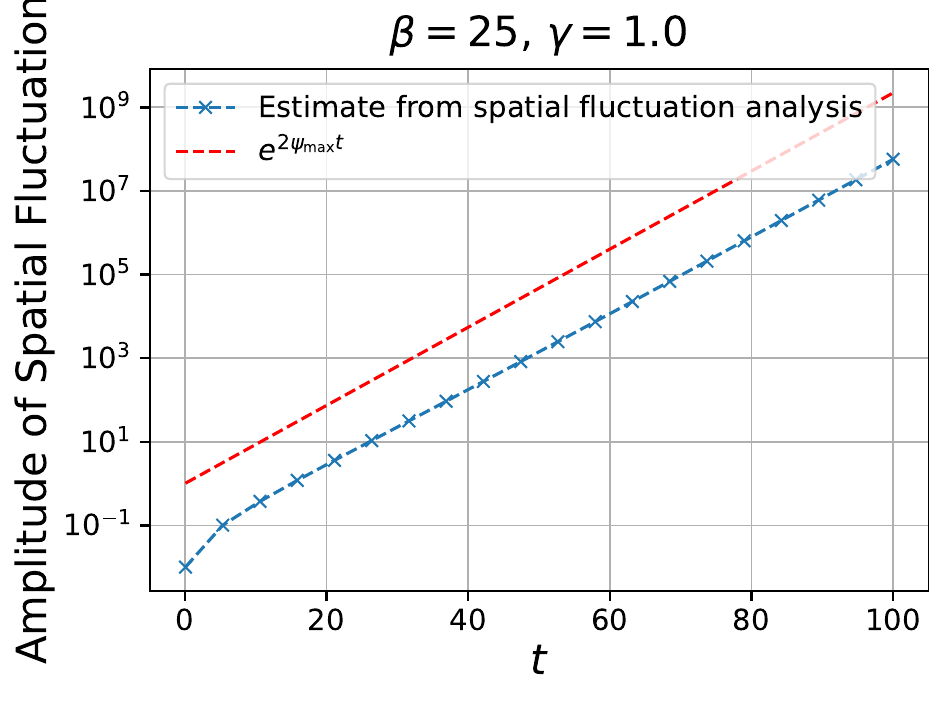} 
\caption{The amplitudes of the oscillations shown in Fig. \ref{fig:gaussian_fluc_snapshots} for varying friction values $\gamma$. These are estimated based on cropping out the Dirac-delta close to  $x=x'$. }
\label{fig:gaussian_amplitude_of_fluc}
\end{figure}

Compared to the overdamped setting (as considered in \cite{garnier2017consensus}), it is difficult to approximate \eqref{eq:fourier_sum} analytically. However we can accurately approximate this numerically, by first truncating the infinite sum and estimating each term {$\mathbb{E}\left[c_0(t,k)\overline{c_0(t,k)}\right]$} by truncating the infinite matrices $\mathbf{A}_{\gamma}(k,\beta)$ and $\mathbf{B}$ and using quadrature methods to estimate the integrals. Via this procedure for a Gaussian interaction potential, we observe in Fig. \ref{fig:gaussian_fluc_snapshots} the same qualitative behaviour as in the overdamped setting \cite{garnier2017consensus}. In particular, for small time the perturbation of the {spatial} density concentrates around $x= x'$ with oscillatory behaviour and with a Dirac-delta at $T = 0$. For large time, the fluctuations converge to a cosine wave centred at $x = x'$ with exponentially increasing amplitude in $T$. From Fig. \ref{fig:gaussian_fluc_snapshots}, we approximate the amplitude of the oscillatory component in the size of fluctuations based on cropping out close to the Dirac-delta at $x=x'$. In Fig. \ref{fig:gaussian_amplitude_of_fluc} we compare the amplitude of the oscillatory component with time, and we observe on a log-scale that the gradient coincides with {$2\psi_{\max}$} (from the linear stability analysis), hence the amplitude of fluctuation scales like {$e^{2\psi_{\max}t}$} in the long-time regime. Under this scaling, analogously to the overdamped setting, the linear system becomes unstable and the central limit theorem breaks down when {$\bar{\mu}(t,x)/\sqrt{N}$, the fluctuation of the spatial perturbation, is larger than $\bar{\rho}_{0}(x) = \frac{1}{L}$, the unperturbed spatial state}. In particular, the central limit theorem breaks down for a time $t_{\textnormal{cl}}>0$ where {$\mathbb{E}[\bar{\mu}(t_{\textnormal{cl}},x)\overline{\bar{\mu}(t_{\textnormal{cl}},x)}]/N \approx 1/L^{2}$, the variance of the initial spatial density $\bar{\rho}_{0}$}. Then, with the scaling observed via the numerical approximation to \eqref{eq:fourier_sum} ({$\mathbb{E}[\bar{\mu}(t_{\textnormal{cl}},x)\overline{\bar{\mu}(t_{\textnormal{cl}},x)}] \propto e^{2\psi_{\max}t}$} for the considered interaction potentials), we have
\begin{equation} \label{eq: t_cl_expression}
t_{\textnormal{cl}} \approx \frac{1}{2\psi_{\max}} \ln{N},
\end{equation}
for $N \gg 1$ sufficiently large.

\section{Numerical simulations}
\label{sec:num_analysis}

\subsection{Numerical integration of kinetic Langevin dynamics}\label{sec:numerical_integrators}
We want to compare the analytical findings of the previous sections with numerical simulations of the $N$-particle IPS, which evolves according to the kinetic (underdamped) Langevin dynamics \eqref{eq: Langevin Dynamics}. For this aim, we need to integrate \eqref{eq: Langevin Dynamics} using suitable numerical schemes. {Although conventional schemes like the standard Euler-Maruyama method can achieve accuracy by reducing the integration stepsize $h$, the corresponding computational workload scales with $h^{-1}$. For an $N$-particle system with pairwise interactions, the computational workload scales as $\mathcal{O}(N^2)$, meaning that scaling up to larger system sizes is only feasible if one can utilize large enough stepsizes. Maximizing the stepsize is also crucial because resolving the phase transition requires simulating long trajectories across thousands of independent realizations to accurately capture statistical averages, a challenge that is further compounded by slow convergence near metastable states. A more detailed discussion of these computational challenges and our stepsize selection criteria is provided in Appendices \ref{sup:sec:IPS_simulation} and \ref{sec: choice of stepsize}. While the standard Euler-Maruyama method introduces severe sampling bias or numerical instabilities at larger stepsizes, specialized splitting integrators remain highly accurate and stable. }  The main integrator we use is based on splitting the SDE \eqref{eq:Langevin} in the following way:
\[
\begin{pmatrix}
d\bm{x} \\
d\bm{v}
\end{pmatrix} = \underbrace{\begin{pmatrix}
0 \\
-\nabla U(\bm{x})dt
\end{pmatrix}}_{\mathcal{B}}+ \underbrace{\begin{pmatrix}
\bm{v}dt \\
0
\end{pmatrix}}_{\mathcal{A}} +\underbrace{\begin{pmatrix}
0 \\
-\gamma \bm{v} dt + \sqrt{2\gamma\beta^{-1}}d\mathbf{B}_t
\end{pmatrix}}_{\mathcal{O}},
\]
where the $\mathcal{B}$, $\mathcal{A}$ and $\mathcal{O}$ parts can be integrated exactly (in the distributional sense) over a time interval of size $h > 0$.
Specifically, define
the maps $\mathcal{A}:(\bm{x},\bm{v},h)\mapsto (\bm{X},\bm{V})$, $\mathcal{B}:(\bm{x},\bm{v},h)\mapsto (\bm{X},\bm{V})$ and
$\mathcal{O}:(\bm{x},\bm{v},h)\mapsto (\bm{X},\bm{V})$, by the formulas
\begin{eqnarray*}
\mathcal{A}(\bm{x},\bm{v},h) & = & (\bm{x}+h\bm{v}, \bm{v}),\\
\mathcal{B}(\bm{x},\bm{v},h) & = & (\bm{x}, \bm{v} -h\nabla U(\bm{x})),\\
\mathcal{O}(\bm{x},\bm{v},h) & = & (\bm{x},\eta_h\bm{v} + \sqrt{(1-\eta_h^2)\beta^{-1}} \bm{\xi}),
\end{eqnarray*}
where, in the last formula, $\eta_h=\exp(-h\gamma)$ and $\bm{\xi} \sim \mathcal{N}(0_d,I_d)$ with $0_d$ and $I_d$ the $d$-dimensional zero and identity matrix, respectively.  Compositions of these maps yield numerical methods with weak approximation properties. For example, the ABO method consists of  the map defined by
\[
\mathcal{O}(\mathcal{B}(\mathcal{A}(\bm{x},\bm{v},h),h),h).
\]
More complicated compositions yield higher order approximations in the weak sense. When a letter appears twice in the definition of an integrator, the assumption is that each part is integrated for a half-timestep. For example, the OBABO method consists of the following sequence of assignments:
\begin{eqnarray*}
(\bm{X},\bm{V}) & := &
\mathcal{O}(\bm{x},\bm{v},h/2),\\
(\bm{X},\bm{V}) & := &
\mathcal{B}(\bm{X},\bm{V},h/2),\\
(\bm{X},\bm{V}) & := &
\mathcal{A}(\bm{X},\bm{V},h),\\
(\bm{X},\bm{V}) & := &
\mathcal{B}(\bm{X},\bm{V},h/2),\\
(\bm{X},\bm{V}) & := &
\mathcal{O}(\bm{X},\bm{V},h/2),
\end{eqnarray*}
the key point being that each part of the vector field is ultimately integrated for a total step of length $h$.

Splitting methods of this type were studied extensively in \cite{LeMa2013,LeMaSt2016}, see also related work \cite{AbViZy2015}.  Among other things it has been shown that these methods have stationary distributions which approximate (in small parameter $h$) the Gibbs-Boltzmann canonical distribution.  Moreover, symmetric splittings yield even-order approximations of the invariant distribution and are thus of weak order at least two.  For details see the above-mentioned references. A particular observation of \cite{LeMa2013,LeMaSt2016} is that the BAOAB scheme, which is normally of weak order two, can achieve a higher order in the approximation of averages with respect to the invariant distribution in the case of configurational observables, specifically when the friction $\gamma$ is large.  This can be seen as a consequence of a fortuitous cancellation of terms in an error expansion derived based on the Baker-Campbell-Hausdorff formula applied to the products of the distributional propagators introduced by the splitting.

By contrast to weakly accurate schemes, so-called strongly accurate methods of order $p$ are pathwise accurate, we provide discussion of the UBU integrator which has improved strong order properties in Sec. \ref{sup:sec:alternative_integrators} of the appendix. We note that both the pathwise and the discretization error depend on the friction parameter $\gamma$. In our simulations, to maintain consistent accuracy across different values of $\gamma$, we choose the stepsize $h =\mathcal{O}\left(\min\left\{\gamma,1/\gamma\right\}\right)$ (see \cite[Assumption A1]{durmus2021uniform}, \cite[Theorem 10]{schuh2024} and \cite[Example 4.5]{LePaWh24} for an explicit computation in the case of a harmonic oscillator).  If one is interested in the pathwise accuracy of the dynamics it would be desirable to use the strongly accurate integrator. However, if the goal is finite-time weak accuracy, both integrators would be order two in the stepsize, and with the appropriate stepsize scaling will approximate the finite-time weak accuracy to the same level of accuracy. Finite-time weak accuracy is likely to be the more useful property in simulating the onset of cluster formation.

\subsection{Simulation and cluster detection} \label{sec:simulation} 
To simulate the IPS, we created \textbf{SimIPS}, an efficient and lightweight C++ code which can be accessed on our \href{https://github.com/SchroedingersLion/Cluster-Formation-in-Diffusive-Systems}{GitHub} page \cite{GitHub_IPS}. The simulation is initialized by placing $N$ particles in a periodic cubic box $[0,L]^d$, where $d\in\{1,2\}$. The initial positions are sampled uniformly. For a given inverse temperature $\beta$, the initial velocity components are sampled independently from their known Gaussian equilibrium distribution $\rho_{\beta,v}\propto e^{-\beta\frac{v^2}{2}}$. Given a stepsize $h$ and friction coefficient $\gamma$, we then use the BAOAB integrator to propagate the system, where various observables as well as configurational snapshots are printed every $n$ steps. For more details on the simulation code, see Sec. \ref{sup:sec:IPS_simulation}. {It also supports the three-dimensional case, $d=3$, and we provide an exemplary video animation in the GitHub repository}.

For $\beta>\beta_{\text{c}}${,} {the linearly stable stationary solutions of the McKean-Vlasov equation \eqref{eq: Fokker-Planck} consist of a family of single-cluster states, unique up to translation.} To numerically detect that state, we {evaluate} the mean particle distance to the system's centre of mass{. In an unbounded Euclidean domain, this is conceptually given by:} 
\begin{equation}\label{eq:dcom}
d_{\text{com}}(t)\coloneqq \frac{1}{N}\sum_{i=1}^{N} \bigg\|x_i(t) - \frac{1}{N}\sum_{j=1}^{N}x_{j}(t)\bigg\|.
\end{equation}
{Because our system resides on a flat torus, the standard arithmetic mean and Euclidean distance are ill-defined. In our simulations, $d_{\text{com}}$ is rigorously computed by substituting the minimum image convention for the distances and employing circular statistics for the centre of mass, following the method of \cite{BaiBreen}.} $d_{\text{com}}$ will be large in the uniform particle distribution and small in the one-cluster state.
We then define the convergence time by 
\begin{equation}\label{eq: d_com_criterion}
t^{*}_{s}\coloneqq \min \{t\ |\ d_{\text{com}}(t)<s \},
\end{equation}
where $s$ is a characteristic length of the interaction potential (e.g., $\sigma$ in the formulas for the Gaussian and GEM-$\alpha$ potentials). As mentioned in the previous section, in order to compare  $t^{*}_{s}$ to theory one needs to average over multiple independent trajectories. Its variance strongly depends on friction $\gamma$ which governs the coupling to the stochastic term in (\ref{eq: Langevin Dynamics}), see Sec. \ref{sup:sec:observable_plots_vs_gamma}.
\\
To measure the time for the onset of clustering, $t_{\text{cl}}$,  the criterion (\ref{eq: d_com_criterion})
is not useful due to the existence of intermediate metastable states of multiple clusters. Instead, we use the DBSCAN algorithm \cite{DBSCAN1} implemented in Scikit-learn \cite{scikit-learn} to detect clusters in the particle configurations. {DBSCAN identifies clusters as contiguous regions of high point density separated by areas of low-density noise. It relies on two hyperparameters: a search radius $\epsilon$ and a minimum number of points $N_{\min}$, which together define a local density threshold for a point to be considered part of a cluster core. Points satisfying this criterion are grouped based on their spatial connectivity, while isolated points are categorized as noise.} We give more details on DBSCAN and how we tuned its hyperparameters in appendix \ref{sup: sec: DBSCAN details}. We analyze a given trajectory frame by frame and detect the time of the first frame that shows a cluster. We also examined various ways in which the mean squared displacement (MSD) or the total potential energy in the system can be used to detect clusters, see Secs. \ref{sup:sec: MSD_details} and \ref{sup:sec:energy_vs_msd}. While more memory and compute efficient than DBSCAN, we found that their success depends more strongly on the system parameters (in particular the friction), and hence refrained from using them in the main experiments.

\subsection{Results} 
We perform extensive numerical simulations to study the behavior of the IPS and compare the results to our theoretical predictions from the previous sections. 
\subsubsection{One-dimensional results} \label{sec:results}
\begin{figure}[htb!]
\includegraphics[width=1.0\textwidth]{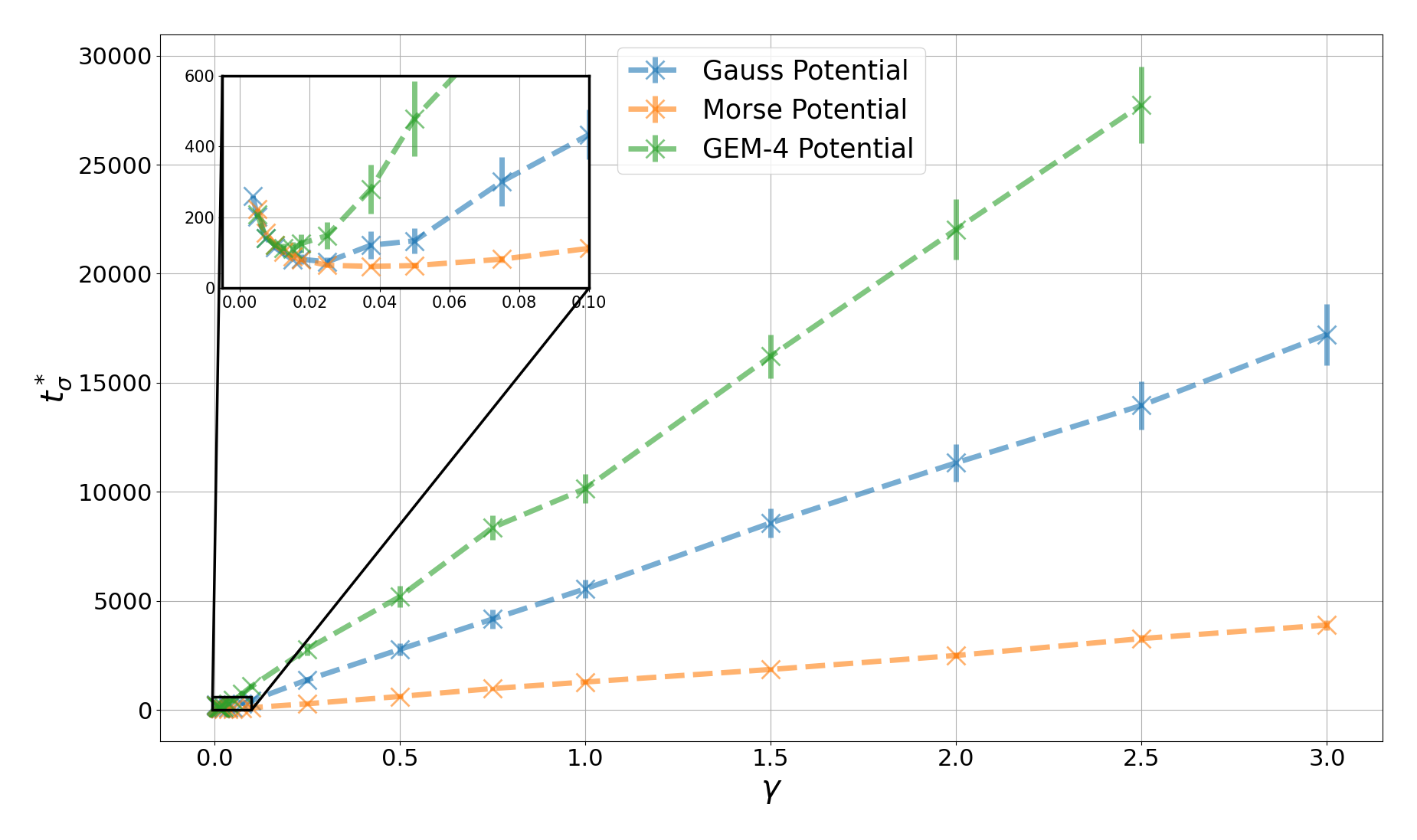}
\caption{\label{fig:convergence_times_1D} Simulation results for the friction-dependent convergence times $t^*_{\sigma}$ to reach the one-cluster equilibrium state when starting from a uniform initial distribution. The vertical bars denote 95\%-confidence intervals.}
\end{figure}

First, we measure the friction-dependent convergence times to reach the one-cluster equilibrium state for $\beta>\beta_c$ when starting from a uniform distribution. We perform experiments with 700 particles at $\beta=25$ for all three interaction potentials. The convergence times are measured with criterion \eqref{eq: d_com_criterion}, where we use $s=\sigma=\sqrt{0.5}$ ($\sigma$ being the width of the Gaussian and GEM potentials, see \eqref{eq:Potentials} or Sec. \ref{sup:sec:potentials} in the appendix).  For each potential and each friction $\gamma$, we average the measured times over a collection of independent trajectories.  Since larger $\gamma$ leads to stronger coupling of the dynamics to the noise term in \eqref{eq: Langevin Dynamics}, the convergence times are subject to higher variance at larger $\gamma$ (see also Sec. \ref{sup:sec:observable_plots_vs_gamma} in the appendix). Therefore, we average over larger numbers of independent trajectories at larger friction values to control the statistical error. We use a stepsize of $h=1$ at $\gamma=1$ and then scale it as $h\sim \min(\gamma, \gamma^{-1})$. For more details on the number of independent trajectories taken for the different runs, that number's influence on the result, and a discussion on how to pick the right stepsize, refer to Sections \ref{sec: choice of stepsize} and \ref{sec: trajectory_average}.  Fig. \ref{fig:convergence_times_1D} shows the results.
All potentials show the expected linear behaviour for intermediate and large frictions. The inset shows a local minimum and an increase in convergence times for $\gamma \to 0$. This limit is not covered by our linear stability analysis, but we provide some discussion in Sec. \ref{sup:sec:gamma=0}. We only note here that the dynamics for $\gamma=0$ no longer samples the canonical distribution. Instead, it reduces to Hamiltonian dynamics, which samples a surface of constant total energy, $H(\boldsymbol{x},\boldsymbol{v})=H(\boldsymbol{x}_0,\boldsymbol{v}_0)=\text{const.}$, where the subscript $0$ denotes the initial conditions.

\begin{figure}[htb!]
\includegraphics[width=1.0\textwidth]{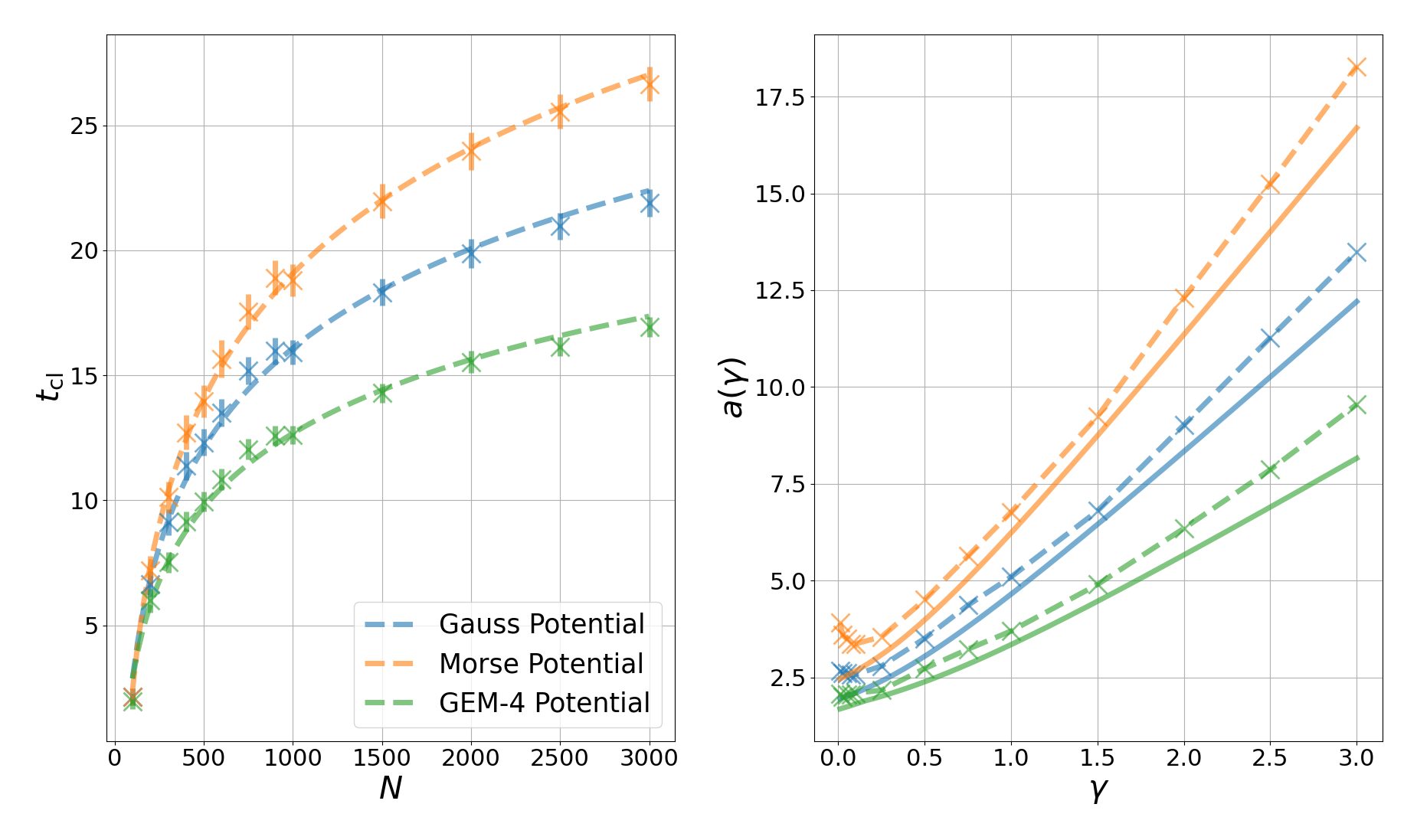}
\caption{\label{fig:prefactor_double_plot_1D} \textbf{Left:} Time to the onset of clustering $t_{\text{cl}}$ as obtained by the simulation at $\gamma=1$. The vertical bars show 95\% confidence intervals as obtained by 200 independent trajectories per particle count $N$. The dashed lines denote fits for $t_{\text{cl}}=a(\gamma, {\beta}) \log(N)+b(\gamma, {\beta})$. \textbf{Right:} $a(\gamma, {\beta})$ as obtained from the fits of the simulation data (dashed lines) and the linear stability analysis (solid lines).}
\end{figure}

For the times $t_{\text{cl}}$ to the onset of clustering, we fix $\beta=25$ as before and run 200 trajectories per friction value $\gamma$ and particle count $N$. We print out the trajectories and use the clustering algorithm DBSCAN to detect the first frame in a trajectory that shows a cluster-like structure. DBSCAN's hyperparameters are set to $\epsilon=0.5$ and $N_{\min}=N_{\min,0}\frac{N}{N_0}$ with $N_{\min,0}=90$ and $N_0=500$, i.e., $N_{\min}$ scales linearly with $N$. \\
Fig. \ref{fig:prefactor_double_plot_1D} (left) shows results for fixed $\gamma$ and varying $N$, ranging from $N=100$ to 3,000. We observe the logarithmic behavior in alignment with the theoretical result \eqref{eq: t_cl_expression}. To compute the prefactors in \eqref{eq: t_cl_expression}, $a(\gamma, {\beta}):=\frac{1}{2\psi_{\max}(\gamma, {\beta})}$, we only consider particle counts $N \geq 600$ for each $\gamma $, because the fluctuation analysis in Sec. \ref{sec:fluc_analysis} considers large particle counts $N$. The obtained results for $a(\gamma, {\beta})$ in Fig. \ref{fig:prefactor_double_plot_1D} (right) show qualitative agreement with the linear stability analysis across all examined potentials except for the smallest $\gamma$. For $\gamma \to 0$, they slightly begin to increase again, consistent with the convergence times $t^{*}_{\sigma}$ in Fig. \ref{fig:convergence_times_1D} (note that the $\gamma$-range in the $a(\gamma, {\beta})$ plots starts at larger values than in the $t^{*}_{\sigma}$ plot).

When looking at the confidence intervals for $t^{*}_{\sigma}$ and $t_{\text{cl}}$ in Figs. \ref{fig:convergence_times_1D} and \ref{fig:prefactor_double_plot_1D}, respectively, we observe that the intervals for $t^{*}_{\sigma}$ are much narrower for the Morse potential than for the other potentials (despite using far fewer trajectories to average over; see Sec. \ref{sec: trajectory_average}). This effect is absent for $t_{\text{cl}}$. We hypothesize that this is because the Morse potential practically has the longest {(non-negligible) tail} (see Fig. \ref{fig:potentials} in Sec. \ref{sup:sec:potentials}). In the intermediate multi-cluster phase, each cluster moves like an independent random walker until it gets close enough to another cluster for the attractive force to become non-negligible, which leads to the two clusters merging \cite{garnier2017consensus}. This process carries most of the variance in the convergence times. If the interaction potential has a longer effective range, it will reduce the average time for two clusters to meet and merge.\footnote{We refer to the sample trajectory animations for the Morse and the Gauss potential in 2D on our \href{https://github.com/SchroedingersLion/Cluster-Formation-in-Diffusive-Systems/tree/main/Animations}{GitHub} page to observe that intermediate clusters start to feel attractive pull at slightly longer distances for the Morse potential compared to the Gauss potential, leading to faster cluster merging.} This also explains why there is no substantial difference in the variances of $t_{\text{cl}}$ across the potentials, as this time measures the onset of clustering. 

\subsubsection{Two-dimensional results.}\label{sec: 2D results} 
\begin{figure}[!htbp]
\includegraphics[width=1.0\textwidth]{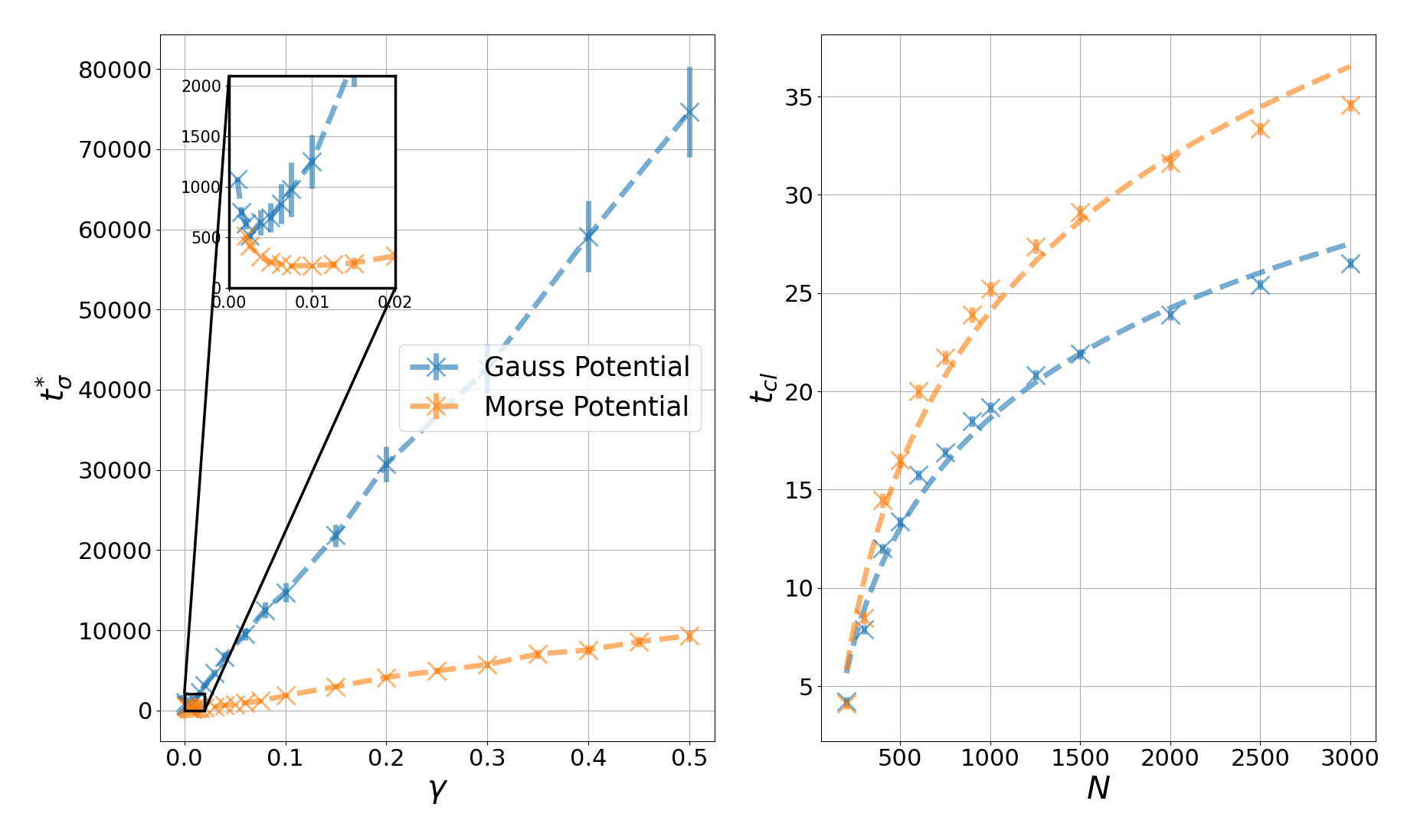}
\caption{\label{fig:results_2D} Simulation results for the two-dimensional IPS. \textbf{Left:} Friction-dependent convergence times $t^*_{\sigma}$ to reach the one-cluster equilibrium state when starting from a uniform initial distribution. \textbf{Right:} Time to the onset of clustering $t_{\text{cl}}$ as obtained by the simulation at $\gamma=0.1$. The dashed lines denote fits for $t_{\text{cl}}=a(\gamma, {\beta}) \log(N)+b(\gamma, {\beta})$. In both plots, the vertical bars denote 95\%-confidence intervals.}
\end{figure}

To provide evidence that our results generalize to higher dimensions, we simulate the IPS in two dimensions at $\beta=150$ for both the Gaussian and the Morse potential. Just like in the 1D case, the $\beta$ value is roughly 5 times larger than the critical $\beta_c$, where the latter was estimated with the formula derived in \cite{martzel2001mean} and via experimentation (see Sec. \ref{sec: critical_temperature}). For fixed particle number $N\!=\!1000$, we measure the friction-dependent times to reach the one-cluster equilibrium state by averaging over a friction-dependent number of independent trajectories, see Table \ref{tab:avg_scheme_2D} in \ref{sup:sec:2D_results_table}. We also obtain the times to the onset of clustering, $t_{\text{cl}}$, at fixed $\gamma=0.1$ and average over 400 independent trajectories (two times as many as in the 1D-case). Fig. \ref{fig:results_2D} shows the results. Comparing with Figs. \ref{fig:convergence_times_1D} and \ref{fig:prefactor_double_plot_1D}, we observe qualitatively similar behavior as in the 1D case: a linear increase of $t^*$ with friction $\gamma$ with a local minimum for $\gamma \approx 0$, and a logarithmic growth of $t_{\text{cl}}$ with particle count $N$.

\subsubsection{Numerical study of the critical temperature.} \label{sec: critical_temperature}
\begin{figure}[!htbp]
\includegraphics[width=1.0\textwidth]{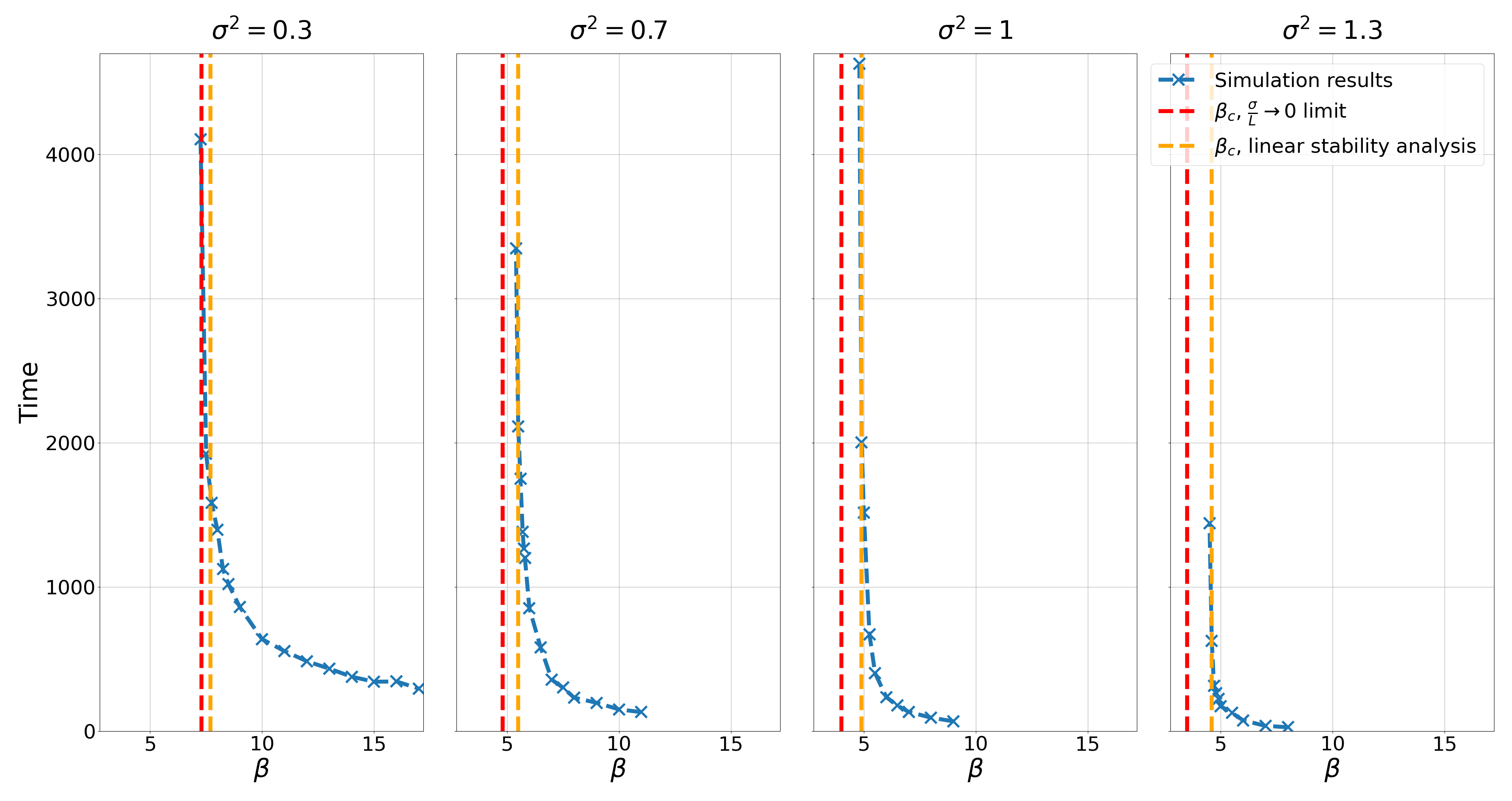}
\caption{\label{fig:critical_temperature} Simulation results for the cluster formation time against (inverse) temperature $\beta$ for an IPS in 1D using Gaussian interaction potentials of width $\sigma$. The red and yellow dashed lines denote the values of $\beta_c$ obtained from the literature (see text) and our linear stability analysis, respectively. Each run was done at $\gamma=0.01$, $1000$ particles, and box length $L=10$.}
\end{figure}

We test the validity of the linear stability analysis in Sec. \ref{sec: linear_stability} by calculating the critical temperature via simulation of the IPS and then comparing it with our analytical prediction. Note that the value of $\beta_c$  naturally depends on the characteristic length of the pairwise interaction potential, e.g., $\sigma$ in the case of a Gaussian kernel. To be more precise, it depends on $\frac {\sigma} {L}$ for torus length $L$, and we consider different values of this fraction here. To experimentally find the critical temperature at a given $\frac {\sigma} {L}$ via particle simulation, we run one-dimensional IPS simulations of 1,000 particles for different $\beta$ values, starting from a $\beta$ way above the analytically estimated $\beta_c$, and then successively decrease $\beta$. For each tested $\beta$, we measure the time until the one-cluster equilibrium state has been reached, defined as the quantity $t^{*}_{s}$ from \eqref{eq: d_com_criterion}, using {a numerical threshold $s=1.43\sigma^2$, where $\sigma^2$ denotes the numerical value of the variance parameter. The factor 1.43 is an empirical choice;} we do not merely use $s=\sigma$ as in our experiments from Fig. \ref{fig:convergence_times_1D}, because the widths of the equilibrium clusters tend to increase for larger temperatures (closer to $\beta_c$). In each setting, we average over 5 trajectories each\footnote{That number of independent trajectories is sufficient as we use a small friction, $\gamma=0.01$, leading to small variances in the convergence times.}. As $\beta$ approaches $\beta_c$ from above, the measured times become larger. The point where they diverge yields $\beta_c$. This is illustrated in Fig. \ref{fig:critical_temperature} for different interaction ranges $\sigma^2$.
The red dashed line in the figure denotes the results for $\beta_c$ derived in \cite{martzel2001mean} for the overdamped limit and $\frac{\sigma}{L}\to 0$, given by 
\begin{equation}
\beta_c = \frac{L^d}{(2\pi)^{\frac{d}{2}}\sigma^d}.
\end{equation}
The agreement of their results with ours (yellow dashed lines) increases for decreasing $\sigma^2$, as expected. 
The reason why the convergence time increases for $\beta \searrow \beta_c$ is because decreasing $\beta$ increases the temperature and hence the thermal energy or noise in the system. This in turn leads to clusters of larger widths and to longer clustering times as more particles are needed to combine their attractive forces to overcome the noise.
The $\beta$-dependent times to equilibrium and the corresponding cluster widths can also be seen in Fig. \ref{sup:fig:cluster_width} in Sec. \ref{sup:sec: temp_dependent_cluster_widths} of the appendix.

\subsubsection{Convergence via metastable states or `vacuum cleaner' dynamics.} \label{sec:trajectory_differences}
\begin{figure}[!htbp]
\includegraphics[width=1.0\textwidth]{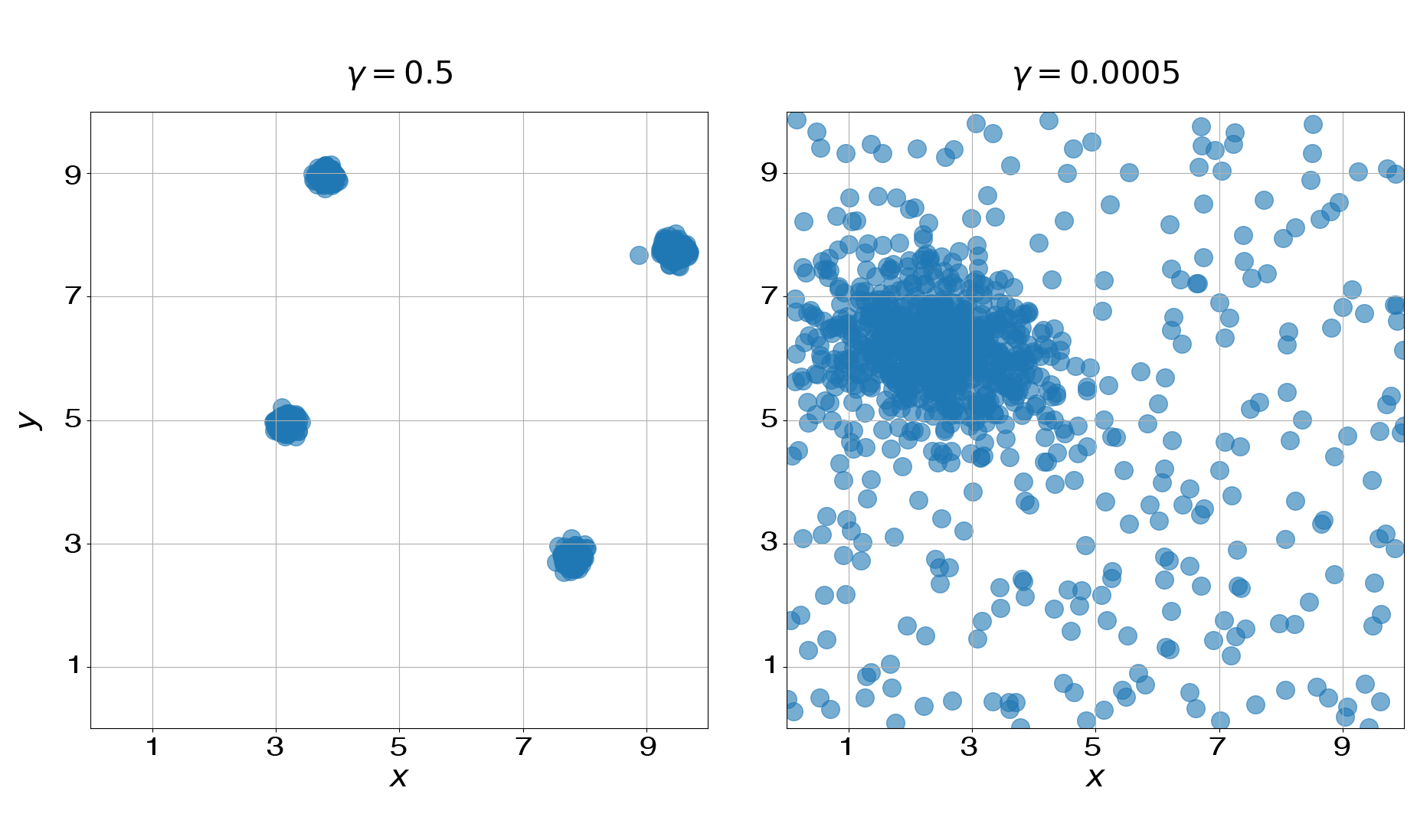}
\caption{\label{fig:metastable_vs_vacuum} Snapshots of a 2D IPS at $\beta>\beta_c$ for two different friction values $\gamma$. The snapshots were taken at the same time with all simulation parameters except $\gamma$ being the same. At large $\gamma$, the system enters metastable multi-cluster states. At small $\gamma$, one observes a single cluster in coexistence with a sea of unbounded particles, which will be successively collected (`vacuumed up') by the cluster.}
\end{figure}

The convergence to the one-cluster equilibrium state for $\beta>\beta_c$ shows qualitatively different behavior for small and large frictions, see Fig. \ref{fig:observbles_vs_gamma} in the appendix. At higher frictions, both the centre-of-mass distance $d_{\text{com}}$ and the mean squared displacement $\delta r$ admit plateaus when the system is in a metastable multi-cluster state. The convergence to the one-cluster state then happens via pronounced merging events, expressing itself as step-like patterns in the two observables. For smaller frictions, one observes a smoother convergence to the one-cluster state with no visible energy plateaus. This implies a slow but steady formation of a single cluster which keeps collecting particles until the steady state has been reached. Fig. \ref{fig:metastable_vs_vacuum} shows exemplary snapshots at small and high frictions at the same time after simulation start. 
In contrast to the high-friction trajectory in which multiple clusters form, the low-friction trajectory forms only a single large cluster with a substantial fraction of the system remaining in an unbound state. We label this type of convergence to equilibrium \textbf{`vacuum cleaner'} dynamics, as the single cluster that has formed slowly `vacuums up' all other particles. The reason why that happens at small frictions is that the dynamics becomes more Hamiltonian in that regime, so particle energies are locally better conserved and over longer timescales. 
For Hamiltonian dynamics, $\gamma=0$, the stationary state is known to be a co-existence phase of a single cluster surrounded by a gas of unbound particles (see Sec. \ref{sup:sec:gamma=0}). `Vacuum cleaner' dynamics at small frictions thus strikes a balance between the convergence to the one-cluster state at $\gamma>0$ and resemblance to the stationary state at $\gamma=0$.  However, it not only occurs at small frictions, but also at high temperatures sufficiently close to the critical temperature, i.e., $\beta \gtrsim \beta_c$. In this regime, it occurs at high frictions as well, but the underlying explanation is different. For large temperatures, rather than the better conservation of particle energies, the strong thermal noise is the reason for only a single cluster forming. At large frictions, the time for this to happen can be substantial, such that the uniform distribution remains stable for long times. An animation of `vacuum cleaner' dynamics at high temperature and with a prolonged uniform phase can be downloaded from \href{https://github.com/SchroedingersLion/Cluster-Formation-in-Diffusive-Systems/tree/main/Animations}{GitHub}.

\FloatBarrier

\subsubsection{Hamiltonian limit.}\label{sup:sec:gamma=0}
\begin{figure}[!htbp]
\includegraphics[width=0.85\textwidth]{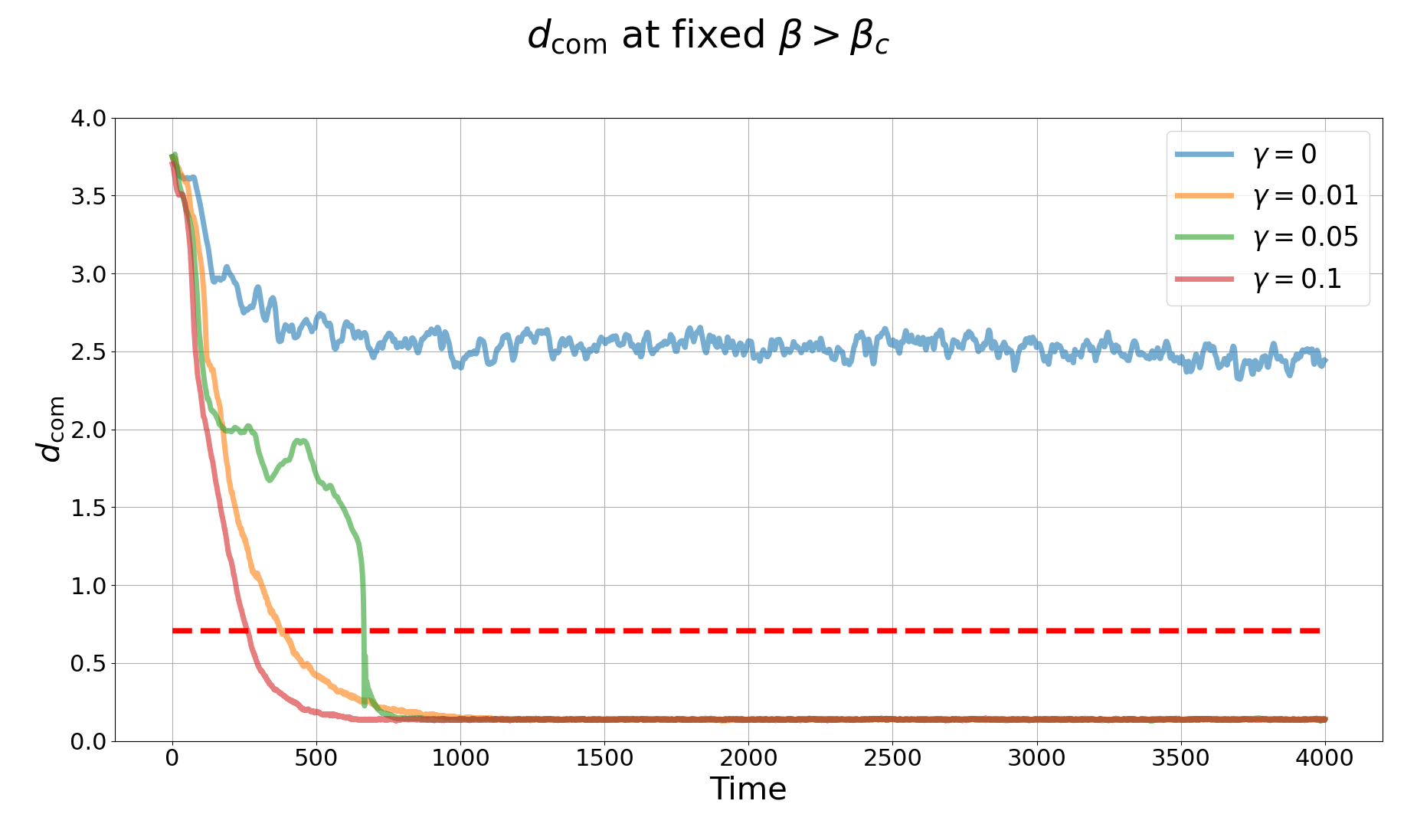}
\caption{\label{fig:cluster_size_vs_gamma} $d_{\text{com}}$ at different friction values $\gamma$ at fixed $\beta=1.5\beta_c$. 2D IPS, Gaussian interaction potential. Dashed line denotes threshold $s$ from (\ref{eq: d_com_criterion}).  }
\end{figure}

In the limit $\gamma\to0$, Langevin dynamics \eqref{eq: Langevin Dynamics} reduces to Hamiltonian dynamics, characterized by the conservation of the Hamiltonian, $\frac{d}{dt}H(\boldsymbol{x}(t), \boldsymbol{v}(t))=0$. The Hamiltonian is the total energy in the system given by the sum of potential energy $U(\boldsymbol{x})$ and kinetic energy $K(\boldsymbol{v})$,
\begin{equation}
H(\boldsymbol{x},\boldsymbol{v})=U(\boldsymbol{x})+K(\boldsymbol{v}).
\end{equation}
In our case, $U(\boldsymbol{x})=\frac{1}{2N}\sum_{i\neq j}W(\boldsymbol{x}_i,\boldsymbol{x}_j)$ and $K(\boldsymbol{v})=\sum_{i=1}^{N}\frac{1}{2}\|\boldsymbol{v}_i\|^2$. In a physical picture, conservation of the Hamiltonian thus means conservation of total energy. This is in contrast to Langevin dynamics for $\gamma>0$, in which the friction models the interaction of the system with a heat bath at constant temperature $T=\beta^{-1}$, exchanging energy with the bath such that the temperature is held constant (apart from statistical fluctuations). For more properties of Hamiltonian dynamics, we refer to \cite{LeimkuhlerReichBook}.
For $\gamma=0$, the dynamics is no longer ergodic with respect to $\rho_{\beta}$, so it is not obvious whether the phase transition persists in that limit. \\
To examine this numerically, we simulate a two-dimensional IPS at fixed $\beta>\beta_c$ and compare the evolution of the mean distance to centre of mass, $d_{\text{com}}$, for various friction parameters $\gamma$. The results can be seen in Fig. \ref{fig:cluster_size_vs_gamma}.

\begin{figure}[!htbp]
\includegraphics[width=1.0\textwidth]{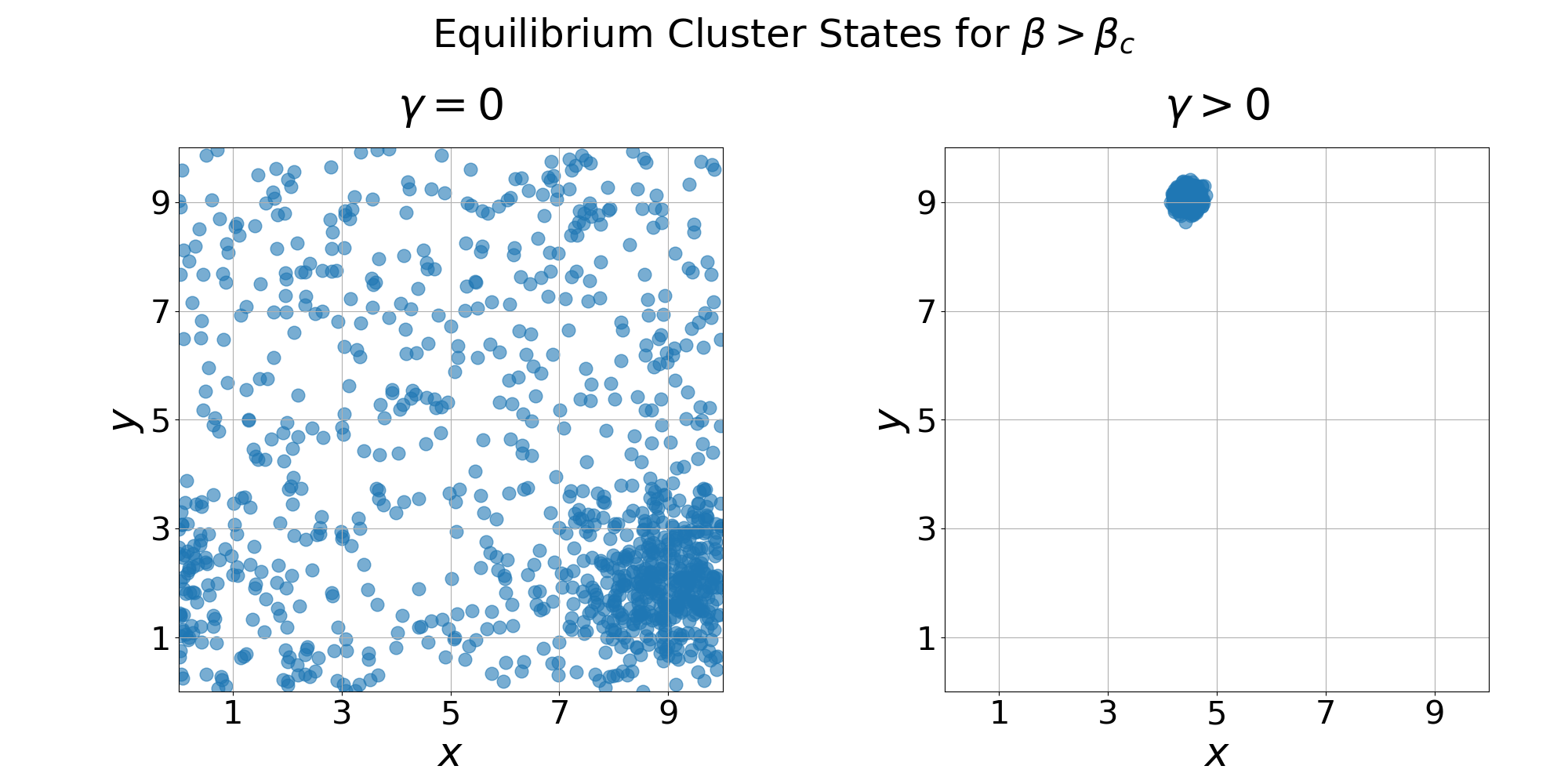}
\caption{\label{fig:final_cluster_state_snapshots} Snapshots of equilibrium configuration at $\beta>\beta_c$ for Gaussian interactions. \textbf{Left:} Hamiltonian dynamics, $\gamma=0$. \textbf{Right:} Langevin dynamics, $\gamma>0$. The snapshots were taken from the curves seen in Fig. \ref{fig:cluster_size_vs_gamma}, at time 1,500. }
\end{figure}

We observe that all trajectories for $\gamma>0$ equilibrate to the same value, well below the threshold that detects the final one-cluster state. This is in line with the fact that $\gamma$ influences the convergence speed to equilibrium, not the equilibrium distribution. The $\gamma=0$ curve reaches an equilibrium as well, but one with much greater $d_{\text{com}}$. The fact that it equilibrates to an intermediate value means that the corresponding state is indeed a semi-bounded state, but it might not be a state in which all particles are part of a cluster. For $\gamma>0$, plateaus of $d_{\text{com}}$ above threshold $s$ would denote intermediate multi-cluster states in which all particles are bound to well-separated clusters, which would then merge to reach the final one-cluster equilibrium (see Fig. \ref{fig:observbles_vs_gamma} in Sec. \ref{sup:sec:observable_plots_vs_gamma}). To confirm that the plateau for $\gamma=0$ is of a different nature, we plot a snapshot of the configuration in Fig. \ref{fig:final_cluster_state_snapshots} together with a snapshot of one of the $\gamma>0$ curves at the same point in time.
The $\gamma=0$ snapshot clearly shows a one-cluster state, but one in which not all particles are bound to the cluster, very different from the usual one-cluster equilibrium state and also from the metastable multi-cluster states. The latter are usually accompanied by a plateau in $d_{\text{com}}$ with $d_{\text{com}}>s$ as well. The particle distribution looks comparable to the one observed during `vacuum cleaner' dynamics seen in Fig. \ref{fig:metastable_vs_vacuum} of Sec. \ref{sec:trajectory_differences}. However, in the latter case all particles slowly but steadily join the cluster, which is accompanied by a steady decline in $d_{\text{com}}$. In contrast, the configuration in Fig. \ref{fig:final_cluster_state_snapshots} is stable. In Sec. \ref{sup:sec:ham_limit} of the appendix we provide additional discussion and investigation of the affect of temperature on this phenomena and on the widths of the forming cluster. Additionally, we provide animations of the $\gamma=0$ Hamiltonian trajectory from the figure on our \href{https://github.com/SchroedingersLion/Cluster-Formation-in-Diffusive-Systems/tree/main/Animations}{GitHub} page.

We note that the observations made here were already discussed as far back as in 1990 in \cite{posch1990dynamics} and the references therein, where it was shown that for Hamiltonian dynamics with the Gaussian interaction potential, the equilibrium state is either a uniform phase or a co-existence phase of a single cluster surrounded by a cloud of free particles. The cluster formation is only possible for small enough total energies (as governed by the initial conditions), and the fraction of trapped particles increases for decreasing energies. Both of these points are consistent with our results in Fig. \ref{fig:final_cluster_state_snapshots} and Fig. \ref{fig:cluster_size_vs_beta} in the appendix. While we did not directly measure the fraction of particles in the cluster, the average particle-to-centre-of-mass distance $d_{\text{com}}$ we plot in Fig. \ref{fig:cluster_size_vs_beta} serves as an equivalent measure. The interested reader may also compare the video animation provided by \cite{posch1990dynamics} (\href{https://www.mediathek.at/atom/018AAB81-08A-0279F-00000484-0189A3E5}{Link}) with the one we provide in our GitHub repository. An additional interesting point discussed in \cite{posch1990dynamics} is that, while the cluster forms, the kinetic energy of the cluster particles and hence the cluster's local temperature is larger than that of the surrounding particle cloud, where temperature alignment only happens in equilibrium.  

Finally, we remark that when measuring the time to full convergence, Fig. \ref{fig:convergence_times_1D} in Sec. \ref{sec:results}, we observe a divergence of the convergence time for $\gamma \to 0$. This divergence is consistent with the observations made here, because the criterion used to detect the final one-cluster state, eq. \eqref{eq: d_com_criterion}, assumes that all particles are bound within one cluster, and the equilibrium state for $\gamma=0$ does not satisfy this condition.

\section{Conclusion}\label{sec:conclusions}

In this paper, we studied cluster formation in interacting kinetic Langevin dynamics. Using linear stability analysis and via an exploration of the fluctuations around the mean-field limit, we established the onset of clustering and the breakdown of the central limit theorem for inverse temperatures $\beta > \beta_c$, {i.e., in the regime where the uniform state of the mean-field limit of the IPS becomes linearly unstable.} To perform the linear stability analysis, we developed an efficient spectral numerical method for solving the linearised kinetic McKean–Vlasov equation. This enabled the study of several attractive interaction potentials. {We also employed extensive numerical simulations to study the Langevin IPS in the weak-friction regime, finding excellent agreement of the experimentally observed critical temperature with the value predicted by our linear stability analysis}. The critical temperature was found to be independent of the friction parameter, while the time scale for cluster formation scales approximately linearly with it.

The present study leaves several open questions that we plan to address in future work. First, although we focused on attractive potentials, our methodology is readily extendable to more general classes of interactions that include repulsive components (such as those typical of molecular simulations \cite{AllenTildesley2017}) where careful initialization may be required. A key open problem is to obtain precise analytic bounds on the real part of the spectrum of the tridiagonal operator arising in our linear stability analysis, which could yield explicit estimates of the critical temperature and cluster formation times. In addition, our simulations revealed intriguing emergent behaviours, such as the ``vacuum cleaner'' dynamics, which warrant further study. In~\cite{gerber2025formationclusterscoarseningweakly} a rigorous analysis of cluster formation was presented for the overdamped dynamics, i.e., in the limit as $\gamma \rightarrow +\infty$, and the link between the dynamics of the clusters and the Massive Arratia flow~\cite{Konarovskyi2017} was established. It would be interesting to obtain similar rigorous results for the kinetic Langevin dynamics, in particular in the weakly underdamped regime. Finally, while we studied the stability of the uniform phase with respect to changing temperatures, one may also examine the stability of the one-cluster state. In \cite{martzel2001mean}, it was already observed for the one-dimensional Gaussian potential and in the high-friction limit that the critical temperature for the transition from clustered phase to uniform phase is larger than for the reverse transition. While we did observe this hysteresis effect in our experiments (not reported here), we leave its rigorous study for future works.

\vspace{1cm}
\paragraph{Acknowledgments}
GP is partially supported by an ERC-EPSRC Frontier Research Guarantee through Grant No. EP/X038645, ERC Advanced Grant No. 247031 and a Leverhulme Trust Senior Research Fellowship, SRF$\backslash$R1$\backslash$241055. \\
BL, RL and PAW acknowledge the support of the Engineering and Physical Sciences Research Council Grant
EP/S023291/1 (MAC-MIGS Centre for Doctoral Training).

\begin{appendix}

\section{Linear Stability Analysis} 
\subsection{The linearized kinetic McKean-Vlasov PDE}\label{sup:sec:linear_stability1}
To perform the linear stability analysis we start with the kinetic Fokker-Planck equation, (\ref{eq: Fokker-Planck}), which is given by
\begin{equation*}
\begin{split}
\frac{\partial \rho}{\partial t} &= -\bm{v}\cdot \nabla_{\bm{x}}\rho +  \nabla_{\bm{v}}\rho \cdot \left[\int_{\mathcal{D}_{\bm{v}'}} \int_{\mathcal{D}_{\bm{y}}} \rho(t,\bm{x}-\bm{y},\bm{v}')\nabla W(\bm{y})d\bm{y}d\bm{v}'\right] + \gamma \nabla_{\bm{v}} \cdot \left(\bm{v}\rho\right) + \gamma \beta^{-1}\Delta_{\bm{v}}\rho,
\end{split}
\end{equation*}
where the integration domains are given by $\mathcal{D}_{\bm{v}'}:=(-\infty, \infty)^d$ and $\mathcal{D}_{\bm{y}}:=\left[-\frac{L}{2}, \frac{L}{2}\right)^d$. Note that despite particle positions being in the domain $\mathbb{T}^d_L=[0,L)^d$, the integration bounds for $\bm{y}$ need to be given by $\mathcal{D}_{\bm{y}}$, since $\bm{y}$ plays the role of a distance between two positions, which due to the minimum image convention lies in $\left[-\frac{L}{2}, \frac{L}{2} \right]^d$.
We insert the ansatz $\frac{1}{L^d}F(\bm{v})+\rho_1(t,\bm{x},  \bm{v})$ for $\rho$ and obtain (using the fact that $F$ is neither time nor space dependent)
\begin{equation*}
\begin{split}
\frac{\partial \rho_1}{\partial t} &= -\bm{v}\cdot \nabla_{\bm{x}}\rho_1 +  \nabla_{\bm{v}}\left(\frac{1}{L^d}F+\rho_1\right) \cdot \left[\int_{\mathcal{D}_{\bm{v}'}} \int_{\mathcal{D}_{\bm{y}}} \left(\frac{1}{L^d}F(\bm{v}')+\rho_1(t,\bm{x}-\bm{y},\bm{v}')\right)\nabla W(\bm{y})d\bm{y}d\bm{v}'\right] \\ &+ \gamma \nabla_{\bm{v}} \cdot \left[\bm{v}\left(\frac{1}{L^d}F+\rho_1\right)\right] + \gamma \beta^{-1}\Delta_{\bm{v}}\left(\frac{1}{L^d}F+\rho_1\right).
\end{split}
\end{equation*}
We can neglect the $\mathcal{O}(\rho_1^2)$ term and use $\int \nabla W(\bm{y})d\bm{y}=\bm{0}$ since $\nabla W$ is odd. We obtain
\begin{equation*}
\begin{split}
\frac{\partial \rho_1}{\partial t} &= -\bm{v}\cdot \nabla_{\bm{x}}\rho_1 +  \frac{1}{L^d}\nabla_{\bm{v}}F \cdot \left[\int_{\mathcal{D}_{\bm{v}'}} \int_{\mathcal{D}_{\bm{y}}} \rho_1(t,\bm{x}-\bm{y},\bm{v}')\nabla W(\bm{y})d\bm{y}d\bm{v}'\right] \\
& \quad + \frac{1}{L^d} \gamma\left[ \nabla_{\bm{v}} \cdot \left(\bm{v}F \right) + \beta^{-1} \Delta_{\bm{v}} F \right] \\
& \quad + \gamma \nabla_{\bm{v}} \cdot \left( \bm{v}\rho_1 \right) + \gamma \beta^{-1} \Delta_{\bm{v}} \rho_1.
\end{split}
\end{equation*}
Since $F$ is stationary for the Ornstein-Uhlenbeck process, we have 
$$ \nabla_{\bm{v}} \cdot \left(\bm{v}F \right) +  \beta^{-1} \Delta_{\bm{v}} F=0,$$ and arrive at
\begin{equation}
\begin{split}
\frac{\partial \rho_1}{\partial t} &= -\bm{v}\cdot \nabla_{\bm{x}}\rho_1 +  \frac{1}{L^d}\nabla_{\bm{v}}F \cdot \left[\int_{\mathcal{D}_{\bm{v}'}} \int_{\mathcal{D}_{\bm{y}}} \rho_1(t,\bm{x}-\bm{y},\bm{v}')\nabla W(\bm{y})d\bm{y}d\bm{v}'\right] \\
& \quad + \gamma \nabla_{\bm{v}} \cdot \left( \bm{v}\rho_1 \right) + \gamma \beta^{-1} \Delta_{\bm{v}} \rho_1,
\end{split}
\end{equation}
which is \eqref{eq: Fokker-Planck-linear} in the main text.

\subsection{Fourier expansion}\label{sup:sec:linear_stability2}
The spatial spectral components of $\rho_1$ are given by
$$\hat{\rho}_1(t,\bm{k},\bm{v}) = \int_{\mathcal{D_{\bm{x}}}}\rho_1(t, \bm{x}, \bm{v}) e^{-i\bm{k}\cdot \bm{x}} d\bm{x}, $$ with domain $\mathcal{D_{\bm{x}}}=\left[0,L\right]^d$ and $\bm{k}\in\left\{\left(\frac{2\pi n_1}{L}, ...\frac{2\pi n_d}{L}\right)\middle| n_1,...,n_d \in \mathbb{Z} \right\}$. Taking the time-derivative and using (\ref{eq: Fokker-Planck-linear}) we obtain 
\begin{equation*}
\begin{split}
\frac{\partial \hat{\rho_1}}{\partial t} &= -\bm{v}\cdot \int_{\mathcal{D}_{\bm{x}}}\nabla_{\bm{x}}\rho_1 e^{-i\bm{k}\cdot \bm{x}}d\bm{x} +  \frac{1}{L^{d}}\nabla_{\bm{v}}F \cdot \left[\int_{\mathcal{D}_{\bm{x}}}\int_{\mathcal{D}_{\bm{v}'}} \int_{\mathcal{D}_{\bm{y}}} \rho_1(t,\bm{x}-\bm{y},\bm{v'})\nabla W(\bm{y})d\bm{y}d\bm{v'}e^{-i\bm{k}\cdot \bm{x}}d\bm{x}\right] \\
&\quad + \gamma \nabla_{\bm{v}} \cdot \left( \bm{v}\int_{\mathcal{D}_{\bm{x}}}\rho_1 e^{-i\bm{k}\cdot \bm{x}}d\bm{x} \right) + \gamma \beta^{-1} \Delta_{\bm{v}} \int_{\mathcal{D}_{\bm{x}}}\rho_1 e^{-i\bm{k}\cdot \bm{x}}d\bm{x}, \\
&= -i\bm{k}\cdot\bm{v} \hat{\rho}_1 +  \frac{1}{L^{d}}\nabla_{\bm{v}}F \cdot \left[\int_{\mathcal{D}_{\bm{v}'}}\int_{\mathcal{D}_{\bm{y}}} \int_{\mathcal{D}_{\bm{x}}} \rho_1(t,\bm{x}-\bm{y},\bm{v'})e^{-i\bm{k}\cdot(\bm{x}-\bm{y})}d\bm{x}\nabla W(\bm{y})e^{-i\bm{k}\cdot \bm{y}}d\bm{y}d\bm{v'}\right] \\
&\quad + \gamma \nabla_{\bm{v}} \cdot \left(\hat{\rho}_1 \bm{v} \right) + \gamma \beta^{-1} \Delta_{\bm{v}} \hat{\rho}_1, \\
&= -i\bm{k}\cdot\bm{v} \hat{\rho}_1 +  \frac{1}{L^d}\nabla_{\bm{v}}F \cdot \left[\int_{\mathcal{D}_{\bm{v}'}} \hat{\rho}_1d\bm{v'}\int_{\mathcal{D}_{\bm{y}}}\nabla W(\bm{y})e^{-i\bm{k}\cdot \bm{y}}d\bm{y}\right] + \gamma \nabla_{\bm{v}} \cdot \left(\hat{\rho}_1 \bm{v} \right) + \gamma \beta^{-1} \Delta_{\bm{v}} \hat{\rho}_1, \\
&=  -i\bm{k}\cdot\left[\bm{v} \hat{\rho}_1 - {\frac{1}{L^{d}}}\nabla_{\bm{v}}F \hat{W} \int_{\mathcal{D}_{\bm{v}'}} \hat{\rho}_1d\bm{v'}\right] + \gamma \nabla_{\bm{v}} \cdot \left(\hat{\rho}_1 \bm{v} \right) + \gamma \beta^{-1} \Delta_{\bm{v}} \hat{\rho}_1,
\end{split}
\end{equation*}
with $\hat{W}(\bm{k}):=\int_{\mathcal{D}_{\bm{y}}} W(\bm{y})e^{-i\bm{k}\cdot \bm{y}}d\bm{y}$ the Fourier transform of $W(\bm{y})$, where $\mathcal{D}_{\bm{y}}=\left[-\frac{L}{2}, \frac{L}{2} \right]^d$. 

\subsection{Hermite polynomials}\label{sup:sec: Hermite_polynomials}
The normalized Hermite polynomials are defined as
\begin{equation}\label{eq:Hermite_functions}
h_n(v)=\frac{1}{\sqrt{n!}}H_n\Big(\sqrt{\beta} v \Big),
\end{equation}
where
\begin{equation}\label{eq:Hermite_polynomials}
H_n(v) = (-1)^n e^{\frac{v^2}{2}}\frac{\text{d}^n}{\text{d}v^n} \Big( e^{-\frac{v^2}{2}} \Big).
\end{equation}
The Hermite polynomials satisfy the recursion relation
\begin{equation}\label{eq:Hermite_recursion}
H_{n+1}(u)=u H_n(u)-n H_{n-1}(u),
\end{equation}
which we make use of in Sec. \ref{sup:sec:linear_stability3}.
The normilized Hermite polynomials $h_n$ are the eigenfunctions of the generator of the Ornstein-Uhlenbeck process~\cite[Sec. 4.4]{PavliotisStochBook}
\begin{equation*}
-\mathcal{L}_{OU}h_n = n h_n, \quad n=0,1,2,...
\end{equation*}
where 
\begin{equation*}
\mathcal{L}_{OU} = -v \frac{\text{d}}{\text{d}v} + \beta^{-1}\frac{\text{d}^2}{\text{d}v^2}.
\end{equation*}
In addition, they form a complete orthonormal basis of the space $L^2(F)$ with $F(v)\text{d}v = \frac{1}{\sqrt{2\pi \beta^{-1}}} e^{-\beta\frac{v^2}{2}} \text{d}v$. In particular, \textbf{a)} they satisfy an orthonormality condition
\begin{equation}\label{eq: Hermite_orhonormality}
\langle h_n, h_m \rangle_{L^2(F)} = \int_{-\infty}^{\infty}h_n(v)h_m(v)F(v)dv=\delta_{n,m},
\end{equation}
and \textbf{b)} any function $g \in L^2(F)$ can be expanded as 
\begin{equation*}
g(v) = \sum_{n=0}^{\infty}c_n h_n(v).
\end{equation*}
For more details on Hermite polynomials in this context, we refer to \cite[Sec.~4.4]{PavliotisStochBook}. 

\subsection{Hermite integral calculation} \label{sup:sec:linear_stability3}
For the linear stability analysis, we need to calculate integrals of the form
\begin{equation*}
\begin{split}
I:=\int_{-\infty}^{\infty}v h_n(v) h_m(v) F(v) dv &= \frac{1}{\sqrt{n!m!}} \int_{-\infty}^{\infty} v H_n(\sqrt{\beta}v) H_m(\sqrt{\beta}v)F(v)dv, \\
&= \frac{1}{\beta \sqrt{n!m!}} \int_{-\infty}^{\infty} u H_n(u)H_m(u)F\left(\frac{u}{\sqrt{\beta}}\right)du.
\end{split}
\end{equation*}
We use the recursion relation for Hermite polynomials \eqref{eq:Hermite_recursion} to obtain
\begin{equation*}
\begin{split}
I &= \frac{1}{\beta \sqrt{n!m!}}\left\{ \int_{-\infty}^{\infty} H_{n+1}(u)H_m(u)F\left(\frac{u}{\sqrt{\beta}}\right)du + n \int_{-\infty}^{\infty} H_{n-1}(u)H_m(u) F\left(\frac{u}{\sqrt{\beta}}\right) du \right\}, \\
&= \frac{1}{\sqrt{2\pi n! m! \beta}} \left\{ \int_{-\infty}^{\infty} H_{n+1}(u)H_m(u)e^{-\frac{u^2}{2}}du + n \int_{-\infty}^{\infty} H_{n-1}(u)H_m(u) e^{-\frac{u^2}{2}} du \right\}, \\
&= \sqrt{\frac{m!}{n!\beta}}\Big\{ \delta_{m,n+1} + (m+1)\delta_{m, n-1} \Big\},
\end{split}
\end{equation*} 
where in the last step we used the orthogonality condition $\int_{-\infty}^{\infty} H_m(u)H_n(u) e^{-\frac{u^2}{2}} du=\sqrt{2\pi}m!\delta_{m,n}$, which follows trivially from the relations \eqref{eq:Hermite_functions} and \eqref{eq: Hermite_orhonormality}.

\section{Additional Simulations and Details}
\subsection{Alternative numerical integrators of kinetic Langevin dynamics.}\label{sup:sec:alternative_integrators}

In Sec. \ref{sec:numerical_integrators} we focus on weakly accurate numerical integrators. By contrast to weakly accurate schemes, so-called {\em strongly accurate} methods of order $p$ are pathwise accurate, meaning that when
a trajectory is defined by iteration of a certain numerical method 
\[
(\bm{x}_{n+1},\bm{v}_{n+1}) = \psi_h(\bm{x}_{n},\bm{v}_{n}),
\]
the solution at a fixed time $T = Nh$ satisfies
\[
\E \|(\bm{x}_N,\bm{v}_N)-(\bm{x}(Nh),\bm{v}(Nh))\| = O(h^p),
\]
where both the discrete and continuous solutions are, importantly, based on the same Wiener noise path.  The splitting methods defined above, even when they have higher weak order, are strong order one.

An alternative splitting first introduced in \cite{skeel1999integration} requires only one gradient evaluation per step and is strong order two (see also \cite{sanz2021wasserstein,BUBthesis,chada2023unbiased}). The alternative splitting of the SDE \eqref{eq:Langevin} is based on the following components
\[
\begin{pmatrix}
d\bm{x} \\
d\bm{v}
\end{pmatrix} = \underbrace{\begin{pmatrix}
0 \\
-\nabla U(\bm{x})dt
\end{pmatrix}}_{\mathcal{B}} +\underbrace{\begin{pmatrix}
\bm{v}dt \\
-\gamma \bm{v} dt + \sqrt{2\gamma\beta^{-1}}d\mathbf{B}_t
\end{pmatrix}}_{\mathcal{U}},
\]
which can be integrated in the weak sense exactly over an interval of size $h>0$. With $\eta = \exp{\left(-\gamma h\right)}$, the operators corresponding to these maps are given by
\begin{equation}\label{eq:Bdef}
\mathcal{B}(\bm{x},\bm{v},h) = (\bm{x},\bm{v} - h\nabla U(\bm{x})),
\end{equation}
and
\begin{equation}\label{eq:Udef}
\begin{split}
\mathcal{U}(\bm{x},\bm{v},h,\bm{\xi}^{(1)},\bm{\xi}^{(2)}) &= \Big(\bm{x} + \frac{1-\exp{\left(-\gamma h\right)}}{\gamma}\bm{v} + \sqrt{\frac{2}{\gamma\beta}}\left(\mathcal{Z}^{(1)}\left(h,\bm{\xi}^{(1)}\right) - \mathcal{Z}^{(2)}\left(h,\bm{\xi}^{(1)},\bm{\xi}^{(2)}\right) \right),\\
& \bm{v}\exp{\left(-\gamma h\right)} + \sqrt{2\gamma\beta^{-1}}\mathcal{Z}^{(2)}\left(h,\bm{\xi}^{(1)},\bm{\xi}^{(2)}\right)\Big),
\end{split}
\end{equation}

where 
\begin{equation}\label{eq:Z12def}
\begin{split}
\mathcal{Z}^{(1)}\left(h,\bm{\xi}^{(1)}\right) &= \sqrt{h}\bm{\xi}^{(1)},\\
\mathcal{Z}^{(2)}\left(h,\bm{\xi}^{(1)},\bm{\xi}^{(2)}\right) &= \sqrt{\frac{1-\eta^{2}}{2\gamma}}\Bigg(\sqrt{\frac{1-\eta}{1+\eta}\cdot \frac{2}{\gamma h}}\bm{\xi}^{(1)} + \sqrt{1-\frac{1-\eta}{1+\eta}\cdot\frac{2}{\gamma h}}\bm{\xi}^{(2)}\Bigg),
\end{split}
\end{equation}
and $\bm{\xi}^{(1)}, \bm{\xi}^{(2)} \sim \mathcal{N}\left(0_{d},I_{d}\right)$ are independent standard normal random variables. Using these maps, the UBU integration scheme with stepsize $h >0 $ is defined by 
\begin{equation}\label{eq:PhUBU}
\begin{split}
\left(\bm{x}_{k+1},\bm{v}_{k+1}\right)&=\mathcal{UBU}\left(\bm{x}_k,\bm{v}_k,h,\bm{\xi}^{(1)}_{k+1},\bm{\xi}^{(2)}_{k+1},\bm{\xi}^{(3)}_{k+1},\bm{\xi}^{(4)}_{k+1}\right)
\\
&=\mathcal{U}\left(\mathcal{B}\left(\mathcal{U}\left(\bm{x}_{k},\bm{v}_{k},h/2,\bm{\xi}^{(1)}_{k+1},\bm{\xi}^{(2)}_{k+1}\right),h\right),h/2,\bm{\xi}^{(3)}_{k+1},\bm{\xi}^{(4)}_{k+1}\right),
\end{split}
\end{equation}
where $\xi^{(i)}_{k+1} \sim \mathcal{N}(0_{d},I_{d})$ for all $i = 1,...,4$ and $k \in \mathbb{N}$. While UBU has strong order 2, both integrators, BAOAB and UBU, have the same order of accuracy in terms of stepsize $h>0$ in their respective invariant measures.

\subsection{Observable time series for different temperatures and frictions.}\label{sup:sec:observable_plots_vs_gamma}
\begin{figure}[!htbp]
\includegraphics[width=1.0\textwidth]{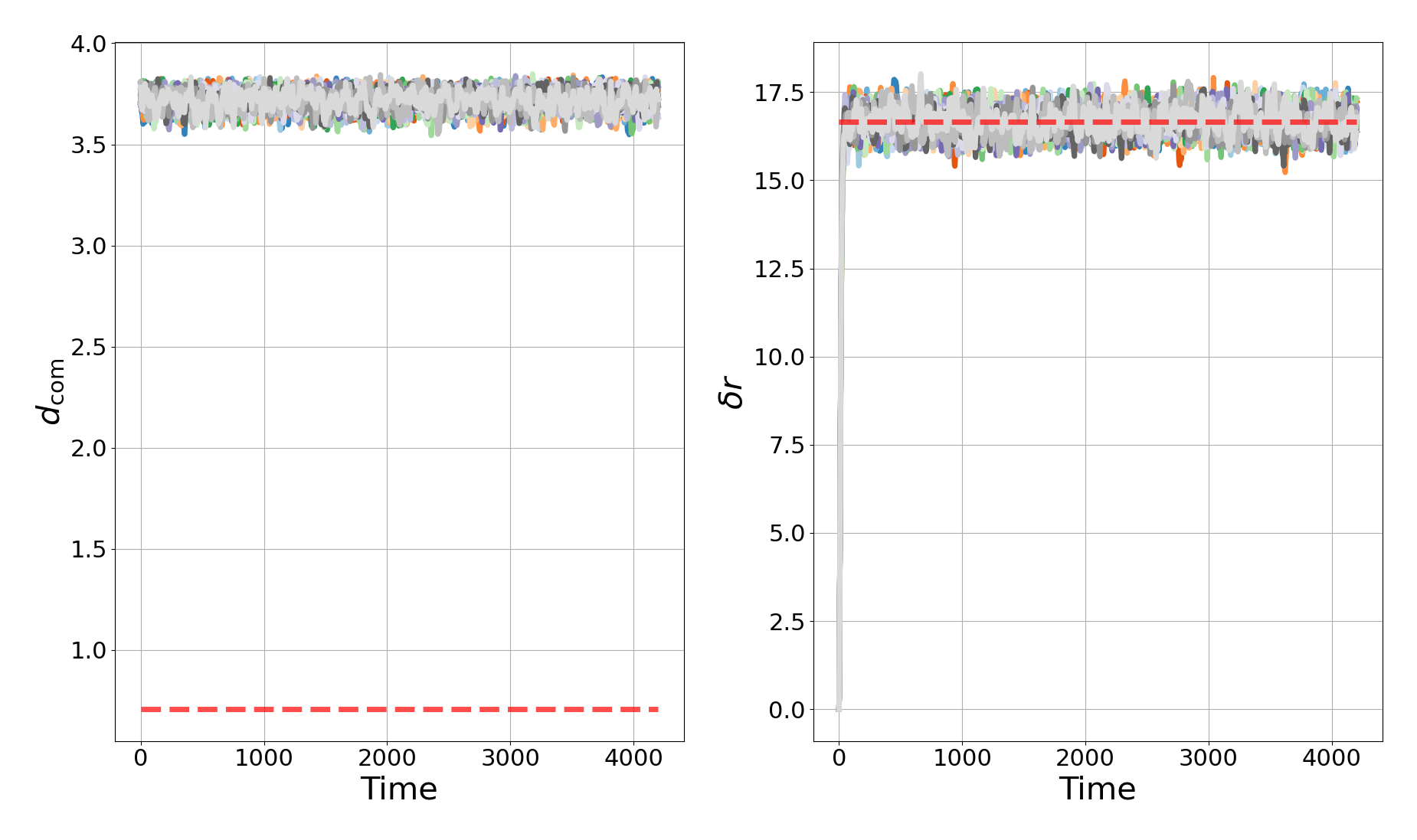}
\caption{\label{fig:observables_highT} Simulation results for a two-dimensional IPS at $\beta<\beta_c$. \textbf{Left:} Mean distance to centre of mass $d_{\text{com}}$. Dashed line denotes $\sigma$ from (\ref{eq: d_com_criterion}).  \textbf{Right:} Mean squared displacement $\delta r$. Dashed line denotes equilibrium expectation. Each curve corresponds to an independent trajectory. }
\end{figure}

\begin{figure}[!htbp]
\includegraphics[width=1.0\textwidth]{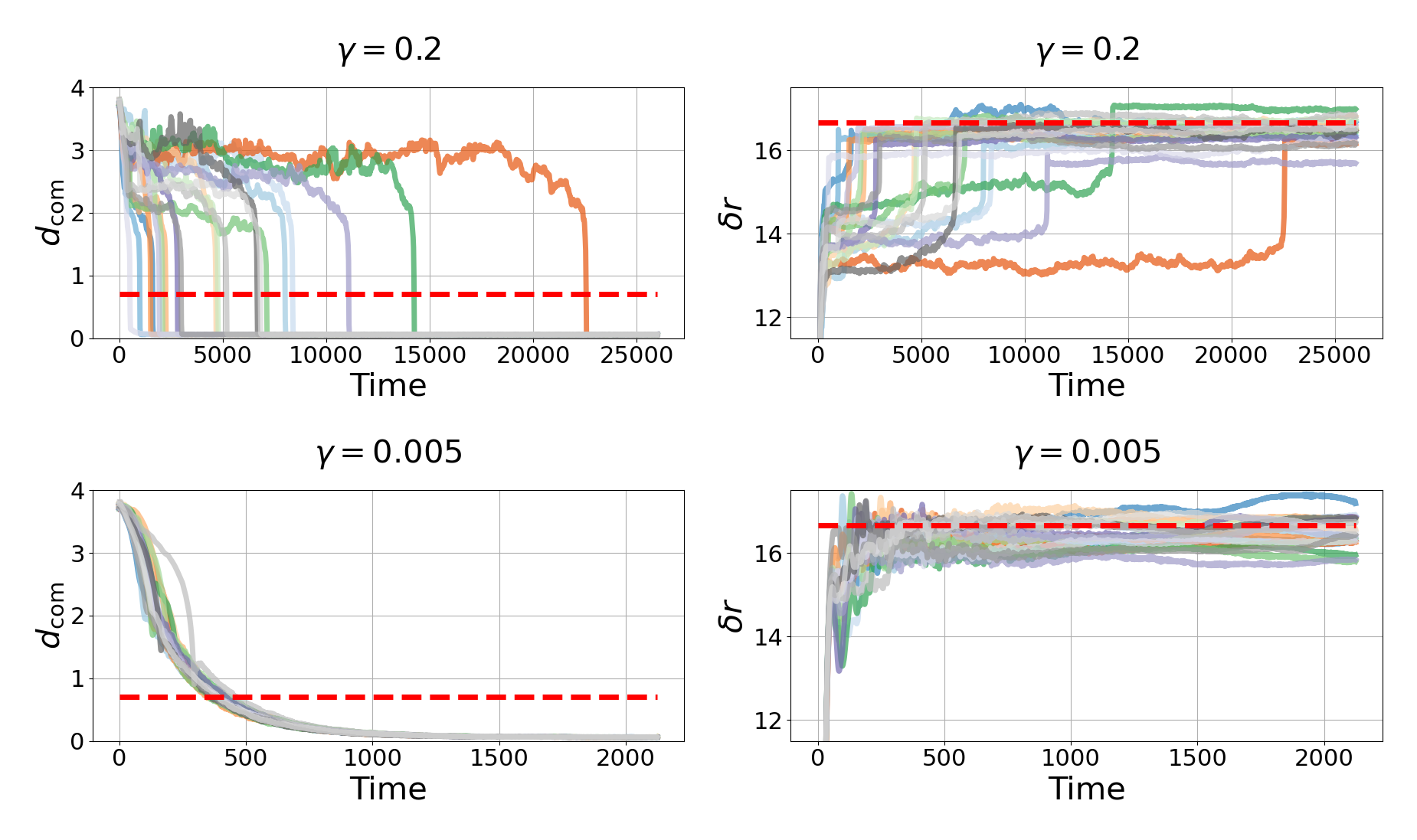}
\caption{\label{fig:observbles_vs_gamma} Simulation results for a two-dimensional IPS at large friction $\gamma$ (\textbf{top}) and small friction (\textbf{bottom}). \textbf{Left:} Mean distance to centre of mass $d_{\text{com}}$. Dashed line denotes $\sigma$ from (\ref{eq: d_com_criterion}).  \textbf{Right:} Mean squared displacement $\delta r$. Dashed line denotes equilibrium expectation. Each curve corresponds to an independent trajectory. }
\end{figure}

We show example plots for the two observables, the mean particle distance $d_{\text{com}}$ to the system's centre of mass from (\ref{eq:dcom}) and the mean squared displacement $\delta r$  (see \eqref{eq: MSD} in Sec. \ref{sup:sec: MSD_details}). We simulate a two-dimensional IPS consisting of $N=1,000$ particles interacting via the Morse potential.
First, we show results for $\beta<\beta_c$, i.e., the gaseous phase, in Fig. \ref{fig:observables_highT}. 
Both observables oscillate around their equilibrium values, with the difference that $d_{\text{com}}$ already starts on that level whereas the $\delta r$ has to converge to there from 0. The red dashed lines denote the threshold value $\sigma$ from criterion (\ref{eq: d_com_criterion}) in the $d_{\text{com}}$ plot and the equilibrium value for $\delta r$ (which is independent of temperature, see Sec. \ref{sup:sec: MSD_details}). We note that the $d_{\text{com}}$ curves remain far away from the threshold $\sigma$  to detect the phase of a single cluster. \\
To demonstrate that these observables can detect cluster formation dynamics, we simulate the same system at $\beta > \beta_c$ for two different friction values $\gamma$ in Fig \ref{fig:observbles_vs_gamma}.  

Unlike the case at high temperature ($\beta < \beta_c$), the $d_{\text{com}}$ curves decay below the threshold, and from the fact that the curves never increase again from there it becomes clear that criterion \eqref{eq: d_com_criterion}, $d_{\text{com}}<\sigma$, is indeed a reliable tool to detect the one-cluster equilibrium state.
Note also that at the larger friction the dynamics of the observables shows more variance. In particular, the times at which $d_{\text{com}}$ falls below the threshold is subject to higher variance, explaining the wider confidence intervals in Fig. \ref{fig:convergence_times_1D} for larger $\gamma$ and the resulting necessity to average over more trajectories (as described in Sec. \ref{sec: trajectory_average}). \\
In the plots for $\delta r$, the plateaus that form for some trajectories below the equilibrium line indicate metastable multi-cluster states. These plateaus are visible in $d_{\text{com}}$ as well, but we often found them to be less pronounced than in $\delta r$.

\FloatBarrier
\subsection{Details on the IPS simulation}\label{sup:sec:IPS_simulation}
To simulate the interacting particle system (IPS), we created \textbf{SimIPS}, an efficient and lightweight C++ code which can be accessed on our \href{https://github.com/SchroedingersLion/Cluster-Formation-in-Diffusive-Systems}{GitHub} page. It can be launched from the command line, where all simulation hyperparameters, including what potential to use, are specifiable via command-line flags. It places $N$ particles in a periodic cubic box $[0,L]^d$, where $d\in\{1,2,3\}$. The initial positions are either sampled uniformly or read from a passed file. For a given inverse temperature $\beta$, the initial velocity components are sampled independently from their known Gaussian equilibrium distribution $\rho_{\beta,v}\propto e^{-\beta\frac{v^2}{2}}$. Given a stepsize $h$, inverse temperature $\beta$, friction coefficient $\gamma$, and desired number of iterations $N_{\text{iter}}$, one of two Langevin integrators are used (either BAOAB or UBU, see Sec. \ref{sec:numerical_integrators}) to propagate the dynamics. Time series data of a collection of observables is printed to an output file by default, but one can also toggle the printing of the entire trajectory. By default, the code collects three scalar observables that are related to the phase of the system (clustered vs. disordered) and that can be used to gauge the stepsize-dependent stability of the simulation: The average particle distance to the centre of mass, $d_{\text{com}}$, the mean squared displacement $\delta r$, and the kinetic temperature $T_{\text{kin}}$. Optionally, the whole trajectory may be printed to a separate file. The GitHub page also offers two Python scripts to visualize the data from the printed files: one to plot the three observables, and one to create a video animation.\footnote{In one dimension, this will just be another static plot.}\\
The code is object-oriented and modular, allowing for straightforward extensions in terms of integrators, interaction potentials, and observables to collect. It is also general in the sense that any $N$-particle system may be treated as long as the interaction function is pairwise additive and radial-symmetric. \\
For more information on how to use and potentially customize it, see the Readme on the GitHub page.\\
\\
From a computational point of view, simulating interacting particle systems is oftentimes challenging. For the particular systems we consider, this is because of three reasons. \textbf{1)} The processes we need to observe, specifically the merging of two clusters, are rare events. In dimensions $d > 1$, this is particularly pronounced. Thus, one needs to be able to run long enough trajectories for the system to actually reach the one-cluster equilibrium state (given $\beta>\beta_c$). \\ 
\textbf{2)} The quantities of interest, mainly the times until cluster formation $t_{\text{cl}}$ or convergence $t^*$, are only meaningful when averaged over multiple independent trajectories. In our experiments, depending on the particular setup and hyperparameters, we needed up to 600 trajectories to resolve the averages and obtain sufficiently small confidence intervals (see, e.g., Fig. \ref{fig:convergence_times_1D}). 
\textbf{3)} For any system interacting via two-body forces (which is not just the case in our setting but also in the majority of molecular dynamics simulations in the scientific community \cite{benMDbook}), standard simulation algorithms are in $\mathcal{O}(N^2)$ due to the necessity to iterate over all particle pairs whenever the force is updated. This quadratic complexity leads to rapidly diverging runtimes when the system size is increased. For these reasons, computational speed may be the most important cornerstone in IPS simulation tools. Using an efficient C++ implementation in conjunction with highly accurate integrators (which allow for larger-than-usual stepsizes) constitutes an important stepping stone toward achieving controlled runtimes. 
SimIPS also employs multithreading via OpenMP to speed up the force computations. For the larger systems we consider in this work ($N\sim 10^3$), we observe an almost perfect scaling of the total runtime with respect to the number of threads up to 10 threads, reducing the runtime of a single trajectory by an order of magnitude. The downside of CPU-based shared-memory parallelization is that it does not scale up to arbitrarily large systems as the quadratic complexity in $N$ is still in place. For this reason, the molecular dynamics community often employs advanced algorithmic approaches like cell lists \cite{AllenTildesley2017} or Verlet lists \cite{Verlet} to improve the complexity and allow for more efficient parallelization. This, however, is only possible in cases where the interaction range is small compared to the spatial domain (which is not given in our case, see Sec. \ref{sup:sec:potentials}), and often only worthwhile for system sizes of $10^4$ or above. Therefore, we have not employed these techniques. Despite the quadratic complexity of the algorithm, we still manage to easily treat systems that are an order of magnitude larger than the ones in similar works ($N\sim10^3$ instead of $N\sim10^2$ as typically used \cite{garnier2017consensus}, \cite{Wang2017}, \cite{EversKolo2016}), and this even in three rather than one dimension.

\FloatBarrier
\subsection{Details on DBSCAN}\label{sup: sec: DBSCAN details}
DBSCAN (\textbf{D}ensity-\textbf{B}ased \textbf{S}patial \textbf{C}lustering of \textbf{A}pplications with \textbf{N}oise \cite{DBSCAN1}) is a prominent algorithm to classify a number of $N$ data points in $\mathbb{R}^d$ in terms of them being part of a cluster or not. The method solves many of the issues associated with older clustering methods. For example, unlike k-means or k-medoid methods, DBSCAN allows for points to be labeled as `noise' (i.e., as not belonging to any cluster), and it performs well even for nonconvex cluster shapes. Importantly, it does not require knowledge of how many clusters there are in the data. These properties are necessary to detect the onset of clustering in the IPS. DBSCAN's complexity is $\mathcal{O}(N^2)$ with $N$ the number of data points/particles in a given frame, but it is often faster than $\Theta(N^2)$ in practice, see the elaborate discussion in \cite{DBSCAN2}. While running the algorithm on ten-thousands of frames (drawn from hundreds of trajectories) is time-consuming, the comparably benign algorithmic complexity and the efficient implementation in Scikit-learn \cite{scikit-learn} still led to a tractable workload in our experiments. 

DBSCAN is based on a formal model of a cluster that closely aligns with human visual intuition: a cluster is a collection of points of high spatial density relative to the lower density of the local environment (consisting of `noise'). This intuition is formalized by various definitions of different types of points (core points, border points, noise) as well as spatial connectivity between them (density reachability/connectivity), see the original work \cite{DBSCAN1} for details. DBSCAN requires two hyperparameters, $N_{\min}$ and $\epsilon$. They determine whether a certain point is a core point or not, where core points are those points in the interior of a cluster. A point $p$ in data set $\mathcal{D}$ is defined to be a core point if $|\{q\in \mathcal{D}| d(q,p)\le \epsilon  \}| \ge N_{\min}$, where $d(\cdot,\cdot)$ is some user-defined metric, which in our case is just the Euclidean norm. A cluster is made up from core points that are `close enough' to one another and of border points that are `close enough' to the core points without being core points themselves (again, refer to \cite{DBSCAN1} for rigorous definitions). In particular, it can be inferred that the parameters $N_{\min}$ and $\epsilon$ define the lowest spatial density a cluster may have, namely $\propto N_{\min}/\epsilon^d$ when using the Euclidean norm in $\mathbb{R}^d$. 

\begin{figure}[!htbp]
\includegraphics[width=0.9\textwidth]{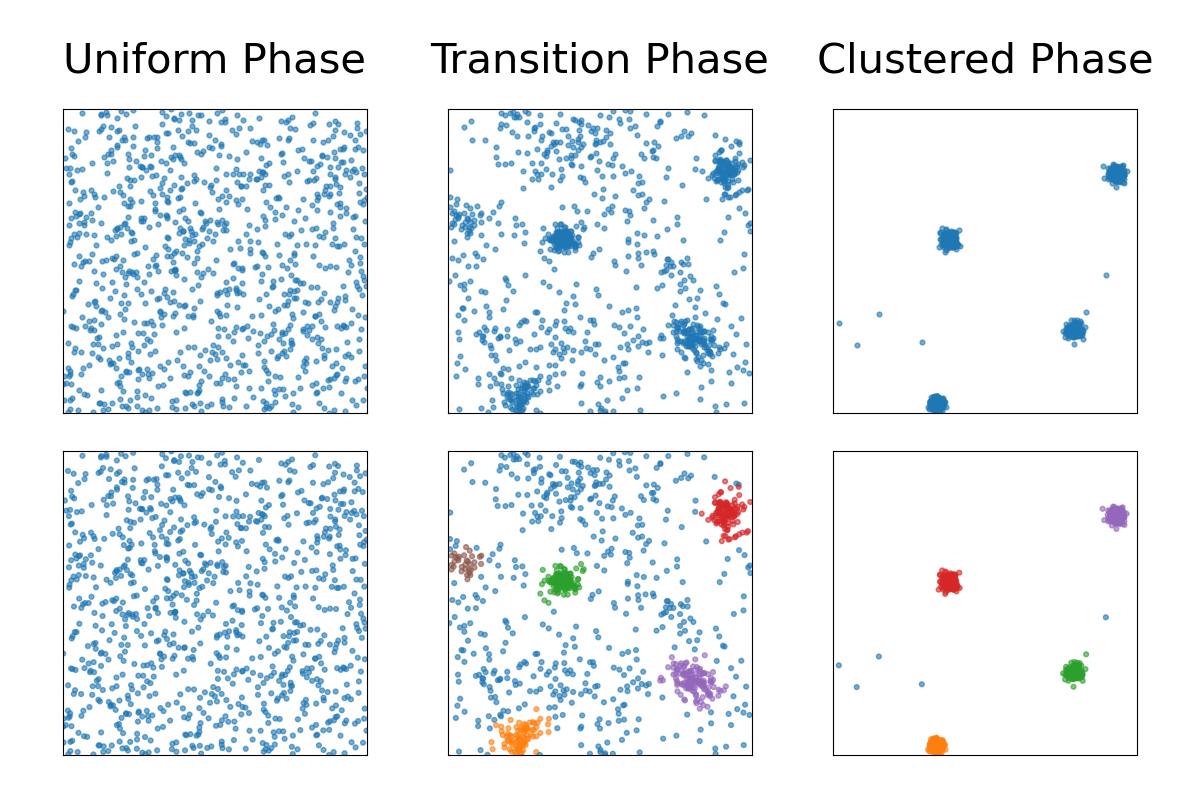}
\vspace{-0.25cm}
\caption{\label{sup:fig:dbscan1}Demonstration of DBSCAN cluster detection on a two-dimensional IPS for 1,000 particles at $\beta>\beta_c$. \textbf{Top:} Raw data. \textbf{Bottom:} Cluster detection results. The frames come from the same trajectory. DBSCAN parameters: $N_{\min}=25$, $\epsilon=0.25$.}
\end{figure}

Selecting the DBSCAN hyperparameters to reliably detect the onset of clustering in the systems we consider is difficult, in particular since the trajectories differ in terms of particle density and friction $\gamma$. Fortunately, the intuitive meaning of the hyperparameters allows for an obvious scaling of one of them with changing particle numbers $N$. Furthermore, while friction $\gamma$ does influence the times until clusters start to form, it does not influence the cluster shapes (with the exception of the $\gamma \to 0$ limit, see Sec. \ref{sup:sec:gamma=0}). This means that there is no need to repeatedly search for new DBSCAN hyperparameters whenever either $N$ or $\gamma$ is changed. In order to detect the onset of clustering in a given trajectory of the IPS, one needs to use values of the two DBSCAN hyperparameters such that no clusters are detected in the initial configuration (the uniform phase), but will be detected in the first clustered phase, as well as during the transition in between. Since the 'onset' of clustering is not rigorously defined, the aim is for DBSCAN to start detecting clusters in frames where human visual inspection starts to detect them as well. A well-working example for the two-dimensional IPS is shown in Fig. \ref{sup:fig:dbscan1}.
The first column shows the uniform phase (the initial configuration) in which DBSCAN does not find any clusters. That means that all particles are classified as noise as desired. The second column shows a frame during the transition from the uniform state to the clustered state. Several clusters are already clearly visible, but a large fraction of the particles are still free. DBSCAN detects the clusters reliably. Comparing this to the final column, which shows the first (metastable) cluster state in which almost all particles belong to a cluster, we observe that the number of detected clusters has decreased from 5 to 4. The leftmost detected cluster from column 2 is absent in column 3. This is not because DBSCAN was too sensitive to the configuration (i.e., detecting clusters where there clearly are none), but because the higher-density region corresponding to that 5th cluster dissolved again before the cluster could fully form. Indeed, from only looking at the raw data in column 2, even a human examiner might be unsure about whether a 5th cluster is about to form or not. Therefore it would not make sense for us to require DBSCAN to find the same number of clusters in the transition state and the first metastable state. This point is also obvious from the fact that we aim to detect the onset of clustering, and not all clusters will start to form at the same time. Note also that the co louring of the clusters in columns 2 and 3 is not identical. The colours are based on the order in which DBSCAN finds a cluster, rather than based on which particles make up the cluster. 

\begin{figure}[!htbp]
\includegraphics[width=1.0\textwidth]{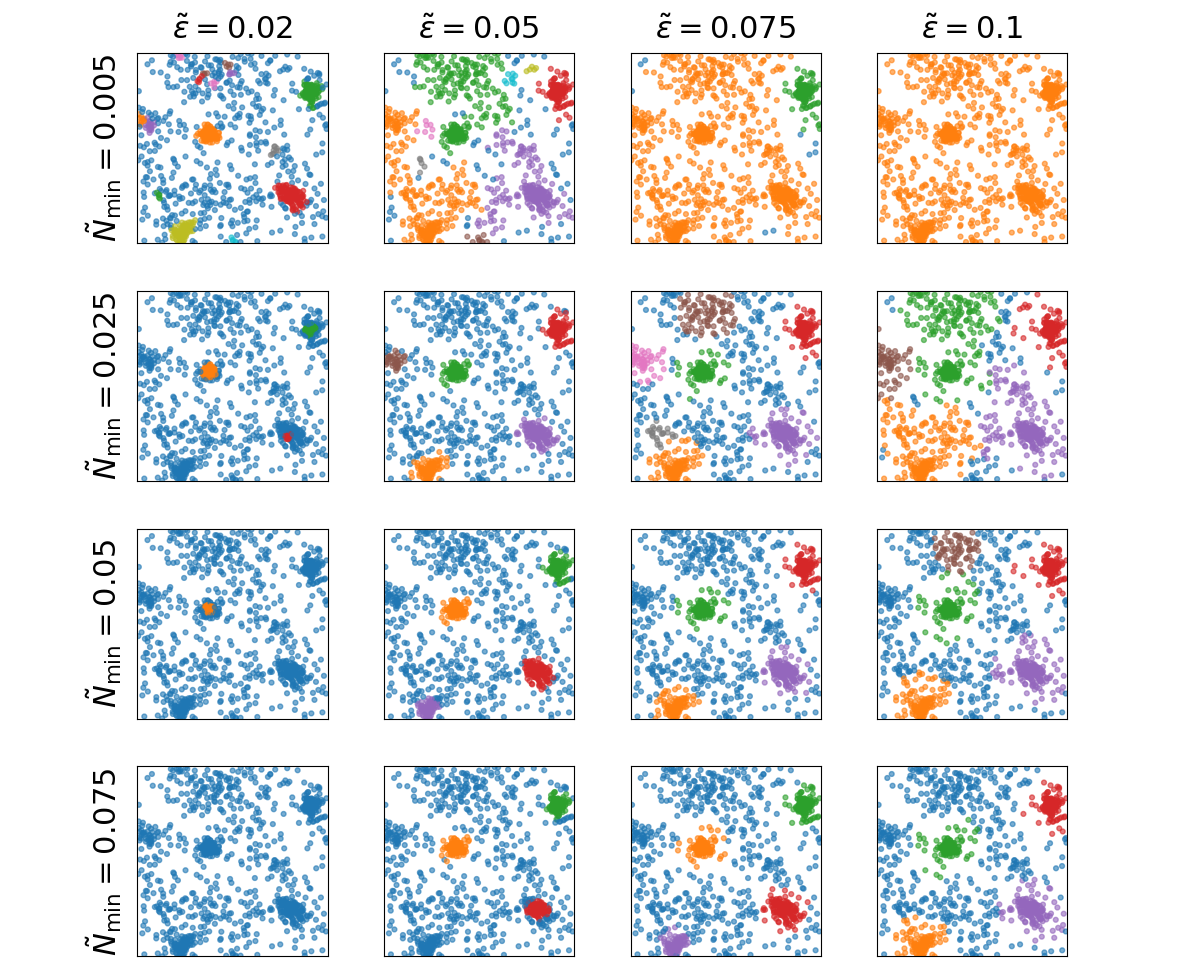}
\vspace{-0.5cm}
\caption{\label{sup:fig:dbscan2} Effect of DBSCAN hyperparameters on cluster detection. All panels show the same frame of 1,000 particles during the transition from the uniform to a metastable cluster phase.}
\end{figure}

In the following, we present our approach to find suitable DBSCAN hyperparameters to detect the onset of clustering in a given trajectory.
In our main experiment, Fig. \ref{fig:prefactor_double_plot_1D}, we considered the one-dimensional IPS.
Here, for the purpose of better visualization, we focus on the two-dimensional case (results for which are given in Fig. \ref{fig:results_2D} in Sec. \ref{sec: 2D results}).
We start by examining the influence of $N_{\min}$ and $\epsilon$ on the cluster detection result for the transition-phase frame used in Fig. \ref{sup:fig:dbscan1}. 
The two parameters have natural dimensionless forms. We set $\tilde{N}_{\min}:=\frac{N_{\min}}{N}$ and $\tilde{\epsilon}:=\frac{\epsilon}{L}$ with number of particles $N$ and box length $L$. The results are shown in Fig. \ref{sup:fig:dbscan2}.

The results vary in terms of the number of clusters\footnote{Note that the word `cluster' is ambiguous in this section. It may refer to an actual cluster formed by the particles or to the detection result of DBSCAN, which labels a certain group of particles a `cluster' (whether reasonable or not).} detected as well as cluster size.
For small $\tilde{\epsilon}$ and large $\tilde{N}_{\min}$, no particle qualifies as a core point, hence no cluster is detected (see the bottom left panel).
Decreasing $\tilde{N}_{\min}$ and/or increasing $\tilde{\epsilon}$ increases the number of detected clusters. It also increases the cluster size, but through different channels. Decreasing $\tilde{N}_{\min}$ for fixed $\tilde{\epsilon}$ leads to particles becoming core points, forming inner points of clusters. The clusters become larger for decreasing $\tilde{N}_{\min}$ since points that neighbour core points may become core points as well. If their distances are smaller than $\epsilon$, they will be part of the same cluster. This process can be neatly observed in the centre two panels of the $\tilde{\epsilon}=0.02$ column. Increasing $\tilde{\epsilon}$ for fixed $\tilde{N}_{\min}$ increases cluster sizes by allowing for larger distances between border points and the interior of a cluster, leading to structures with smoother cluster-to-noise transitions. If changing the DBSCAN hyperparameters leads to cluster growth or new clusters arising between two already detected clusters, the clusters may merge to form bigger clusters. Therefore, while decreasing $\tilde{N}_{\min}$ and increasing $\tilde{\epsilon}$ tend to increase the number of detected clusters at first, beyond certain threshold values, that number decreases again due to cluster merging. The extreme case is seen in the top right panel, in which all particles in the system are read as one large cluster.

\begin{figure}[!htbp]
\includegraphics[width=0.75\textwidth]{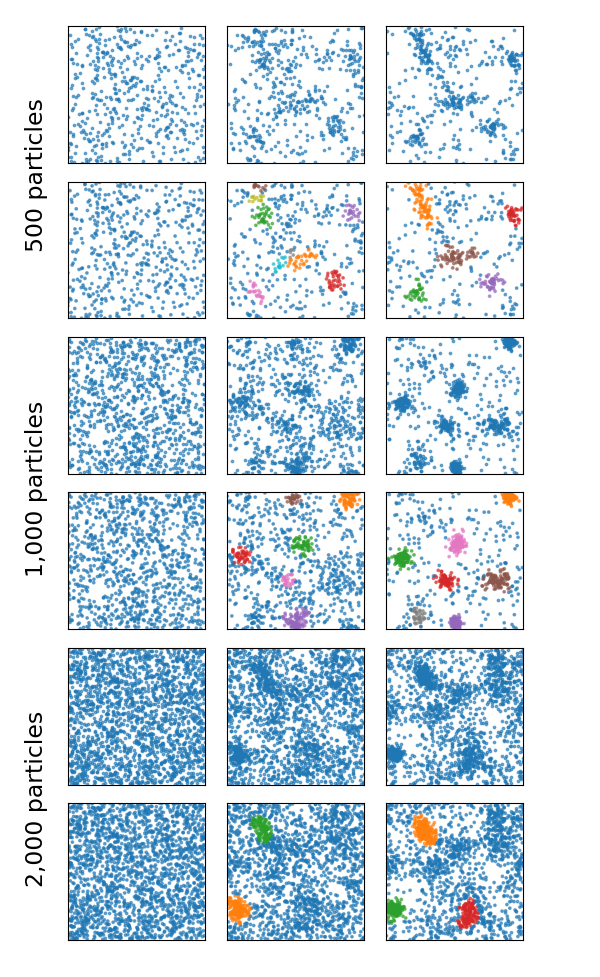}
\vspace{-0.5cm}
\caption{\label{sup:fig:dbscan3} DBSCAN results for different particle numbers $N$ (with fixed box size) when using the linear scaling $N_{\min}=N_{\min,0}\frac{N}{N_0}$ for $N_{\min,0}=25$ and $N_0=1000$, such that $\tilde{N}_{\min}=0.025$ remains constant. $\tilde{\epsilon}=0.05$ for all runs. The first row for each particle number shows raw data, the second row DBSCAN detection results. }
\end{figure}

\begin{figure}[!htbp]
\includegraphics[width=1.0\textwidth]{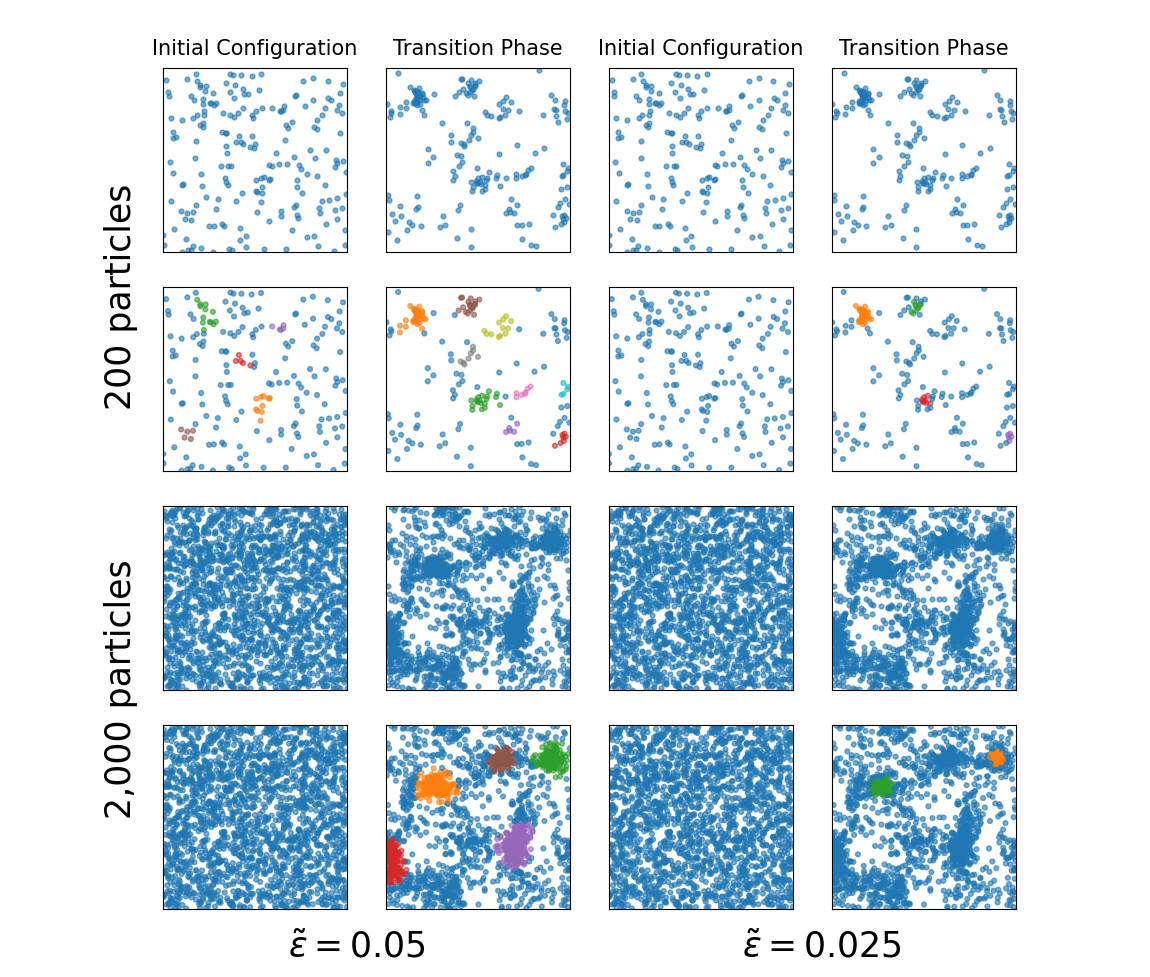}
\vspace{-0.5cm}
\caption{\label{sup:fig:dbscan4} DBSCAN results for two different particle numbers $N$ and two different $\tilde{\epsilon}$. $\tilde{N}_{\min}$ was scaled linearly as in Fig. \ref{sup:fig:dbscan3}.  The first row per particle density shows raw data, the second row DBSCAN results. }
\end{figure}

Judging by Fig. \ref{sup:fig:dbscan2}, the setting $(\tilde{N}_{\min},\, \tilde{\epsilon})=(0.025, \, 0.05)$ leads to good agreement of the DBSCAN detection results with the actual particle clusters visible in the frame and would therefore be a good candidate to use for further experiments. Indeed, the usage of relatively small $\tilde{\epsilon}$ values aligns with the recommendation found in the literature \cite{DBSCAN2}. Note that in our main experiments, we needed to run DBSCAN on frames with different particle densities. Hence, in order to keep $\tilde{N}_{\min}$ constant across different particle counts, we rescaled $N_{\min}$ linearly with $N$, i.e., $N_{\min}=N_{\min,0}\frac{N}{N_0}$, for base values $N_{\min,0}$ and $N_0$. In Fig. \ref{sup:fig:dbscan3}, we experimentally verify that this linear scaling leads to reasonable results on various particle densities. For each particle density examined, we show DBSCAN results for three frames: the initial uniform configuration (left column), an early stage of the transition phase (centre column), and a late stage of the transition phase (right column). The frames were chosen via visual inspection of the trajectory. We require DBSCAN to detect only noise in the first frame, and clusters in the second and third frame, respectively. This ensures that the first frame DBSCAN detects a cluster in  corresponds to the 'onset' of clustering for a given trajectory.
While Fig. \ref{sup:fig:dbscan3} shows that parameters $(\tilde{N}_{\min},\, \tilde{\epsilon})=(0.025, \, 0.05)$ together with the proposed linear scaling for $N_{\min}$ work reasonably well, one should consider multiple trajectories per particle density and also include the most extreme particle densities of interest to find the most suitable DBSCAN setting. We proceeded as follows. For particle numbers $N\in \{200, 400, 500, 1000, 2000, 4000 \}$, we simulated 5 trajectories per $N$. As before, we visually inspected each trajectory and extracted two test frames: one from an early stage of the transition phase, one from a later stage. We ran DBSCAN on each of the obtained frames as well as the corresponding initial configuration and counted the number of clusters detected. The aim was to find values for $(\tilde{N}_{\min},\, \tilde{\epsilon})$ such that no cluster is ever detected in the initial configuration, but for each trajectory at least one cluster is detected in at least one of the two transition-phase frames. Doing that revealed that the values used in Fig. \ref{sup:fig:dbscan3} were too sensitive for small particle numbers, for which clusters were detected even in the initial configurations. Decreasing $\tilde{\epsilon}$ to 0.025 fixed that issue at the cost of being more likely to miss out on clusters in the transition phase (in particular for larger particle numbers), see Fig. \ref{sup:fig:dbscan4}.
Since the actual `onset' of clustering is not rigorously defined anyway, we deem it more important to prevent the algorithm from seeing clusters in the uniform phase instead of ensuring that all actual particle clusters are detected in the transition phase. Hence, for our 2D experiments presented on the right-hand side of Fig. \ref{fig:results_2D} in Sec. \ref{sec: 2D results}, we use $(\tilde{N}_{\min},\, \tilde{\epsilon})=(0.025, \, 0.025)$. We used the same process to find suitable DBSCAN parameters for the 1D experiments presented in the main text. Finally, we remark that we did no take proper account of periodic boundaries when analysing a frame for clusters, which we do not assume to make a statistically significant difference in the obtained results of the main experiments.

While DBSCAN is able to reliably determine not just clustered states but even the onset of clustering across all trajectories and system parameter ranges that we examined in this work, it comes with a few drawbacks. First, it requires storing of entire trajectories. Given that for each $N$ and $\gamma$ we need to run hundreds of independent trajectories, the memory footprint of this approach is rather large. Despite making sure that the trajectories are as short as possible (long enough so that all trajectories at given $N$ and $\gamma$ reach the first clustered state, but not much longer), the experiments to detect the onset of clustering already produced 500GB of data. Second, the compute time to run DBSCAN on all resulting frames is substantial as well. In our case, a run with \textbf{already-finetuned} hyperparameters to create the plots in this work takes multiple hours. A third point is the already mentioned challenging selection of hyperparameters. While the hyperparameters are robust to changes in $\gamma$ and $N$ (apart from scaling, see \ref{sup: sec: DBSCAN details}), they are not robust to changes in temperature or the interaction potential, as both of them modify the shape and density of clusters. Since we used constant temperature for our experiments and three potentials with comparable width and depth, this did not become an issue for us. 
Finally, like most machine-learning algorithms, DBSCAN is a black-box approach. Its determination of the onset of clustering does not rely on the underlying physics of the system. That means in particular that the obtained times are arbitrary in the sense that different DBSCAN hyperparameters lead to different times. Therefore, the actual time values cannot be expected to completely reproduce theoretical predictions. 
In order to circumvent many of these issues, it would be desirable to use a dynamical, optimally scalar quantity to detect the onset of clustering, similar to the criterion (\ref{eq: d_com_criterion}) for the time until full convergence. We examine various ways in which the mean squared displacement (MSD) or the total potential energy in the system can be used to detect clusters in Secs. \ref{sup:sec: MSD_details} and \ref{sup:sec:energy_vs_msd}.
While more memory and compute efficient than DBSCAN, we found that their success depends more strongly on the system parameters (in particular the friction). We leave the question of the best way to detect the onset of clustering for future research.

\FloatBarrier

\subsection{Cluster detection with the mean squared displacement.} \label{sup:sec: MSD_details}

\begin{figure}[!htbp]
\includegraphics[width=1.0\textwidth]{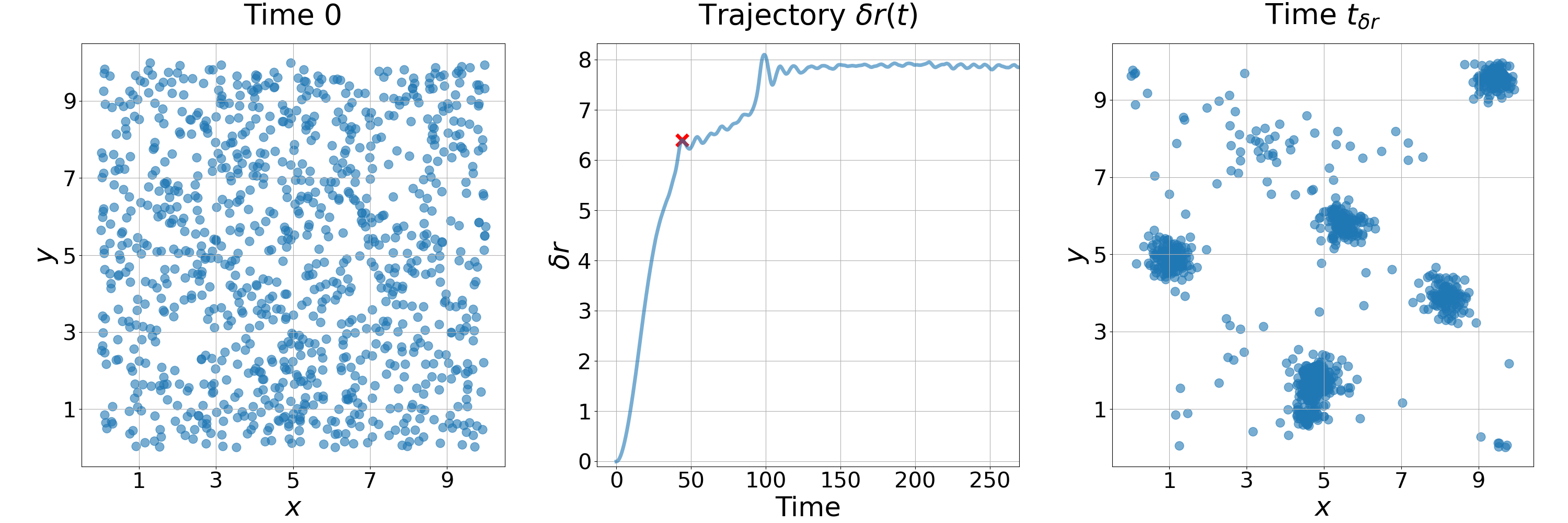}
\caption{\label{sup:fig:msd_demo1} Measuring the time of the first clustered state of a 2D IPS using $t_{\delta r} $ from \eqref{sup:eq:msd_derivative_criterion}. \textbf{Left:} Initial uniform configuration. \textbf{Centre:} Time series of the mean squared displacement. The marker denotes the point at time $t_{\delta r}$. \textbf{Right:} Configuration at time $t_{\delta r}$. \textbf{Simulation parameters:} $\gamma=0.1$, $\beta=150\approx 5 \beta_c$. Gaussian interaction potential.}
\end{figure}

\begin{figure}[!htbp]
\includegraphics[width=1.0\textwidth]{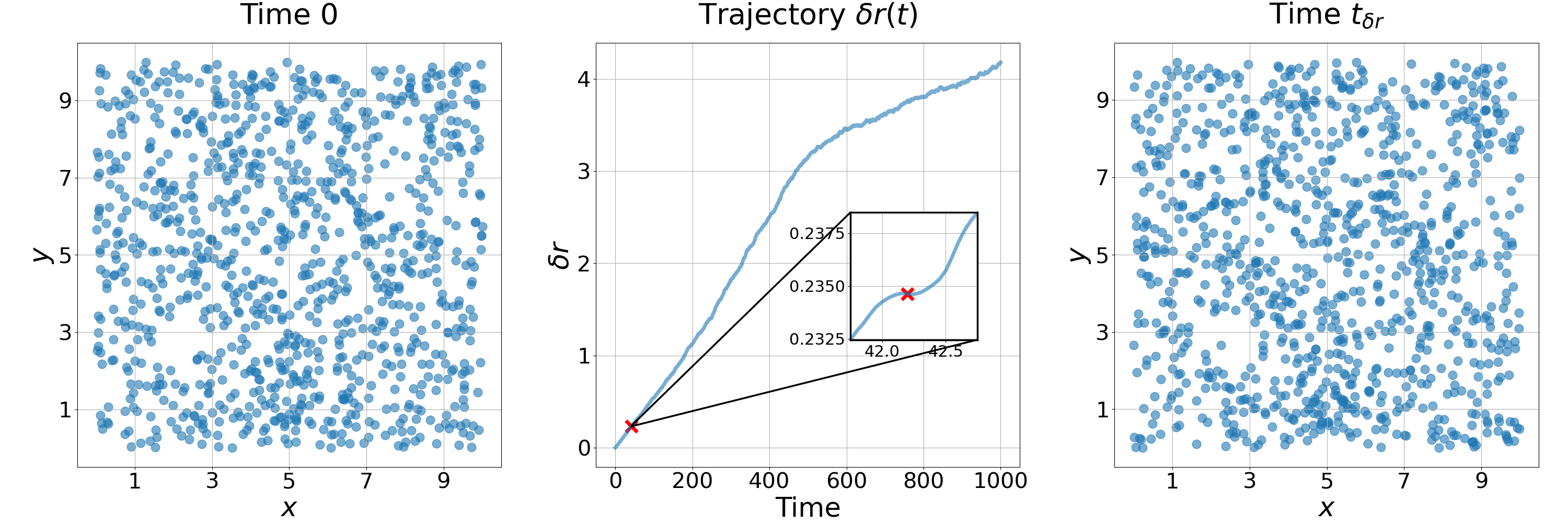}
\caption{\label{sup:fig:msd_demo2} Measuring the time of the first clustered state of a 2D IPS using $t_{\delta r} $ from \eqref{sup:eq:msd_derivative_criterion}. \textbf{Left:} Initial uniform configuration. \textbf{Centre:} Mean squared displacement. The marker denotes the point at time $t_{\delta r}$, which is detected due to noise rather than an actual plateau (see inset). \textbf{Right:} Configuration at time $t_{\delta r}.$ \textbf{Simulation parameters:} $\gamma=5$, $\beta=150\approx 5 \beta_c$. Gaussian interaction potential.}
\end{figure}

\begin{figure}[!htbp]
\includegraphics[width=1.0\textwidth]{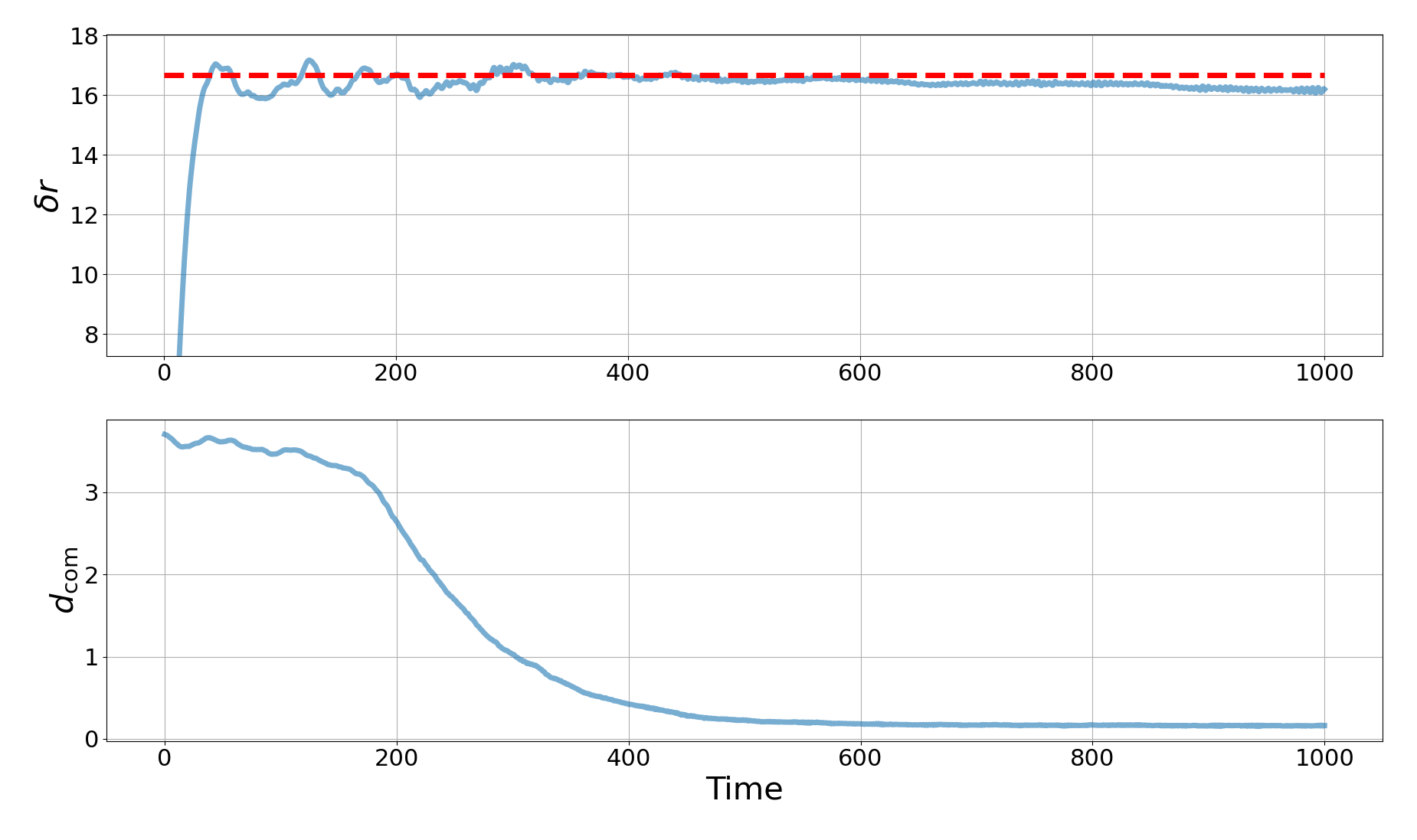}
\caption{\label{sup:fig:msd_demo3} Comparing $\delta r$ (\textbf{top}) with $d_{\text{com}}$ (\textbf{bottom}) at high temperatures below the critical temperature, $\beta \approx 1.06 \beta_c$. The red dashed line denotes the equilibrium value of $\delta r$ given by $\frac{L^2}{6} = 16.67$ (see text).}
\end{figure}

A quantity that is often used to detect qualitative changes in the dynamics of $N$-particle systems is the mean squared displacement (MSD). Examples mainly come from the physical sciences, e.g., \cite{msd_ref1,msd_ref2,msd_ref3,msd_ref4}, but these works do not consider clustering processes. The MSD is given by 
\begin{equation} \label{eq: MSD}
\delta r (t)\coloneqq \frac{1}{N}\sum_{i=1}^{N} \big(x_i(t) - x_i(0)\big)^2.
\end{equation}
Importantly, this is the MSD along a single trajectory (not to be confused with the MSD obtained via averaging over independent trajectories as is typically done when simulating equilibrium states). After an initial steady growth, the MSD adopts pronounced plateaus whenever a metastable clustered state is reached, jumping from one plateau to a higher one when two clusters merge, see Fig. \ref{fig:observbles_vs_gamma}. Therefore, one way of detecting the time of first clustered state is to measure
\begin{equation}\label{sup:eq:msd_derivative_criterion}
 t_{\delta r}:=\min\{t\ |\  \frac{\mathrm{d}}{\mathrm{dt}} \delta r(t)<0 \}.
\end{equation}
 This is illustrated in Fig. \ref{sup:fig:msd_demo1} where the right-hand panel shows a snapshot of the configuration at $t_{\delta r}$. 
The configuration has clearly entered a clustered state, despite many particles still being free.
From looking at $\delta r$ plotted in the centre panel, the time $t_{\delta r}$ is obtained at the beginning of the first `plateau'. The plateau is not actually horizontal as it is short-lived with $\delta r$ quickly jumping to a higher plateau. This implies that two or more of the visible clusters in the frame on the right-hand side merged shortly after $t_{\delta r}$. The ability to spot the beginning of the first clustered state even if it is a transient one makes criterion \eqref{sup:eq:msd_derivative_criterion} look promising at first glance. However, it has a few drawbacks. First, while it can be used to detect the beginning of the first clustered state, it would be a stretch to call the configuration on the right-hand side of Fig. \ref{sup:fig:msd_demo1} the `onset' of clustering. Clearly, the state is already more clustered than uniform. Second, the reliability of using \eqref{sup:eq:msd_derivative_criterion} depends on the hyperparameters of the simulation. In Fig. \ref{sup:fig:msd_demo1} we used small friction $\gamma$ and temperature $\beta^{-1}$. For larger frictions, the $\delta r$ dynamics becomes more noisy which can lead to cases where $\frac{\mathrm{d}}{\mathrm{dt}} \delta r(t)\!<\!0$ is merely due to noise, not because a clustered state has been reached, see Fig. \ref{sup:fig:msd_demo2}.
\\
Another issue arises from the fact that the MSD $\delta r$ is bounded in a bounded domain. Even in the case of $\beta \! < \!\beta_c$, where no clustering would occur, $\delta r$ will reach a plateau at one point. The value of this plateau can be estimated:
We assume $\beta \ll \beta_c$ for simplicity. For the force we then have $\nabla W \! \approx \! \boldsymbol{0}$, and all particles are become i.i.d. Brownian walkers. Assume the one-dimensional case. Taking expectation of \eqref{eq: MSD},
\begin{equation}\label{sup:eq:E_msd}
\begin{split}
E\Big( \frac{1}{N}\sum_{i=1}^{N} \big(X_i(t) - X_i(0)\big)^2 \Big) &= E\Big( \big(X_k(t) - X_k(0)\big)^2 \Big),
\end{split}
\end{equation}
where index $k$ is arbitrary.
For long times, and due to the periodic boundaries, $X_k(t)-X_k(0)$ becomes uniformly distributed, $X_k(t)-X_k(0)\sim \mathcal{U}([-\frac{L}{2}, \frac{L}{2}]$.  Thus, we have 
\begin{equation}
\mathbb{E}\Big( \big(X_k(t) - X_k(0)\big)^2 \Big) = \text{Var}\Big( X_k(t) - X_k(0) \Big) = \frac{L^2}{12}, 
\end{equation}
which in 2D becomes $\frac{L^2}{6}$ due to independency of dimensions. This result also holds for smaller temperatures, even below the critical temperature. This is because the stationary state at $\beta > \beta_c$ is the one-cluster state, where the cluster's centre of mass evolves like a Brownian walker for which the same calculation applies. Together with the fact that individual particles will remain in the vicinity of the cluster's centre of mass, the same equilibrium value of $\delta r$ is obtained (demonstrated experimentally in Fig. \ref{sup:fig:msd_demo3} and also on the right-hand side of Fig. \ref{fig:observbles_vs_gamma}).\\
At $\beta \gtrsim \beta_c$, the beginning of cluster formation can take prolonged times, as more particles are necessary to randomly get close to one another for the attractive forces to become relevant enough. In this case, $\delta r$ might reach its equilibrium value before clustering starts, leading to incorrect results for $t_{\delta r}$. 
This situation can be observed in Fig. \ref{sup:fig:msd_demo3}, which compares $\delta r$ with $d_{\text{com}}$ from Sec. \ref{sec:simulation}, the quantity that can be used to detect the final one-cluster state. It can be seen that $\delta r$ reaches its equilibrium value long before $d_{\text{com}}$ shows any significant decline. Indeed, from looking at $\delta r$ only, the only evidence of clusters forming is the transition from noisy and chaotic oscillations in the curve to much smaller and more systematic oscillations at time $\mathord{\sim}350$, denoting a stronger correlation of the movements of the individual particles in the system.
Therefore, criterion \eqref{sup:eq:msd_derivative_criterion} will fail to properly predict the onset of clustering at large frictions and temperatures. Still, given the importance of dynamical quantities like the MSD to detect subdiffusive behavior in the physical sciences (see again the references at the start of this section), we believe it should be possible to make use of MSD-like quantities instead of machine learning-based algorithms to detect the transition from the uniform to the clustered state in our systems. We leave this to future research.

\FloatBarrier

\subsection{Cluster detection with the total potential energy.}\label{sup:sec:energy_vs_msd}
\begin{figure}[!htbp]
\includegraphics[width=1.0\textwidth]{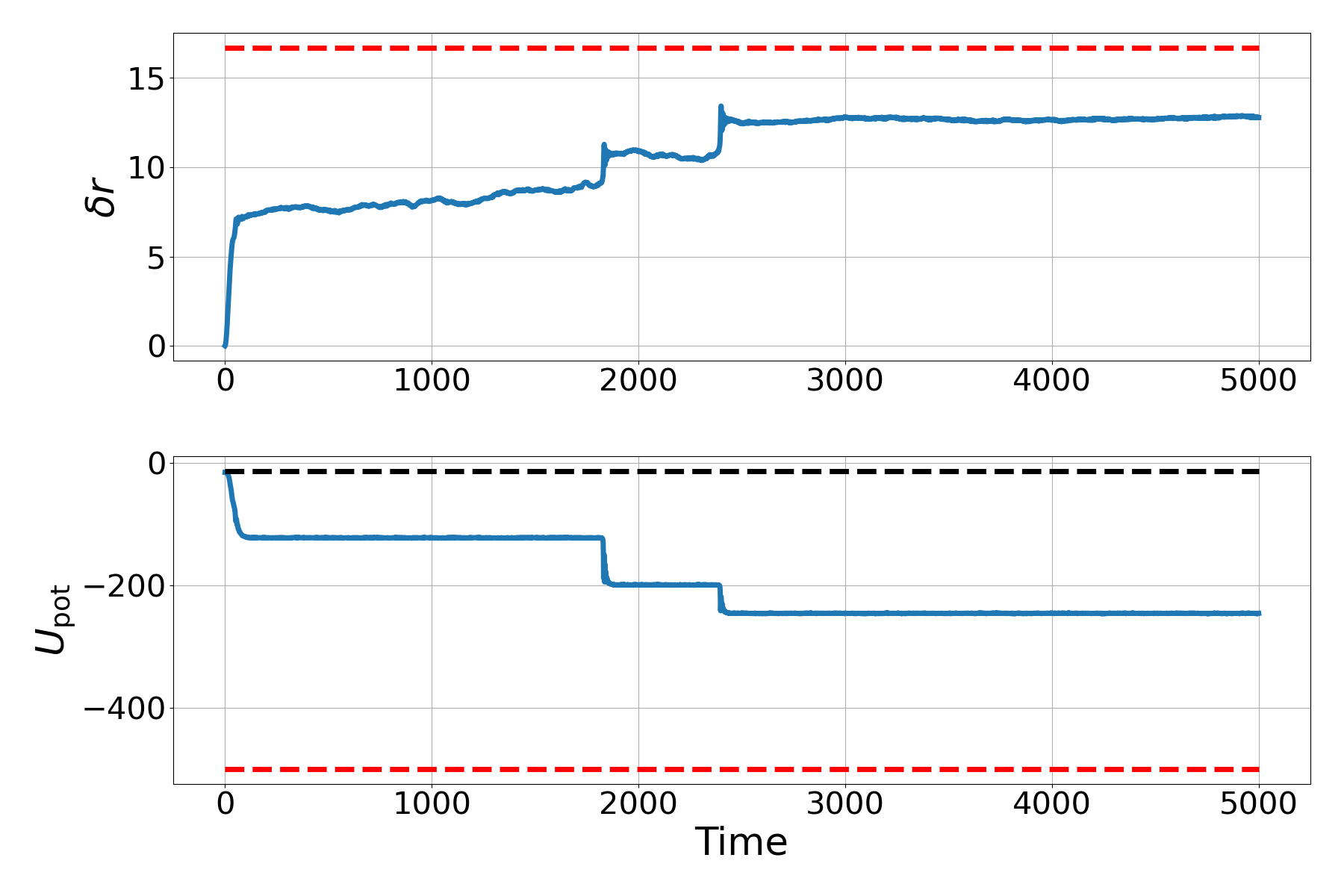}
\caption{\label{sup:fig:energy_vs_msd} The mean squared displacement $\delta r$ (\textbf{top}) and the total potential energy $U_{\text{pot}}$ (\textbf{bottom}) of a 2D particle system. The red dashed line denotes $U_{\text{pot}}^{\min}$ and the equilibrium value of $\delta r$, respectively. The black dashed line denotes the estimated initial potential energy $U_{\text{pot}}^0$ (see text). Simulation parameters: 1,000 particles, Gaussian interaction with $\sigma^2=0.5$, $\gamma=0.1$, $\beta=150>\beta_c$. }
\end{figure}

In Sec. \ref{sup:sec: MSD_details} we have considered the mean squared displacement $\delta r$ (MSD) as a dynamical quantity potentially useful to detect the beginning of the clustering process. Here, we briefly highlight another quantity of interest for that aim, the potential energy $U_{\text{pot}}$ of the system. Using the notation from \eqref{eq: Langevin Dynamics}, the total potential energy in the system is given by 
\begin{equation}
U_{\text{pot}}:= \frac{1}{2N} \sum_{i \neq j} W(\boldsymbol{x}_i, \boldsymbol{x}_j), \label{sup:eq:Upot}
\end{equation}
with pairwise interaction potential $W(\boldsymbol{x}_i, \boldsymbol{x}_j)$. The potential energy is known as a so-called `order parameter' for these types of systems, showing distinct differences between the uniform phase and the various metastable clustered states (see, e.g., Fig. 7 in \cite{martzel2001mean} or Fig. 1 in \cite{Wang2017} for a simpler model). The potential energy of the systems we consider is bounded from above and below (similar to the MSD $\delta r$ and the mean distance to the centre of mass $d_{\text{com}}$ on a bounded simulation domain). For the three choices of $W$ we consider in this work (see Sec. \ref{sup:sec:potentials}), the lowest-energy state is given for the configuration in which all $\boldsymbol{x}_i$ are identical, given by $U_{\text{pot}}^{\min}\!=\!-0.5 (N-1)$. The highest-energy state is less obvious to determine. It likely consists of a crystalline configuration. However, we can write down an estimate for the potential energy of the initial conditions we use, i.e., the uniform configuration. Denote by $\mathbb{E}^d_0$ the expectation with respect to the uniform measure on $[0,L]^d$. For the Gaussian potential $W(\boldsymbol{x}, \boldsymbol{y}) \!=\! - e^{-\frac{\|\boldsymbol{x}-\boldsymbol{y} \|^2}{2\sigma^2}}$, due to the independency of dimensions, we have
\begin{equation} \label{sup:eq:expectation_relation}
\mathbb{E}^d_0\big(W(\boldsymbol{x}, \boldsymbol{y})\big)=-\Big[\mathbb{E}^1_0\Big(e^{-\frac{|x-y|^2}{2\sigma^2}}\Big)\Big]^d,
\end{equation}
with scalar positions $x$ and $y$. For the total potential energy in the uniform configuration, using \eqref{sup:eq:expectation_relation} and the fact that all particles are i.i.d, we then have 
\begin{align}
U_{\text{pot}}^0:=\mathbb{E}^d_0(U_{\text{pot}}) &= \mathbb{E}^d_0\Big( \frac{1}{2N} \sum_{i \neq j} W(\boldsymbol{x}_i, \boldsymbol{x}_j) \Big), \\
&= \frac{N-1}{2} \mathbb{E}^d_0\big(W(\boldsymbol{x}_i, \boldsymbol{x}_j)\big), \\
&= -\frac{N-1}{2} \Big[\mathbb{E}^1_0\Big(e^{-\frac{|x-y|^2}{2\sigma^2}}\Big)\Big]^d, \\
&= -\frac{N-1}{2} \bigg(\int_{[0,L]}\int_{[0,L]}\rho_0(x_i,x_j)e^{-\frac{|x-y|^2}{2\sigma^2}}\text{d}x_i\text{d}x_j \bigg)^d,
\end{align}
with uniform measure $\rho_0(x_i, x_j)\!:=\!1/L^2$. The integral can be evaluated via quadrature. For $d\!=\!2$, $\sigma^2\!=\!0.5$, $N\!=\!1000$, and $L\!=\!10$, we obtain $U_{\text{pot}}^0\!=\!-13.97$. Note that we have neglected periodic boundary conditions in this calculation, but the introduced relative error is negligible for small $\sigma/L$.
Fig. \ref{sup:fig:energy_vs_msd} plots the evolution of the potential energy at these system parameters in comparison to the MSD at $\beta \!> \!\beta_c$. 
We see that $U_{\text{pot}}$ roughly evolves like a step function whose plateaus are aligned with the plateaus in $\delta r$. From this, we conclude that $U_{\text{pot}}$ can also be used to detect clustered states, where each plateau corresponds to one metastable multi-cluster state. As clusters merge, $U_{\text{pot}}$ drops to a lower plateau (in contrast to the MSD which jumps to a higher plateau). The stationary one-cluster state would correspond to the lowest plateau (assuming the trajectory has been simulated long enough). The one-cluster plateau would still sit at a level above $U_{\text{pot,min}}$ (red dashed line in the figure), since the thermal noise will always ensure that the clusters have finite width. From this it becomes clear that $U_{\text{pot}}$, just like $\delta r$, cannot easily be used to detect the one-cluster phase, at least not nearly as good as $d_{\text{com}}$. Similar to $\delta r$, however, it may be useful to detect the onset of clustering. On one hand, the figure shows that it might be the more useful quantity than $\delta r$, as its plateaus are overall less noisy. On the other hand, using $U_{\text{pot}}$ to detect the onset of clustering runs into similar issues we discussed for $\delta r$ in Sec. \ref{sup:sec: MSD_details}, and hence we refrain from using it in our main experiments.

\FloatBarrier

\subsection{Number of trajectories for the two-dimensional results.}\label{sup:sec:2D_results_table}
\par\noindent
\begin{table}[H]
    \centering
    \begin{minipage}[t]{0.45\textwidth}
    \centering
    \caption*{Gauss Potential}
    \vspace{-0.25cm}
    \begin{tabular}{|c|c|}
        \hline
        $\gamma$-range                 & $N_{\text{traj}}$      \\ \hline \hline
        $0.1 < \gamma$                 & 400                    \\ \hline
        $0.015 \leq \gamma \leq 0.1$   & 300                    \\ \hline
        $0.0025 \leq \gamma <0.015 $   & 100                    \\ \hline
        $ \gamma < 0.0025$             & 50                     \\ \hline
    \end{tabular}
    \end{minipage}
    \begin{minipage}[t]{0.45\textwidth}
    \centering
    \caption*{Morse Potential}
    \vspace{-0.25cm}
    \begin{tabular}{|c|c|}
        \hline
        $\gamma$-range                & $N_{\text{traj}}$ \\ \hline \hline
        $ 0.075 < \gamma $            & 400               \\ \hline
        $0.02 \leq \gamma \leq 0.075$ & 100               \\ \hline
        $0.00375 \leq \gamma<0.02 $   & 50                \\ \hline
        $\gamma < 0.00375$            & 10                \\ \hline
    \end{tabular}
    \end{minipage}
    \caption{The number of trajectories $N_{\text{traj}}$ that were averaged to generate the $t^*_{\sigma}$-data in Fig.~\ref{fig:results_2D}, depending on friction value $\gamma$.}
    \label{tab:avg_scheme_2D}
\end{table}

\FloatBarrier

\subsection{Details on stepsize selections}\label{sec: choice of stepsize}
\begin{figure}[!htbp]
\centerline{\includegraphics[width=1\textwidth]{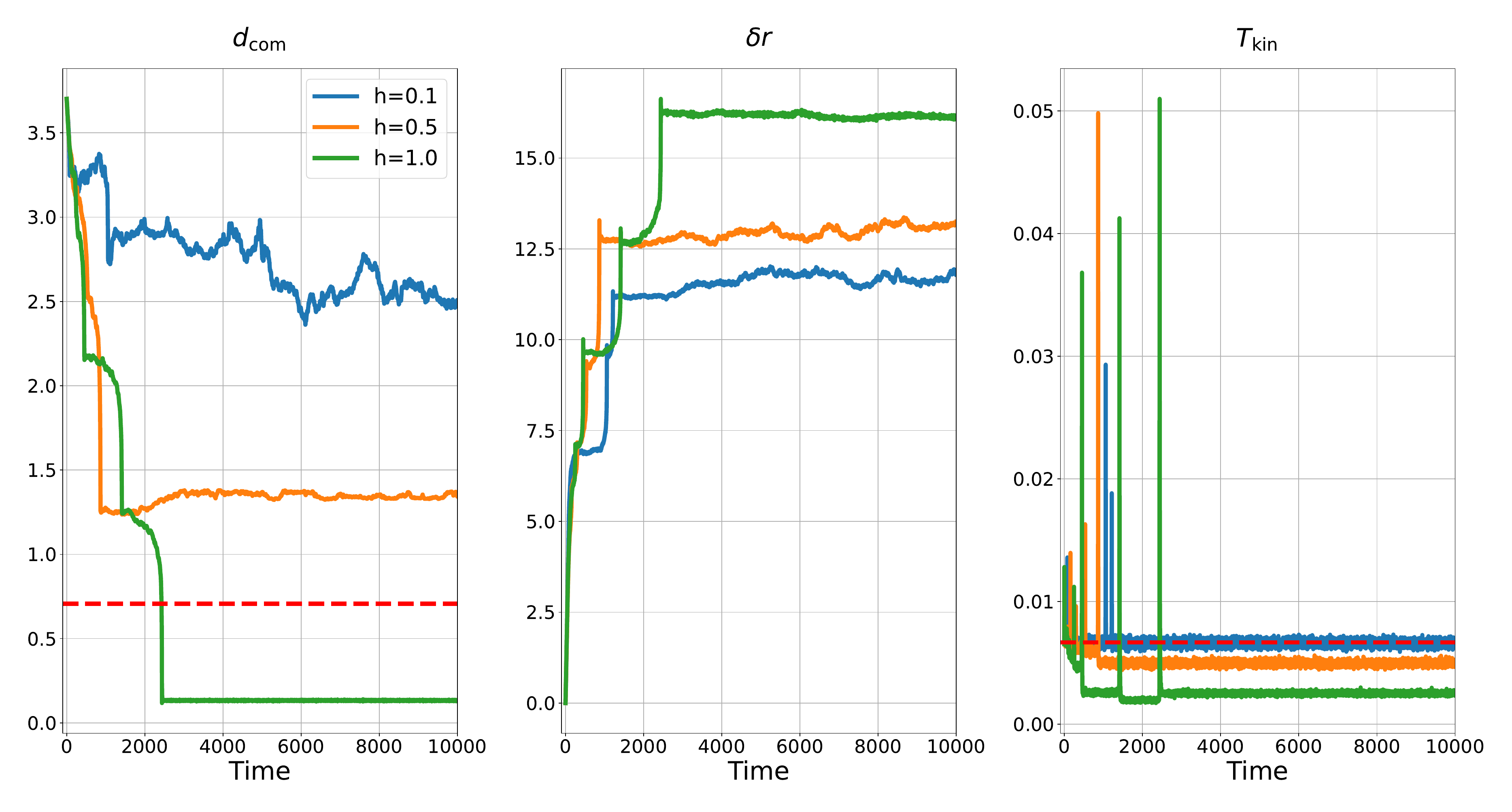}}
\caption{2D simulation results for $d_{\text{com}}$, $\delta r$, and $T_{\text{kin}}$ for three different stepsizes $h$. Dashed lines denote convergence threshold $\sigma=\sqrt{0.5}$ in the $d_{\text{com}}$ plot and the system temperature $T=\frac{1}{\beta}=0.0067$ in the $T_{\text{kin}}$ plot. Simulation settings: Morse potential, $N=1000$, $\gamma=0.5$, $\beta=150$.}
\label{fig:effect_of_stepsize_observables}
\end{figure}

When integrating dynamical systems, the choice of stepsize influences simulation accuracy, stability, and efficiency (see, e.g., \cite{LeimkuhlerReichBook}). 
In IPS simulations and molecular dynamics, the processes of interest are often slow or rare in occurrence, such that the times that need to be simulated are long \cite{HMR_MD,Lagardere2019}. A large stepsize reduces the number of integrator iterations needed to simulate these times. In particular, it requires fewer compute intensive force calculations. In our case, we are interested in the formation of clusters starting from an initial, unstable uniform distribution. Depending on the system parameters (mainly temperature and interaction range) the formation of the first clusters can take arbitrarily long times. Furthermore, the convergence to equilibrium (the one-cluster stationary state) can be arbitrarily slow due to the existence of metastable multi-cluster states. In such a state, multiple clusters randomly meet and merge to either reach a new multi-cluster state (with one cluster less) or the one-cluster equilibrium state. Since clusters behave as independent Brownian walkers (given sufficient distance between them \cite{garnier2017consensus}), depending on the system parameters, merge events can be arbitrarily rare. This point becomes more important in higher spatial dimensions, since the chance of two random walkers meeting decreases with increasing dimension. Given the $\mathcal{O}(N^2)$ complexity of the integration of an $N$-particle system with pairwise interactions and given the fact that we typically have to average over an ensemble of independent trajectories per setting (per friction $\gamma$ or particle count $N$ in our experiments) to obtain reliable statistical estimates, the required compute time is substantial. Increasing the integrator stepsize can reduce that compute time linearly. However, an increasing stepsize leads to less accurate integration, which in the case of SDEs induces a bias in the sampled distribution (see, e.g., Ch. 7 in \cite{benMDbook}). The choice of integrator matters for that reason, as higher-quality integrators like BAOAB or UBU (see Sec. \ref{sec:numerical_integrators}) suffer less from this effect than standard ones. Furthermore, large stepsizes can lead to unstable behavior either in the physical sense (unphysical trajectories) or the computational sense (numerical overflows). It is therefore important to strike a suitable balance between compute work and simulation accuracy and stability. \\
To illustrate this, we run a simulation in 2D at $\beta>\beta_c$ for three different stepsizes and plot the observables presented in Sec. \ref{sec:simulation}: the mean distance to the centre of mass $d_{\text{com}}$, the mean squared displacement $\delta r$, and the kinetic temperature $T_{\text{kin}}$. Fig. \ref{fig:effect_of_stepsize_observables} shows the results.

\begin{figure}[!htbp]
\centerline{\includegraphics[width=1\textwidth]{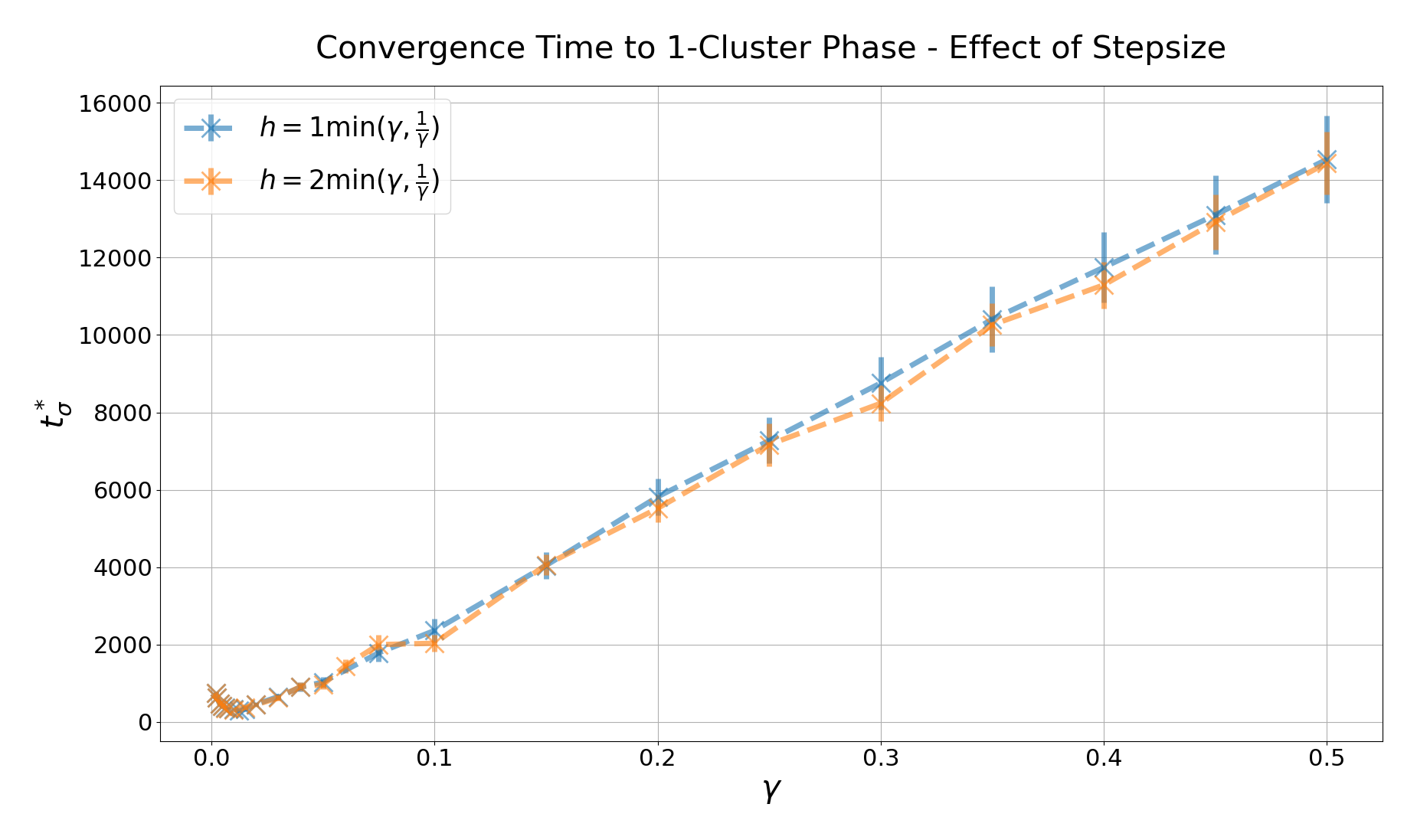}}
\caption{2D simulation results for convergence time $t_{\sigma}^*$ for two different stepsize settings. Simulation settings: Morse potential, $N=1000$, $\beta=150$. The number of trajectories that were averaged at each $\gamma$ corresponds to Tab. \ref{tab:avg_scheme_2D} for the curve at smaller stepsizes. The other curve uses either the same numbers or 50 to 100 more trajectories at each $\gamma$.}
\label{fig:effect_of_stepsize_convergence_times}
\end{figure}

We observe that all stepsizes lead to different dynamical paths (as is normal when simulating LD). The $\delta r$ plot shows that all trajectories undergo cluster formation and -merging until time $\sim$2,500 (refer to Sec. \ref{sup:sec: MSD_details} for a discussion of the mean squared displacement). From looking at $d_{\text{com}}$, we see that the $h=1$ trajectory even reaches the one-cluster stationary state, whereas the other two seem to be stuck in a metastable multi-cluster state. All trajectories, however, remain numerically stable and show the expected physical behavior, and from looking at the first two observables alone it is not clear which stepsize should be preferred. The third observable, the kinetic temperature $T_{\text{kin}}$, is special in the sense that it can measure integration accuracy. For Langevin dynamics (\ref{eq: Langevin Dynamics}), the kinetic temperature $T_{\text{kin}}=\frac{\|\boldsymbol{v} \|^2}{d}$, with dimensionality $d$ and Euclidean norm $\|.\|$, satisfies 
\begin{equation}\label{sup:eq:Tkin}
\mathbb{E}\big(T_{\text{kin}}\big) = T = \frac{1}{\beta},
\end{equation}
because $\boldsymbol{v}\sim \mathcal{N}(\boldsymbol{0},\beta^{-1}\mathbf{I})$ in equilibrium. When simulating LD with a discretization stepsize $h$, the deviation of $T_{\text{kin}}$ (rather, its moving time average) from $T$ is a measure of the $h$-dependent bias introduced by the discretization. Fig. \ref{fig:effect_of_stepsize_observables} shows that the bias decreases for decreasing stepsize as expected, and becomes negligible at $h=0.1$. However, as mentioned before, there is incentive to use a large stepsize to reduce compute work. In our experiments in the main text, at the particular $\gamma$ used in Fig. \ref{fig:effect_of_stepsize_observables}, we used the intermediate of the shown stepsizes, $h=0.5$. While there is a small visible bias left at this stepsize, the integrator we used for our main experiments, BAOAB (see Sec. \ref{sec:numerical_integrators}), tends to have excellent accuracy on relevant configurational observables even when there is a small bias left on the kinetic temperature (see \cite{ben_charlie_robust} or Fig. 7.7 in \cite{benMDbook}). Therefore, if all three shown stepsizes still lead to physically correct behavior, allowing for a small bias in $T_{\text{kin}}$ is a suitable approach. How large the bias may be strongly depends on the application and the particular quantities of interest. For example, from just looking at $T_{\text{kin}}$ in Fig. \ref{fig:effect_of_stepsize_observables} and the discussion until now, it is not clear why one can't go with the largest of the three stepsizes. Only through further experimentation can the choice against this be justified. Fig. \ref{fig:effect_of_stepsize_convergence_times} shows the results for the convergence time $t_{\sigma}^*$ for the Morse potential in 2D (see also Sec. \ref{sec: 2D results}) in comparison to the same results when the stepsize at each $\gamma$ is twice as large. Note that for $\gamma=0.5$ the two stepsizes correspond to the two larger ones used in Fig. \ref{fig:effect_of_stepsize_observables}.

We observe that the smaller stepsizes lead to more stable results. This is particularly noteworthy as the confidence intervals for the larger stepsizes are either comparable or narrower than the ones for the smaller stepsizes as they were obtained by averaging over a larger number of trajectories. This shows that the largest stepsize in Fig. \ref{fig:effect_of_stepsize_observables} would indeed have been too large for the main experiments (at the obtained statistical resolution).\\
Finally, as explained in Sec. \ref{sec:numerical_integrators} in the main text, we scale the stepsize with $\gamma$ as $h=h_1\min(\gamma, \frac{1}{\gamma})$. For an experiment such as the one in Fig. \ref{fig:effect_of_stepsize_convergence_times}, for a given range of $\gamma$ to examine, one thus needs to choose $h_1$ such that the smallest resulting stepsize still leads to reasonable compute time at the corresponding $\gamma$. Similarly, the largest resulting stepsize should still lead to stable trajectories with reasonable bias. When the `production run' of the experiment consists of running hundreds of trajectories per $\gamma$, doing some early test runs to determine the right value of $h_1$ may ultimately decrease the overall (compute) work.

\FloatBarrier

\subsection{Number of independent trajectories.}\label{sec: trajectory_average}
To generate the plot of the convergence times to the one-cluster equilibrium phase in Fig. \ref{fig:convergence_times_1D}, we picked the following averaging scheme. For both the Gaussian and the Morse potential, we averaged over 100 trajectories for $\gamma<0.1$, 200 trajectories for $0.1 \leq \gamma <1$, and 400 trajectories for $\gamma \geq 1$. For the GEM4 potential, we follow the same scheme for $\gamma<0.75$ and use 600 trajectories for $\gamma \geq 0.75$ as the corresponding measurements for $t_{\sigma}^*$ adopted the highest variance. 

\begin{figure}[!htbp]
\centerline{\includegraphics[width=1\textwidth]{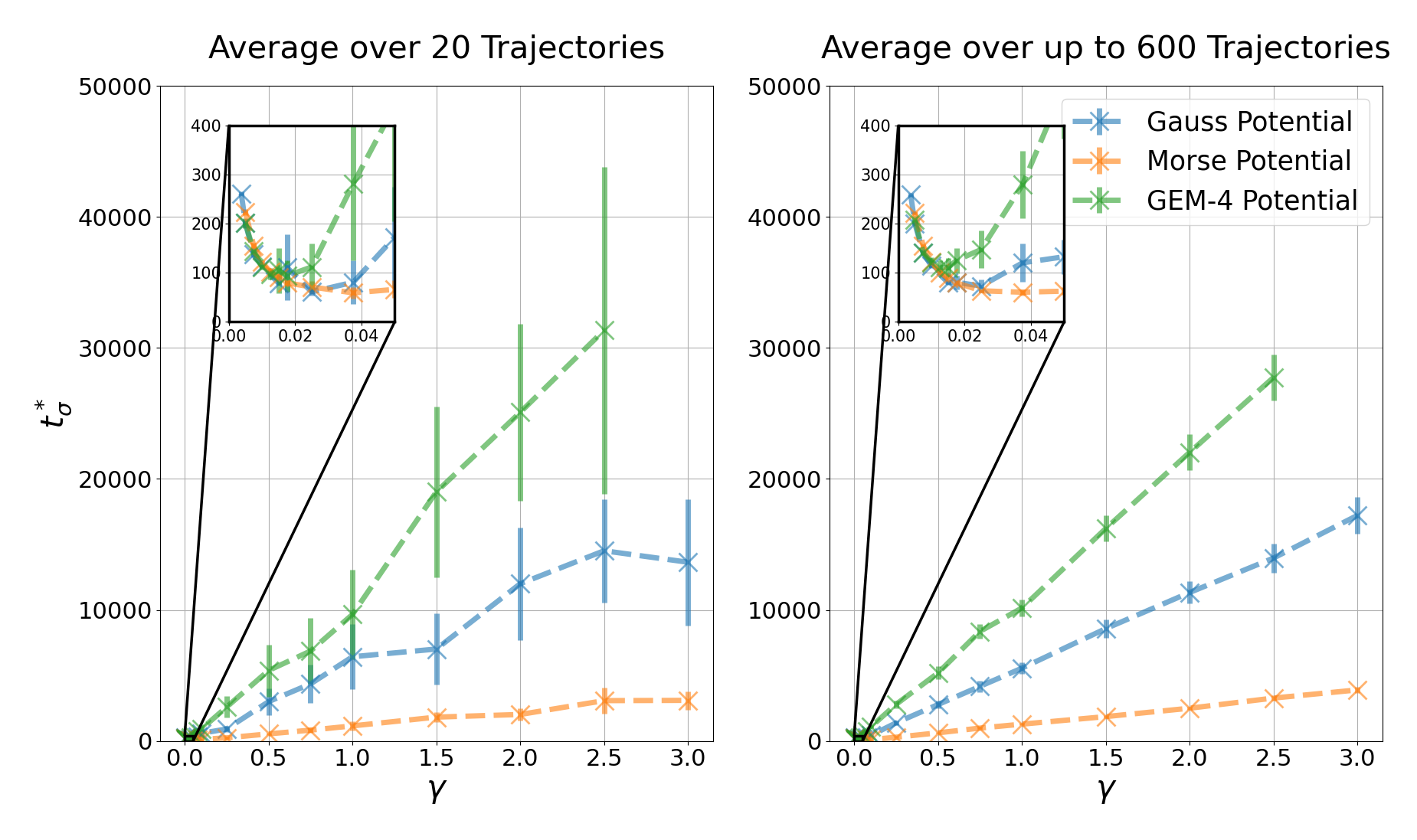}}
\caption{Simulation results for the friction-dependent times $t_{\sigma}^*$ to convergence to the stationary one-cluster configuration for $\beta>\beta_c$ when averaged over 20 independent trajectories per $\gamma$ vs. over up to 600 trajectories. Vertical bars denote 95\%-confidence intervals.}
\label{fig:influence of trajectory avg}
\end{figure}

Fig. \ref{fig:influence of trajectory avg} shows a comparison of the measured convergence times $t_{\sigma}^*$ when averaged over only 20 trajectories per friction value $\gamma$ compared to using the averaging scheme described above (up to 600 trajectories).  We observe that the linear-growth regimes are better resolved for larger numbers of independent samples, in particular for the Gaussian potential, and that the statistical variance is substantially reduced. The Morse potential already shows small confidence intervals for only 20 independent trajectories (in fact, they are even smaller than the intervals for the other potentials when using up to 600 trajectories). As mentioned in Sec. \ref{sec:results} in the main text, this is due to the overall wider interaction range of the Morse system, leading to faster merging of intermediate clusters.

\FloatBarrier

\subsection{More details on the Hamiltonian limit}\label{sup:sec:ham_limit}
\begin{figure}[!htbp]
\includegraphics[width=1.0\textwidth]{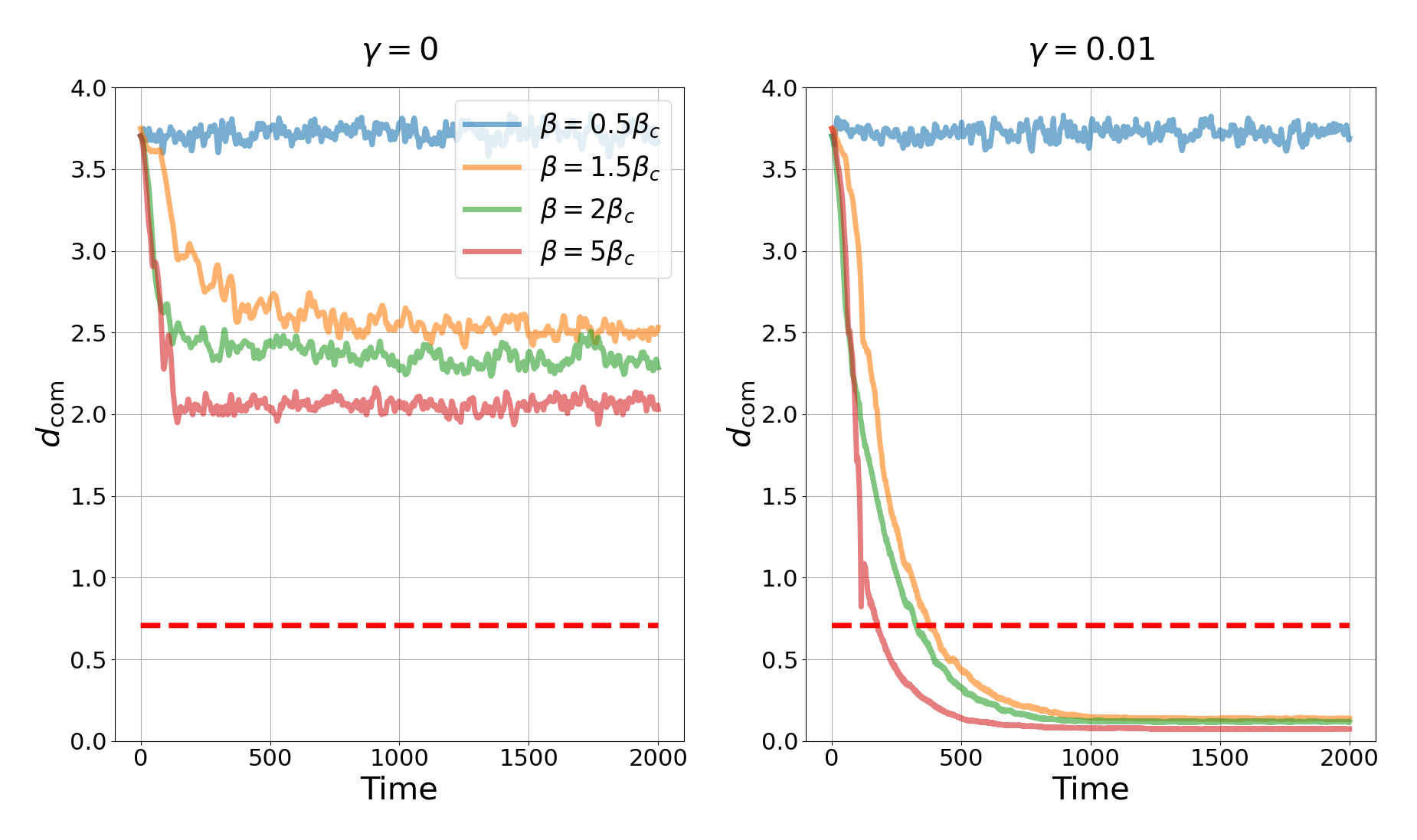}
\caption{\label{fig:cluster_size_vs_beta} $d_{\text{com}}$ at various temperatures for two different frictions $\gamma$. 2D IPS, Gaussian interaction potential. \textbf{Left:} Hamiltonian dynamics, $\gamma=0$. \textbf{Right:} Langevin dynamics at $\gamma=0.01$. Dashed line denotes $s$ from (\ref{eq: d_com_criterion}).}
\end{figure}

An equilibrium state of a cluster in an environment of unbound particles reminds of phase coexistence in statistical mechanics (e.g., during the liquid-gas transition).
Conceptually, what happens at $\gamma=0$ in terms of clustering is this: The attractive force naturally tries to pull the system towards a clustered state, in fact, even to a state with $x_1=x_2=...=x_N$. However, the total energy is conserved for $\gamma=0$, which means that only those states are admissible that satisfy $H(\boldsymbol{x},\boldsymbol{v})=H(\boldsymbol{x}_0, \boldsymbol{v}_0)$, the subscript 0 denoting initial conditions. States of minimal potential energy (fully-clustered states) would thus require the particles to have maximal kinetic energy. While these states are admissible, there are far more admissible points in phase space with less extreme energy distributions, and since all phase-space points in a microcanonical ensemble are equally likely, it follows that the system will most often be found in non-fully-clustered configurations. This is equivalent to saying that clustered states are energetically favourable but entropically unfavourable. From this intuition, it follows that if the total energy in the system is increased, clustered states will become less likely and/or the cluster width will increase, resembling less-tightly bound particles. We can easily verify this by running simulations at $\gamma=0$ for different values of $\beta$. While $\beta^{-1}$ does not describe the temperature of a thermostat in this setting, it does influence the total energy in the system because the initial momenta are still drawn from a Gaussian with variance $\beta^{-1}$.  We have already observed an increase in cluster width with decreasing $\beta$ for $\gamma>0$ in Fig. \ref{sup:fig:cluster_width} of Sec. \ref{sec: critical_temperature}. Here, we also examine it for $\gamma=0$, using the same 2D IPS as in Figs. \ref{fig:cluster_size_vs_gamma} and \ref{fig:final_cluster_state_snapshots}. The resulting $d_{\text{com}}$ curves are plotted in Fig. \ref{fig:cluster_size_vs_beta}, together with trajectories for $\gamma>0$ for comparison.
In both cases, $\gamma=0$ and $\gamma>0$, the uniform distribution remains stable at $\beta<\beta_c$.
For larger $\beta$, the $\gamma=0$ curves equilibrate to smaller values, supporting the reasoning above. The results for $\gamma>0$ are consistent with the behavior we observed before: the systems enter fully-clustered states, where the cluster widths also decrease for larger $\beta$. \\
It should be stressed, however, that the mechanism behind the decrease of cluster width for larger $\beta$ is different for the two cases. In both cases, the initial momenta are drawn from $\mathcal{N}(\boldsymbol{0},\beta^{-1}\mathbf{I})$, but this leads to different temperatures (as measured by the mean kinetic temperature, see Sec. \ref{sec: choice of stepsize}) in the canonical compared to the microcanonical ensemble. In the canonical ensemble, the temperature will remain stable at $\beta^{-1}$, whereas in the microcanonical ensemble it will be determined by the initial total energy in the system and not necessarily coincide with $\beta^{-1}$. This is equivalent to saying that canonical sampling preserves the momenta distribution $\mathcal{N}(\boldsymbol{0},\beta^{-1}\mathbf{I})$, whereas microcanonical sampling does not. So in the microcanonical case, $\beta$ enters only through the initial conditions, with smaller values injecting more energy into the system, hence wider clusters. In the canonical case, $\beta$ will govern the size of clusters independently of the initial conditions.

\FloatBarrier

\subsection{Interaction potentials}\label{sup:sec:potentials}
\begin{figure}[!htbp]
\centerline{\includegraphics[width=0.8\textwidth]{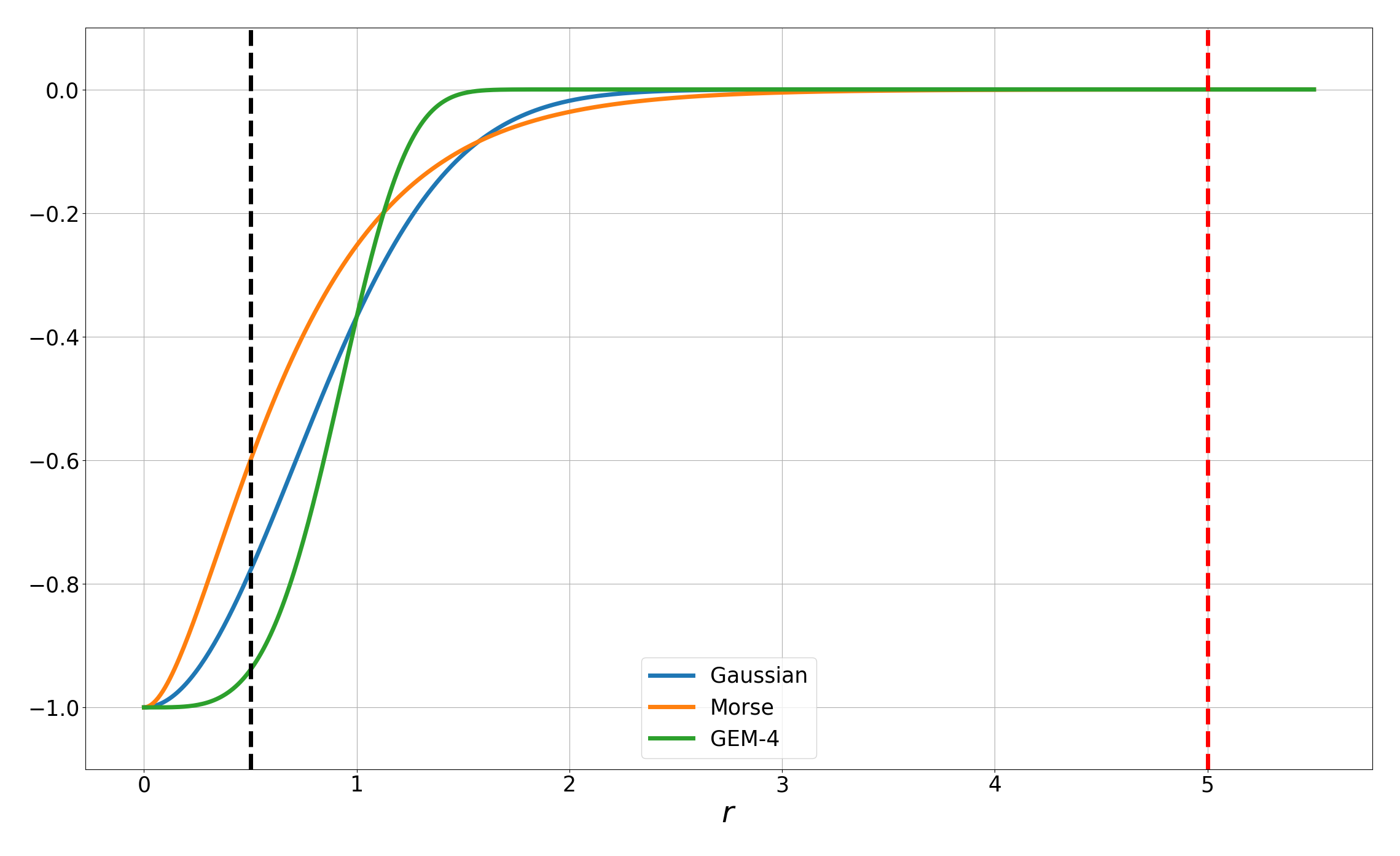}}
\vspace{-0.5cm}
\caption{The three pairwise interaction potentials used in this work. The dashed lines denote the points $r=0.5=\sigma^2$ (black) and  $r=5=\frac{L}{2}$ (red), see text.}
\label{fig:potentials}
\end{figure}

The three pairwise interaction potentials we consider in this work are given by functions $W(\boldsymbol{x})$, with $\boldsymbol{x} \in \mathbb{R}^d$, $d\in\{1,2,3\}$, defined as follows.
\begin{equation}\label{sup:eq:Potentials}
\begin{aligned} 
W(\boldsymbol{x}) &= - e^{-\frac{\|\boldsymbol{x} \|^2}{2\sigma^2}}  & \text{Gaussian Potential}, \\
W(\boldsymbol{x}) &= D_e \Big(e^{-2a\|\boldsymbol{x}\|} -2e^{-a\|\boldsymbol{x}\|}\Big)  & \text{Morse potential}, \\
W(\boldsymbol{x}) &= -e^{-\Big(\frac{\|\boldsymbol{x}\|}{\sqrt{2\sigma^2}}\Big)^\alpha}  & \text{GEM-$\alpha$},
\end{aligned}
\end{equation}
where in practice $\bm{x}$ will be given by the difference of two particle positions, e.g., $\bm{x}:=\bm{x}_{i,j}:=\bm{x}_i-\bm{x}_j$. The value of $W(\bm{x})$ just depends on the norm of the vector. We set $\sigma^2=0.5$, $a=2$, $D_e=1$, and $\alpha=4$. Note that the generalized exponential model (GEM-$\alpha$) reduces to the Gaussian potential for $\alpha=2$ and becomes supergaussian for $\alpha=4$. With the parameters chosen that way, the potentials have comparable width and equal depth, which is important to relate differences in the numerical results to the actual shape of the potentials, rather than vastly different range or depth. It also leads to them having similar critical temperatures (Tab. \ref{tab:beta_crit} in the main text) and friction ranges of interest. This reduces the need for hyperparameter tuning when moving from one potential to the next, which is particularly important for the cluster detection via DBSCAN, see the discussion in Sec. \ref{sup: sec: DBSCAN details}.\\
Fig. \ref{fig:potentials} shows the potentials $W(\boldsymbol{x})$ plotted against  $r:=\|\boldsymbol{x}\|. $
The black dashed line denotes the threshold distance that is used to detect the one-cluster state via criterion (\ref{eq: d_com_criterion}). For our main experiments, it is given by $s=\sigma^2=0.5$, i.e., a characteristic length of the potentials.
The red dashed line denotes the distance $\frac{L}{2}=5$, where $L=10$ is the edge length of the cubic simulation box which is the domain of the particle positions. Due to the usage of periodic boundary conditions and the minimum-image convention, the potential parameters need to be chosen such that the potentials (actually, their derivatives) are negligible for $r\geq\frac{L}{2}$, which is satisfied in our case as evident by the figure.\\

\FloatBarrier

\subsection{Temperature-dependent convergence and cluster width}\label{sup:sec: temp_dependent_cluster_widths}
\begin{figure}[!htbp]
\includegraphics[width=1.0\textwidth]{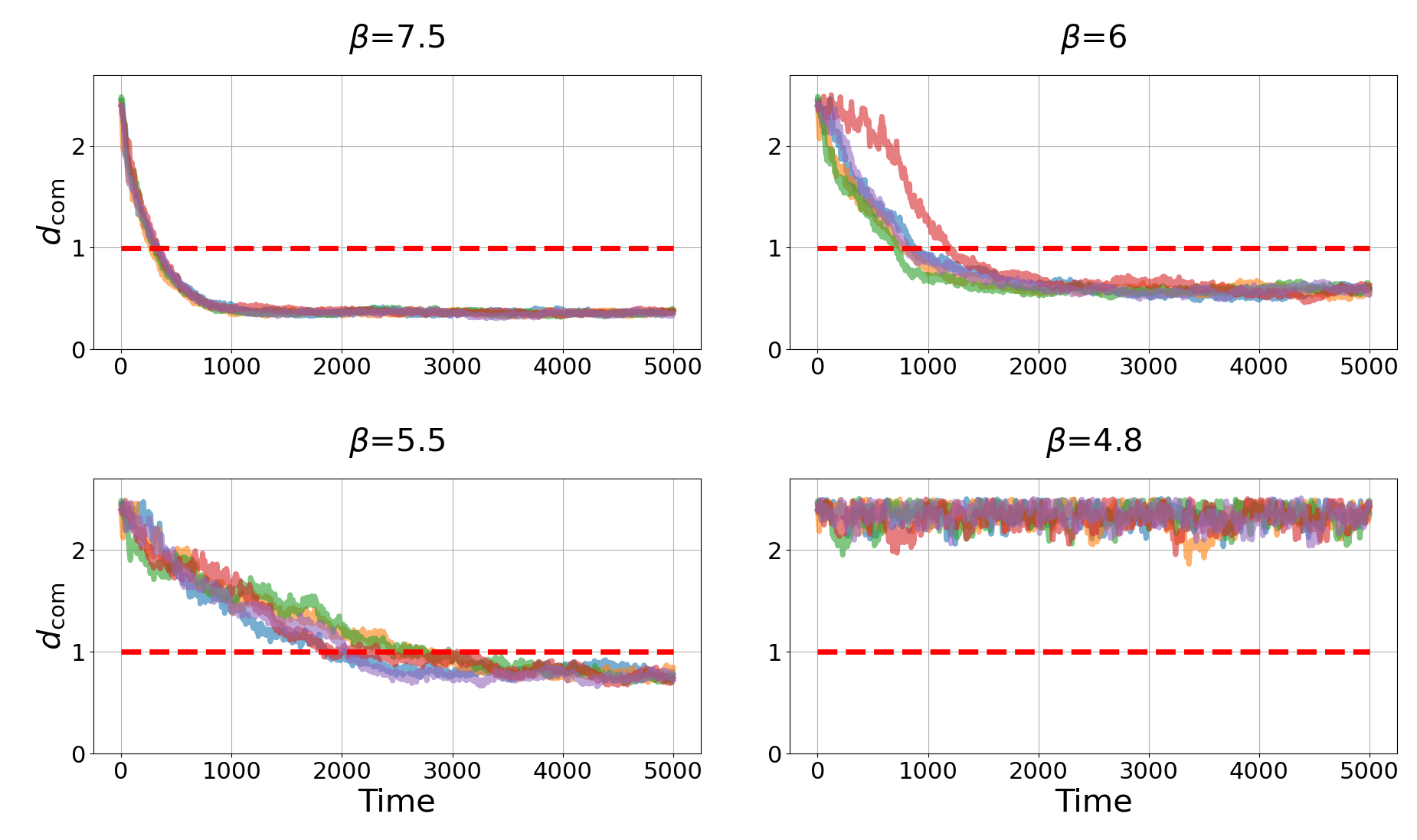}
\caption{\label{sup:fig:cluster_width} Mean distance to centre of mass $d_{\text{com}}$ for four different $\beta$-values for the same system as in Fig. \ref{fig:critical_temperature}, using $\sigma^2=0.7$. Each panel shows 5 independent trajectories. The red dashed line denotes the convergence threshold. The critical $\beta$ value is given by $\beta_c \approx 5$. }
\end{figure}

We examine the convergence to equilibrium for a one-dimensional IPS for various temperature parameters $\beta$, where we experimentally estimated $\beta_c \approx 5$. Fig. \ref{sup:fig:cluster_width} shows the results for the observable $d_{\text{com}}$. For all shown $\beta > \beta_c$, the times to equilibrium, given by the times at which $d_{\text{com}}$ falls below $s$ (red dashed lines), increase with increasing temperature $\beta^{-1}$. The same is true for the magnitude of $d_{\text{com}}$ in equilibrium, which is an estimate for the cluster widths. For $\beta < \beta_c$, $d_{\text{com}}$ remains at its initial height, denoting stability of the uniform phase.

\FloatBarrier

\end{appendix}

\printbibliography
\end{document}